\documentclass[10pt]{article}  \usepackage{amsmath}
\usepackage{amscd}
\usepackage{amssymb}
\usepackage{amsthm} 
\usepackage{graphicx}      
\newtheorem{thm}{Theorem}[section]
\newtheorem{prop}[thm]{Proposition}
\newtheorem{lem}[thm]{Lemma}
\newtheorem{cor}[thm]{Corollary}  \theoremstyle{definition}
\newtheorem{df}[thm]{Definition}   \theoremstyle{definition}
\newtheorem{ques}[thm]{Question}
\newtheorem{prob}[thm]{Problem}
\newtheorem{rem}[thm]{Remark}                \theoremstyle{plain}
 \theoremstyle{definition}
\newtheorem{ex}[thm]{Example}   \def\CC{\Bbb{C}}
\def\RR{\Bbb{R}}  
\def\ZZ{\Bbb{Z}}
\def\CCI{\hat{\CC}}        \def\NN{\Bbb{N}} 
\def\B1{{\rm\kern.32em\vrule    width.12em       height1.4ex
depth-.05ex\kern-.28em 1}}

\begin{document}
\title{Interaction cohomology
of\\ forward or backward self-similar systems\footnote{Date: June 26, 2009. Published in 
Adv. Math., 222 (2009) 729--781.  
2000 Mathematical Subject Classification: 37F05, 37F20.  
Keywords: Self-similar systems, iterated function systems, cohomology, complex dynamics, 
rational semigroups, random iteration, 
Julia set, fractal geometry.}} 
\author{Hiroki Sumi\
\\   Department of Mathematics, Graduate School of Science,\\ 
Osaka University\\   
1-1, Machikaneyama,\ Toyonaka,\  Osaka,\ 560-0043,\ 
Japan\\ E-mail:    sumi@math.sci.osaka-u.ac.jp\\ 
http://www.math.sci.osaka-u.ac.jp/$\sim $sumi/}
\date{}
\maketitle
\begin{abstract}


We investigate the dynamics of forward or backward self-similar 
systems (iterated function systems) and the topological structure of their invariant sets. 
 We define a new cohomology theory (interaction cohomology)
  for  forward or backward  self-similar systems. 
We show that under certain conditions, 
the space of connected components of the invariant set is 
isomorphic to the inverse limit of the spaces of 
connected components of the realizations of the nerves of 
finite coverings ${\cal U} $ of the invariant set,  
where each ${\cal U}$ consists of 
(backward) images of the invariant set under elements of  
finite word length. 
We give a criterion 
for the invariant set to be connected. 
Moreover, we give a sufficient condition for the first 
  cohomology group to 
  have infinite rank. As an application, 
  we obtain many results on the dynamics of 
semigroups of polynomials. 
 Moreover, we define postunbranched systems and we investigate the interaction cohomology groups 
of such systems. Many examples are given.      
\end{abstract}
\section{Introduction} 
The theory of iterated function systems has been widely and deeply investigated in fractal geometry (\cite{Ha,F,K, MU, Ka1, Ka2}). 
It deals with systems ${\frak L}=(L,(h_{1},\ldots ,h_{m}))$, where 
$L$ is a non-empty compact metric space and  
$h_{j}:L\rightarrow L$ is a continuous map for each $j=1,\ldots ,m$, 
such that 
$L=\bigcup _{j=1}^{m}h_{j}(L).$ In this paper,  such a system $(L,(h_{1},\ldots ,h_{m}))$ 
is called a forward 
self-similar system (Definition~\ref{d:fss}). 
For any two forward self-similar systems 
${\frak L}_{1}=(L_{1},(h_{1},\ldots, h_{m}))$ and 
${\frak L}_{2}=(L_{2},(g_{1},\ldots ,g_{n}))$, 
a pair $\Lambda =(\alpha ,\beta )$, 
where $\alpha : L_{1}\rightarrow L_{2}$ is a continuous map and $\beta :\{ 1,\ldots ,m\} \rightarrow 
\{ 1,\ldots ,n\} $ is a map, is called a morphism of ${\frak L}_{1}$ to 
${\frak L}_{2}$ 
if $g_{\beta (j)}\circ \alpha =\alpha \circ h_{j}$ on $L$  for each $j=1,\ldots ,m.$ 
If $\Lambda $ is a morphism of ${\frak L}_{1}$ to ${\frak L}_{2}$,  we write $\Lambda :{\frak L}_{1}\rightarrow {\frak L}_{2}.$
If $\Lambda _{1}=(\alpha _{1},\beta _{1}):{\frak L}_{1}\rightarrow {\frak L}_{2}$ 
and $\Lambda _{2}=(\alpha _{2},\beta _{2}):{\frak L}_{2}\rightarrow {\frak L}_{3}$ are 
such morphisms, then $\Lambda _{2}\circ \Lambda _{1}:=(\alpha _{2}\circ \alpha _{1}, \beta _{2}\circ \beta _{1})$ 
is a morphism of ${\frak L}_{1}$ to ${\frak L}_{3}.$ Moreover, for each system ${\frak L}=(L,(h_{1},\ldots ,h_{m}))$, 
the morphism $Id _{{\frak L}}=(Id _{L}, Id):{\frak L}\rightarrow {\frak L}$ is called the identity morphism.  
With these notations, we have a category. This is called the category of forward self-similar systems (see Definition~\ref{d:fsscat}).        
For any forward self-similar system 
${\frak L}=(L,(h_{1},\ldots ,h_{m}))$, the set $L$ is called the invariant set of the system and 
each $h_{j}$ is called a generator of the system. In many cases, 
the invariant set is quite complicated. For example, 
the Hausdorff dimension of the invariant set may not be an integer (\cite{F, MU}). 

 Another famous subject in fractal geometry is the study of Julia sets (where the dynamics are unstable) of rational maps 
 on the Riemann sphere $\CCI $. (For an introduction to complex dynamics, see \cite{B,Mi}.) The Julia set can be defined for a rational semigroup, i.e., a semigroup 
of rational maps on $\CCI $ (\cite{HM, GR}). 
For a rational semigroup $G$, we denote by $F(G)$ the largest open subset  
of $\CCI $ on which the family of analytic maps $G$ is equicontinuous with respect to the spherical distance.  
The set $F(G)$ is called the {\bf Fatou set} of $G$, and the complement $J(G):=\CCI \setminus F(G)$ 
is called the {\bf Julia set} of $G.$ 
In \cite{S0}, it was shown that 
for a rational semigroup $G$ which is generated by finitely many elements $\{ h_{1},\ldots ,h_{m}\} $, 
the Julia set $J(G)$ of $G$ satisfies the following backward self-similarity property  
$J(G)=\bigcup _{j=1}^{m}h_{j}^{-1}(J(G))$ (Lemma~\ref{openjfinlem}). 
(For additional results on rational semigroups, 
see \cite{St, SS, SU0, SU1, S3,S1, S01, S5, S4, S10, S11, SdpbpI, SdpbpII, SdpbpIII}. For a software to draw graphics of 
the Julia sets of rational semigroups, see \cite{Co}.)  We also remark that the study of rational semigroups is 
directly and deeply related to that of random complex dynamics. (For results on random complex dynamics, 
see \cite{FS, BBR, Br1, GQL,  S11,S10, SdpbpIII, Srcd}.)
Based on the above point of view, 
it is natural to introduce the following ``backward self-similar systems.'' In this paper, 
${\frak L}=(L,(h_{1},\ldots ,h_{m}))$ is called a backward self-similar system if $L$ is a compact subset of 
a metric space $X$, $h_{j}:X\rightarrow X $ is a continuous map for each $j=1,\ldots ,m$, 
$L=\bigcup _{j=1}^{m}h_{j}^{-1}(L)$, and for each $z\in L$ and each $j$, $h_{j}^{-1}(\{ z\} )\neq \emptyset $ 
(see Definition~\ref{d:bss}).  
The category of backward self-similar systems is defined 
in a similar way to that of forward self-similar systems (Definition~\ref{d:bsscat}). 
For a topological manifold $M$, we investigate how the coordinate neighborhoods  
overlap to obtain topological or geometric information about $M.$ 
On the other hand, for the invariant set $L$ of a forward (resp. backward) self-similar system 
${\frak L}=(L,(h_{1},\ldots ,h_{m}))$,  
we do not have such good coordinate neighborhoods that are homeomorphic 
to open balls in Euclidian space anymore. However, 
we have small ``copies'' (images) $h_{w_{1}}\cdots h_{w_{k}}(L)$ (resp. 
$h_{w_{1}}^{-1}\cdots h_{w_{k}}^{-1}(L) $) of $L$ under finite word elements $h_{w_{1}}\cdots h_{w_{k}}.$ 
These small copies contain important information on the topology of the invariant set $L.$ For example, 
we have the following well-known result: 
\begin{thm}[a weak form of (Theorem 4.6 in \cite{Ha}) or (Theorem 1.6.2 in \cite{K})]
\label{weakconnthm} 
\ \\ 
Let ${\frak L}=(L,(h_{1},\ldots ,h_{m}))$ be a forward self-similar system such that 
for each $j=1,\ldots ,m$, $h_{j}:L\rightarrow L$ is a contraction. Then, 
$L$ is connected if and only if   
for each $i,j\in \{1,\ldots ,m\} $, there exists a sequence 
$\{ i_{t}\} _{t=1}^{s}$ in $\{ 1,\ldots ,m\} $ such that 
$i_{1}=i,i_{s}=j$, and $h_{i_{t}}(L)\cap h_{i_{t+1}}(L)\neq \emptyset $ for each 
$t=1,\ldots ,s-1.$ 
\end{thm}
One motivation of this paper is to generalize and further develop the essence of Theorem~\ref{weakconnthm}. 
The following is a natural question:
\begin{ques}
\label{ques1}
For a fixed $k\in \NN $, we ask in what fashion do the small images $h_{w_{1}}\cdots h_{w_{k}}(L)$ (resp. 
$h_{w_{1}}^{-1}\cdots h_{w_{k}}^{-1}(L) $) of $L$ under $k$-words $h_{w_{1}}\cdots h_{w_{k}}$ 
overlap? How does this vary as $k$ tends to $\infty $? 
\end{ques} 
Here are some other natural questions:
\begin{ques}
\label{ques2}
What can we say about the topological aspects of the invariant set $L$? 
How many connected components does $L$ have? What about the number of connected components of 
the complement of $L$ when $L$ is embedded in a larger space?  
\end{ques} 
\begin{ques}
\label{ques3}
How can we describe the dynamical complexity of these (forward or backward) self-similar systems? 
How can we describe the interaction of different kinds of dynamics inside a single 
(forward or backward) self-similar system? 
How can we classify the isomorphism classes of forward or backward self-similar systems?  
How are these questions related to Question~\ref{ques1} and \ref{ques2}?
\end{ques}
These questions are profoundly related to the dynamical behavior of the systems ${\frak L}.$ 
In this paper, to investigate the above questions,  
we introduce a new kind of cohomology theory for such systems, which we call ``interaction cohomology.'' 
 We do this as follows. 
For each $m\in \NN $, let $\Sigma _{m}^{\ast }:= \bigcup _{n=1}^{\infty }\{ 1,\ldots ,m\} ^{n}$ (disjoint union). 
For each $w=(w_{1},\ldots ,w_{n})\in \{ 1,\ldots ,m\} ^{n}$, we set 
$|w|:=n$ and $\overline{w}:=(w_{n},\ldots ,w_{1}).$  
Let ${\frak L}=(L,(h_{1},\ldots ,h_{m}))$ be a forward (resp. backward) self-similar system. 
For each $w=(w_{1},\ldots ,w_{n})\in \Sigma _{m}^{\ast }$, we set 
$h_{w}:=h_{w_{n}}\circ \cdots \circ h_{w_{1}}.$ 
Let 
${\cal U}_{k}={\cal U}_{k}({\frak L})$ be the finite covering of $L$ 
defined as ${\cal U}_{k}:=\{ h_{\overline{w}}(L)\mid 
w\in \Sigma _{m}^{\ast }, |w|=k\} $ (resp. 
${\cal U}_{k}:=\{ h_{w}^{-1}(L)\mid w\in \Sigma _{m}^{\ast }, |w|=k\} $).  
Note that for each $k\in \NN $, ${\cal U}_{k+1}$ is a refinement of ${\cal U}_{k}.$  
Let $R$ be a $\ZZ $ module. Let $N_{k}=N_{k}({\frak L})$ be the nerve of 
${\cal U}_{k}$. Thus $N_{k}$ is a simplicial complex such that 
the vertex set is equal to $\{ w\in \Sigma _{m}^{\ast }\mid |w|=k\} $ and 
mutually distinct $r$ elements $w^{1},\ldots ,w^{r}\in \Sigma _{m}^{\ast }$ 
with $|w^{1}|=\cdots =|w^{r}|=k$ make an $(r-1)$-simplex of $N_{k}$ if and only if 
$\bigcap _{j=1}^{r}h_{\overline{w^{j}}}(L)\neq \emptyset $ (resp. $\bigcap _{j=1}^{r}h_{w^{j}}^{-1}(L)\neq \emptyset $).  
 Let $\varphi _{k}: N_{k+1}\rightarrow N_{k}$ be the simplicial map 
defined as $(w_{1},\ldots ,w_{k+1})\mapsto (w_{1},\ldots ,w_{k})$ for each 
$(w_{1},\ldots ,w_{k+1})\in \{ 1,\ldots ,m\} ^{k+1}.$ 
(For an example of $N_{k}$, see Example~\ref{ex:SG1}, Figure~\ref{fig:sgasket}, 
Figure~\ref{fig:sierpn1n2-2}.)  
 We consider the cohomology groups $H^{r}(N_{k};R).$ Note that 
 $\{ \varphi _{k}^{\ast }: H^{r}(N_{k};R)\rightarrow H^{r}(N_{k+1};R)\} _{k\in \NN }$ 
 makes a direct system of $\ZZ $ modules. 
The interaction cohomology 
groups $\check{H}^{r}({\frak L};R)$ are defined to be the direct limits $\varinjlim _{k}H^{r}(N_{k};R)$ 
(see Definition~\ref{bssicdf}, 
Definition~\ref{fssicdf}). 
Note that $\check{H}^{r}({\frak L};R)\cong \check{H}^{r}(\varprojlim _{k}|N_{k}|;R)$ (see \cite{W}), 
where for each simplicial complex $K$, we denote by $|K|$ the realization of $K $ (\cite[p.110]{Sp}).     
Note also that ${\frak L}\mapsto H^{\ast }(N_{k}({\frak L});R)$ and ${\frak L}\mapsto \check{H}^{\ast }({\frak L};R)$ 
are contravariant functors from the category of forward (resp. backward) self-similar systems 
to the category of $\ZZ $ modules (Remark~\ref{r:functor}). 
In particular, 
if ${\frak L}_{1}\cong {\frak L}_{2}$, then $H^{\ast }(N_{k}({\frak L}_{1});R)\cong H^{\ast }(N_{k}({\frak L}_{2});R)$ 
and $\check{H}^{\ast }({\frak L}_{1};R)\cong \check{H}^{\ast }({\frak L}_{2};R).$ 
Thus the isomorphism classes of $H^{\ast }(N_{k}({\frak L});R)$ and $\check{H}^{\ast }({\frak L};R)$ are 
invariant under the isomorphisms of forward (resp. backward) self-similar systems.  
We have a natural homomorphism $\Psi $ from the interaction cohomology groups of a system ${\frak L}$ 
to the \v{C}ech cohomology groups $\check{H}^{\ast }(L;R)$ of the invariant set $L$ of 
the system ${\frak L}$ (see Remark~\ref{nathomrem}). 
Note that by the Alexander duality theorem (\cite{Sp}), 
for a compact subset $K$ of an oriented $n$-dimensional 
manifold $X$, there exists an isomorphism 
$\check{H}^{p}(K;R)\cong H_{n-p}(X,X\setminus K;R)$ (hence if $X=\RR ^{n}$ 
then $\check{H}^{p}(K;R)\cong \tilde{H}_{n-p-1}(X\setminus K;R)$, where 
$\tilde{H}_{\ast }$ denotes the reduced homology).  
For a forward self-similar system ${\frak L}=(L,(h_{1},\ldots ,h_{m}))$ such that 
each $h_{j}:L\rightarrow L$ is a contraction, $\Psi $ is an isomorphism (see Remark~\ref{r:Psiiso}).  
However, $\Psi $ is not 
an isomorphism in general. In fact, $\Psi $ may not even be a monomorphism (see Proposition~\ref{rgeq2exprop}). 
In this paper, we show the following result:
\begin{thm}[see Theorems~\ref{precompojth} and \ref{ssprecompojth}]
\label{introprecompojth}
Let ${\frak L}=(L,(h_{1},\ldots ,h_{m}))$ be a forward (resp. backward) self-similar system. 
Suppose that for each $x\in \{ 1,\ldots ,m\} ^{\NN }$, $\bigcap _{j=1}^{\infty }h_{x_{1}}\cdots h_{x_{j}}(L)$ 
(resp. $\bigcap _{j=1}^{\infty }h_{x_{1}}^{-1}\cdots h_{x_{j}}^{-1}(L)$) is connected. Then, we have the following.
\begin{itemize}
\item[{\em (1)}]
There exists a bijection {\em Con}$(L)\cong \varprojlim_{k}\mbox{{\em Con}}(|N_{k}|)$, where for each topological 
space $X$, we denote by {\em Con}$(X)$ the set of all connected components of $X.$ 
\item[{\em (2)}]
$L$ is connected if and only if $|N_{1}|$ is connected, that is, 
for each $i,j\in \{1,\ldots ,m\} $, there exists a sequence 
$\{ i_{t}\} _{t=1}^{s}$ in $\{ 1,\ldots ,m\} $ such that 
$i_{1}=i,i_{s}=j$, and $h_{i_{t}}(L)\cap h_{i_{t+1}}(L)\neq \emptyset $ 
(resp. $h_{i_{t}}^{-1}(L)\cap h_{i_{t+1}}^{-1}(L)\neq \emptyset $) for each 
$t=1,\ldots ,s-1.$ 
\item[{\em (3)}] Let $R$ be a field. Then,  
 $\sharp \mbox{{\em Con}}(L)<\infty $  if and only if 
$\dim _{R}\check{H}^{0}({\frak L};R)<\infty .$ 
If  $\sharp \mbox{{\em Con}}(L)<\infty $, then  
$\Psi :\check{H}^{0}({\frak L};R)\rightarrow \check{H}^{0}(L;R)$ is an isomorphism. 
\end{itemize}
\end{thm}
Note that Theorem~\ref{introprecompojth} (2) generalizes Theorem~\ref{weakconnthm}. 
Moreover, note that until now, no research has investigated the space of connected components of the 
invariant set of such a system; Theorem~\ref{introprecompojth} 
gives us new insight into the topology of the invariant sets of such systems.

Furthermore, a sufficient condition for the rank of the first interaction cohomology 
groups to be infinite is given (Theorem~\ref{preicmainthm}, \ref{sspreicmainthm}). 
More precisely, we show the following result: 
\begin{thm}[Theorem~\ref{preicmainthm}]
\label{intropreicmainthm}
Let $\frak{L}=(L,(h_{1},\ldots ,h_{m}))$ be a backward self-similar
 system. 
 Let $R$ be a field. 
We assume all of the following conditions {\em (a),...,(d)}: 
\begin{itemize}
\item[{\em (a)}]$|N_{1}|$ is connected.   
\item[{\em (b)}]$(h_{1}^{2})^{-1}(L)\cap (\bigcup _{i: i\neq 1}h_{i}^{-1}(L))
=\emptyset .$  
\item[{\em (c)}]There exist mutually distinct elements 
$j_{1}, j_{2}, j_{3}\in \{ 1,\ldots ,m\} $ such that 
$j_{1}=1$ and  such that for each $k=1,2,3,$ 
$h_{j_{k}}^{-1}(L)\cap h_{j_{k+1}}^{-1}(L)
\neq \emptyset ,$ where $j_{4}:=j_{1}.$   
\item[{\em (d)}]For each $s,t\in \{ 1,\ldots ,m\} $, 
if $s,t,1$ are mutually distinct, 
then $h_{1}^{-1}(L)\cap h_{s}^{-1}(L)\cap h_{t}^{-1}(L)=\emptyset .$  
\end{itemize}
Then, $\dim _{R}\check{H}^{1}({\frak L};R)=
\infty .$
\end{thm}
A similar result is given for forward self-similar systems ${\frak L}$ (Theorem~\ref{sspreicmainthm}). 

Using Leray's theorem (\cite{G}), 
we also find a sufficient condition for the natural homomorphism $\Psi $ to be a 
monomorphism between the first cohomology groups (Lemma~\ref{icfundlem}). 

 The results in the above paragraphs are applied to the study of the dynamics of polynomial semigroups (i.e., semigroups 
of polynomial maps on $\CCI $).
For a polynomial semigroup $G$, 
we set 
$P(G):= \overline{\bigcup _{g\in G}\{ \mbox{all critical values of }g:\CCI \rightarrow \CCI \} }.$ 
We say that a polynomial semigroup $G$ is postcritically bounded if $P(G)\setminus \{ \infty \} $ is 
bounded in $\CC .$ 
For example, if $G$ is generated by 
a subset of $\{ h(z)=cz^{a}(1-z)^{b}\mid 
a,b\in \NN , c>0, c(\frac{a}{a+b})^{a}(\frac{b}{a+b})^{b}\leq 1\} $, 
then $G$ is postcritically bounded (see Remark~\ref{pgrem} or \cite{SdpbpI}). 
Regarding the dynamics of postcritically bounded polynomial semigroups, 
there are many new and interesting phenomena which cannot hold in the 
dynamics of a single polynomial (\cite{SdpbpI, SdpbpII, SdpbpIII, S10}). Combining Theorem~\ref{introprecompojth} (Theorem~\ref{precompojth})  
with potential theory, 
we show the following result:
\begin{thm}[Theorem~\ref{compojth}]
\label{introcompojth}
Let $m\in \NN $ and for each $j=1,\ldots ,m$, 
let $h_{j}:\CCI \rightarrow \CCI $ be a polynomial map with $\deg (h_{j})\geq 2.$ 
Let $G$ be the polynomial semigroup generated by $\{ h_{1},\ldots ,h_{m}\} .$ 
Suppose that $G$ is postcritically bounded.  
Then, for the backward self-similar system 
${\frak L}=(J(G),(h_{1},\ldots ,h_{m}))$, 
all of the statements {\em (1),(2)}, and {\em (3)} in Theorem~\ref{introprecompojth} hold.  
\end{thm}
Moreover, 
combining Theorem~\ref{intropreicmainthm} (Theorem~\ref{preicmainthm}), 
Theorem~\ref{introcompojth} (Theorem~\ref{compojth}), 
the Riemann-Hurwitz formula (\cite{B,Mi}), Leray's theorem (\cite{G}), 
and the Alexander duality theorem (\cite{Sp}), we give a sufficient condition for the Fatou set (where 
the dynamics are stable) of a postcritically bounded polynomial semigroup $G$ 
to have infinitely many connected components (Theorem~\ref{icmainthm}). 
More precisely, we show the following result:
\begin{thm}[Theorem~\ref{icmainthm}]
\label{introicmainthm}
Let $m\in \NN $ and for each $j=1,\ldots ,m$, 
let $h_{j}:\CCI \rightarrow \CCI $ be a polynomial map with $\deg (h_{j})\geq 2.$ 
Let $G$ be the polynomial semigroup generated by $\{ h_{1},\ldots ,h_{m}\} .$ 
Suppose that $G$ is postcritically bounded.  
Moreover, regarding the backward self-similar system ${\frak L}=(J(G),(h_{1},\ldots ,h_{m}))$, 
suppose that all of the conditions {\em (a), (b), (c),}  and {\em (d)} in the assumptions of Theorem~\ref{intropreicmainthm} 
hold. Let $R$ be a field. Then, we have that 
$\dim _{R}\check{H}^{1}({\frak L}; R)=
\dim _{R}\Psi (\check{H}^{1}({\frak L};R))=\infty $, 
$\Psi : \check{H}^{1}({\frak L}; R)
\rightarrow \check{H}^{1}(J(G);R)$ 
is a monomorphism, and the Fatou set $F(G)$ of $G$ has infinitely many connected components.
\end{thm}
Moreover, we give an example 
of a finitely generated postcritically bounded polynomial semigroup $G=\langle h_{1},\ldots ,h_{m}\rangle $ 
such that the backward self-similar system ${\frak L}=(J(G),(h_{1},\ldots ,h_{m}))$ 
satisfies the assumptions of Theorem~\ref{introicmainthm} and   
the rank of the first interaction cohomology group of ${\frak L}$ 
 is infinite (Proposition~\ref{pbpch1inftyexprop}, Figure~\ref{fig:h1infty}).  

Theorem~\ref{introprecompojth} and Theorem~\ref{introcompojth} have many applications.  
In fact, using the connectedness criterion for the Julia set of 
a postcritically bounded polynomial semigroup (Theorem~\ref{introcompojth}), 
we investigate the space of 
postcritically bounded polynomial semigroups having $2$ generators  
(\cite{S10}). 
As a result of this investigation, 
we can obtain numerous results on random complex dynamics. 
Indeed, letting $T_{\infty }(z)$ denote the probability of 
the orbit under a seed value $z\in \CCI $ tending to $\infty $ under the random walk 
generated by the application of randomly selected polynomials from the set 
$\{ h_{1}, h_{2}\} $, we can show that in some parameter space,  
 the function $T_{\infty }$  
is continuous on $\CCI $ and varies only on the Julia set $J(G)$ of the 
corresponding polynomial semigroup $G$ generated by $\{ h_{1}, h_{2}\} .$  
In this case, the Julia set $J(G)$ is a very thin fractal set. 
Moreover, we can show that in some parameter region $\Lambda $, the Julia set $J(G)$ has 
uncountably many connected components,  
and in the boundary $\partial \Lambda $, the Julia set $J(G)$ is connected. 
This implies that the function $T_{\infty }$ on $\CCI $ is a complex analog of the Cantor function or 
Lebesgue's singular function. (These results have been announced in \cite{S10,S11}. See also \cite{Srcd}.)   

When we investigate a random complex dynamical system, it is important to know 
the topology of the Julia set and the Fatou set of the associated semigroup. 
Indeed, setting $C(\CCI ):=\{ \varphi : \CCI \rightarrow \CC \mid \varphi \mbox{ is continuous} \}$, 
for a general random complex dynamical system ${\cal D}$, under certain conditions, 
any unitary eigenvector $\varphi \in C(\CCI )$ of the transition operator $M$ of ${\cal D}$ is locally constant 
on the Fatou set of the associated semigroup (see \cite{Srcd}). Thus Theorem~\ref{introicmainthm} provides us 
important information of unitary eigenvectors of $M.$ Moreover, by \cite{Srcd}, 
the space $V$ of all finite linear combinations of unitary eigenvectors of $M$ is finite-dimensional, and  
for any $\varphi \in C(\CCI )$, $\{ M^{n}(\varphi ) \} _{n}$ tends to the finite-dimensional subspace $V.$

 Another area of interest in forward or backward self-similar systems ${\frak L}=(L,(h_{1},\ldots h_{m}))$ is 
the structure of the cohomology groups $H^{r}(N_{k};R)$ of the 
nerve $N_{k}$ of ${\cal U}_{k}$ and the growth rate $g^{r}({\frak L})$ of the rank $a_{r,k}$ of 
$H^{r}(N_{k};R)$ as $k$ tends to $\infty $, where $R$ is a field. (See Definition~\ref{bssicdf}, Definition~\ref{fssicdf}, Definition~\ref{d:cohcomp}.)  
 The above invariants are deeply related to the 
dynamical complexity of ${\frak L}.$ In section~\ref{Postunbranched}, we introduce 
``postunbranched'' systems (see Definition~\ref{d:pu}), and we show the following result:
\begin{thm}[for the precise statement, see Theorem~\ref{puthm}]
\label{introputhm}
Let ${\frak L}=(L,(h_{1},\ldots ,h_{m}))$ be a forward or backward self-similar system. 
Suppose that ${\frak L}$ is postunbranched. When ${\frak L}$ is a forward self-similar system, 
we assume further that $h_{j}:L\rightarrow L$ is injective for each $j=1,\ldots ,m.$ 
Let $R$ be a field. 
Then, we have the following. 
\begin{itemize}
\item[{\em (1)}] For each $r\geq 2$,  
there exists an exact sequence of $R$ modules:
$$
0\longrightarrow H^{r}(N_{1};R)\longrightarrow \check{H}^{r}({\frak L};R)
\longrightarrow \bigoplus _{j=1}^{m}\check{H}^{r}({\frak L};R)\longrightarrow 0.
$$
\item[{\em (2)}] If $r\geq 2$, or if $r=1$ and $|N_{1}|$ is connected, 
then $a_{r,k+1}=ma_{r,k}+a_{r,1}$ for each $k\in \NN .$ 
\item[{\em (3)}] $a_{0,k+1}=ma_{0,k}-m+a_{0,1}-a_{1,1}-ma_{1,k}+a_{1,k+1}$ for each $k\in \NN .$  
\item[{\em (4)}] $ma_{1,k}\leq a_{1,k+1}\leq ma_{1,k}+a_{1,1}$ and 
$ma_{0,k}-m+a_{0,1}-a_{1,1}\leq a_{0,k+1}\leq ma_{0,k}-m+a_{0,1}$ for each $k\in \NN .$    
\item[{\em (5)}] {\em (a)} If $r\geq 1$, then $g^{r}({\frak L})\in \{ -\infty ,\log m\} .$ 
{\em (b)} $g^{0}({\frak L})\in \{ 0,\log m\} .$
\item[{\em (6)}] Let $r\geq 1$. Then,
 $\dim _{R}\check{H}^{r}({\frak L};R)$ is either $0$ or $\infty .$
\item[{\em (7)}] Suppose $m\geq 2.$  Then, \\ 
$\dim _{R}\check{H}^{0}({\frak L};R)\in  \{ x\in \NN \mid a_{0,1}\leq x\leq \frac{1}{m-1}
(m-a_{0,1}+a_{1,1})\} \cup \{ \infty \} .$ 
\end{itemize} 
\end{thm}
Moreover, for any  $n\in \NN \cup \{ 0\} $, we give an example of a postunbranched backward  
self-similar system ${\frak L}=(L,(h_{1},\ldots ,h_{n+2}))$ such that $L\subset \CC $ and 
the rank of the $n$-th interaction cohomology group $\check{H}^{n}({\frak L};R)$ 
of ${\frak L}$ is equal to $\infty $ (Proposition~\ref{rgeq2exprop}).  
In this case, if $n\geq 2$, 
the natural homomorphism $\Psi : \check{H}^{n}({\frak L};R)\rightarrow \check{H}^{n}(L;R)$ is not a monomorphism, 
since for each $l\geq 2$ the \v{C}ech cohomology group $\check{H}^{l}(L;R)$ of $L$ is equal to zero. 
For any $n\in \NN \cup \{ 0\} $, 
we also give an example of a postunbranched 
forward self-similar system ${\frak L}=(L,(h_{1},\ldots ,h_{n+2}))$ such that $L\subset \RR ^{3}$, 
each $h_{j}:L\rightarrow L$ is injective, and the rank of the $n$-th interaction cohomology group of ${\frak L}$ is 
equal to $\infty$ (Proposition~\ref{rgeq2exprop}). In this case, if $n\geq 3$, 
the natural homomorphism $\Psi : \check{H}^{n}({\frak L};R)\rightarrow \check{H}^{n}(L;R)$ is not a monomorphism, 
since for each $l\geq 2$ the \v{C}ech cohomology group $\check{H}^{l}(L;R)$ of $L$ is equal to zero. 
We remark that these examples imply that 
the interaction cohomology groups of ${\frak L}=(L,(h_{1},\ldots ,h_{m}))$ 
may contain more (dynamical) information than the \v{C}ech cohomology 
groups of the invariant sets $L.$ 
Thus interaction cohomology groups of self-similar systems tell us 
information of dynamical behavior of the systems as well as 
the topological information of the invariant sets of the systems. 

Furthermore, we give many ways to construct examples of postunbranched 
systems (Lemmas~\ref{pssubsyslem}, \ref{puiteratelem}, \ref{psiterateblem}, \ref{psiterateflem}). From these, we see that 
if $L$ is one of the Sierpi\'{n}ski gasket, the snowflake, the pentakun, the heptakun, the octakun, and so on (\cite{K}), 
then there exists a postunbranched forward self-similar system ${\frak L}=(L,(h_{1},\ldots ,h_{m}))$ such that 
each $h_{j}:L\rightarrow L$ is an injective 
contraction (Examples~\ref{SGex}, \ref{PSex}). Moreover, we also see that  
for each $n\in \NN $, any subsystem of an $n$-th iterate of the above ${\frak L}$ is a 
postunbranched forward self-similar system (Examples~\ref{SGex}, \ref{PSex}).  

 We summarize the purpose and the virtue to introduce interaction cohomology groups 
for the study of self-similar systems ${\frak L}=(L,\{ h_{1},\ldots ,h_{m}\} )$ as follows. 
\begin{itemize}
\item[(1)] We can get information about the dynamical behavior 
of the system ${\frak L}$ and the interaction of different maps in the system.   
The cohomology groups $\check{H}^{r}({\frak L};R)$, 
the cohomology groups $H^{r}(N_{k};R)$ of the nerve of ${\cal U}_{k}$, and 
the growth rate $g^{r}({\frak L})$ of the rank $a_{r,k}$ of $H^{r}(N_{k};R)$ are  
new invariants for the dynamics of self-similar systems. 
These invariants reflect the dynamical behavior and the complexity of the systems. 
Under certain conditions,  
we can show, by using these invariants, that two self-similar systems are not isomorphic, 
even when we cannot show this by using \v{C}ech cohomology groups of 
the invariant sets of the systems (e.g., Examples~\ref{SGcalex}, 
\ref{Pencalex},  \ref{Sfcalex}, \ref{SGsub7ex},
Proposition~\ref{rgeq2exprop}, and the proof of Proposition~\ref{rgeq2exprop}).      
Moreover, the interaction cohomology groups $\check{H}^{r}({\frak L};R)$ 
are new invariants for the dynamics of finitely generated semigroups of continuous maps. 
 (See Remark~\ref{r:semigrinv}.) 
\item[(2)] By using the natural homomorphism 
$\Psi : \check{H}^{r}({\frak L};R)\rightarrow \check{H}^{r}(L;R)$, 
we can get information about the \v{C}ech cohomology groups 
$\check{H}^{r}(L;R)$ of the invariant sets $L$ (Lemma~\ref{icfundlem}, 
Theorems~\ref{precompojth},  \ref{ssprecompojth}, 
\ref{preicmainthm}, \ref{sspreicmainthm}, \ref{compojth}, 
\ref{icmainthm}, Proposition~\ref{pbpch1inftyexprop}, Theorem~\ref{puthm}, 
Proposition~\ref{notinjmuexprop}, Examples~\ref{SGcalex}, \ref{Pencalex}, 
\ref{Sfcalex}, \ref{SGsub7ex}, \ref{SGsubcalex}).  
By using the Alexander duality theorem, the \v{C}ech cohomology groups of $L$ 
tell us information of the (reduced) homology groups of the complement of $L$ 
in the bigger space (see Theorem~\ref{icmainthm}, Proposition~\ref{pbpch1inftyexprop}).  
Under certain conditions, $\Psi $ is a monomorphism (Lemma~\ref{icfundlem}).  
If ${\frak L}$ is a forward self-similar system and if each $h_{j}: L\rightarrow L$ 
is contractive, then $\Psi : \check{H}^{r}({\frak L};R)\rightarrow \check{H}^{r}(L;R)$ 
is an isomorphism (Remarks~\ref{nathomrem}, \ref{r:Psiiso}). Moreover, under the same condition, 
for each $w\in \Sigma _{m}$ and $x_{0}\in L$ such that 
$x_{0}\in h_{w_{1}}\cdots h_{w_{k}}(L)$ for each $k$, 
the interaction homotopy groups $\check{\pi }_{r}({\frak L},w)$ (see Definition~\ref{bssicdf}) of ${\frak L}$
 are isomorphic to 
the shape groups $\check{\pi }_{r}(L,x_{0})$ of the invariant set $L$ (see Remark~\ref{r:Psiiso}). 
 (For the definition of 
shape groups and shape theory, see \cite{MS}.) 
\item[(3)]  
Interaction cohomology groups $\check{H}^{r}({\frak L};R)$ and $H^{r}(N_{k};R)$ may have 
more dynamical information of the systems than 
the \v{C}ech cohomology groups or shape groups of the invariant sets.    
The natural homomorphism $\Psi :\check{H}^{r}({\frak L};R)\rightarrow \check{H}^{r}(L;R)$ is 
not an isomorphism in general. Similarly, 
$\check{\pi }_{r}({\frak L},w)$ and $\check{\pi }_{r}(L,x_{0})$ 
($x_{0}\in h_{w_{1}}\cdots h_{w_{k}}(L)$ for each $k\in \NN $)  
are not isomorphic in general. It may happen that 
$\check{\pi }_{1}({\frak L},w)$ is not trivial, $\check{H}_{1}({\frak L};\ZZ )\neq 0$ and  
$\check{H}^{1}({\frak L}; \ZZ )\neq 0$, even though 
$\check{\pi }_{r}(L,x_{0}), \check{H}_{r}(L;\ZZ )$, and $\check{H}^{r}(L;\ZZ )$ are trivial for all $r\geq 1$ and 
for all $(w,x_{0})$ such that $x_{0}\in h_{w_{1}}\cdots h_{w_{k}}(L)$ for each $k\in \NN $ 
(see  Example~\ref{ex:finsetic}).   
Moreover, for each $n\geq 4$, there are many examples ${\frak L}=(L,(h_{1},\ldots ,h_{m}))$ 
(of postunbranched systems)  
such that $L\subset \RR ^{3}$ and 
$\check{H}^{r}(L;R)=0$ for each $r\geq 3$, but 
the interaction cohomology group $\check{H}^{n}({\frak L};R)$ 
is not zero (Proposition \ref{rgeq2exprop}). In these examples, 
since $\check{H}^{\ast }({\frak L};R)\cong \check{H}^{\ast }(\varprojlim _{k}|N_{k}({\frak L})|;R)$, 
the above statement means that the dimension of  $\varprojlim _{k}|N_{k}({\frak L})|$ is larger than 
that of $L.$ In other words, the manner of overlapping of small images of $L$ is ``more higher-dimensional'' 
than the invariant set $L.$ These phenomena reflect the complexity of the dynamics of the self-similar systems. 
We remark that as illustrated in Remark~\ref{r:jisbss}, Remark~\ref{r:minsetimp}, and examples in section~\ref{Preliminaries},  
it is important to investigate self-similar systems whose generators may not be contractive.   
\item[(4)]
For any two self-similar systems ${\frak L}_{1}=(L_{1},(h_{1},\ldots ,h_{m}))$ 
and ${\frak L}_{2}=(L_{2},(g_{1},\ldots ,g_{n}))$, 
interaction cohomology groups $\check{H}^{r}({\frak L}_{1};R)$ and 
$\check{H}^{r}({\frak L}_{2};R)$ may not be isomorphic even when  
$L_{1}$ and $L_{2}$ are homeomorphic 
(see Example~\ref{ex:finsetic}, Example~\ref{ex:fintriv}, 
and Remark~\ref{r:icdlh}).  
\item[(5)] 
For any forward self-similar system ${\frak L}=(L,(h_{1},\ldots ,h_{m}))$ whose generators $h_{j}$ 
are contractions, 
the interaction cohomology groups $\check{H}^{r}({\frak L}; R)$ are isomorphic 
to the \v{C}ech cohomology groups $\check{H}^{r}(L; R)$ (Remark~\ref{r:Psiiso}). Thus in this case 
the interaction cohomology groups $\check{H}^{r}({\frak L}; R)$ are not new invariants.  
However, even for such systems, $H^{r}(N_{k};R)$ and $g^{r}({\frak L})$ are new invariants. 
Given two forward self-similar systems ${\frak L}_{1}$ and ${\frak L}_{2}$ whose generators are  
contractions,  
 $H^{r}(N_{k}({\frak L}_{1});R)$ and $H^{r}(N_{k}({\frak L}_{2}); R)$ may not be isomorphic,  
and $g^{r}({\frak L}_{1})$ and $g^{r}({\frak L}_{2})$ may be different, and these invariants 
$H^{r}(N_{k};R)$ and $g^{r}$  may tell us that ${\frak L}_{1}$ and ${\frak L}_{2}$ are not isomorphic, 
even when the \v{C}ech cohomology groups of their invariants are isomorphic (Examples~\ref{SGcalex}, 
\ref{Pencalex}, 
\ref{Sfcalex}, \ref{SGsub7ex}).  
\item[(6)]
Under certain good conditions (e.g., postunbranched condition (Definition~\ref{d:pu})),   
$\check{H}^{r}({\frak L};R)$, $H^{r}(N_{k};R)$ and $g^{r}({\frak L})$ can be exactly calculated 
for each $r\geq 0$ (Theorem~\ref{puthm}). In particular, 
there are inductive formulae for  $H^{r}(N_{k};R)$ with respect to 
$k$ (with an exception when $r=1$ and $|N_{1}|$ is disconnected, in which  
the situation is more complicated, see Theorem~\ref{puthm}-\ref{puthmakvalues1}, Proposition~\ref{notinjmuexprop} and Remark~\ref{r:nonly1}).  
These results are applicable to many famous examples (e.g., the snowflake, the Sierpi\'{n}ski gasket, 
the pentakun, the hexakun, the heptakun, the octakun, etc., and any subsystems of their iteration systems,  
see Examples~\ref{SGex}, 
\ref{PSex}, \ref{SGcalex}, \ref{Pencalex}, \ref{Sfcalex}, \ref{SGsub7ex}) and 
many new examples (see Propositions~\ref{rgeq2exprop}, \ref{notinjmuexprop}).    
For one of the keys in the proof of the above results, see the diagrams (\ref{eq:coseqp1}), (\ref{eq:coseqp2}), 
which come from the long exact sequences of cohomology groups.  
\item[(7)] We can apply the results on the interaction cohomology 
to the self-similar systems whose generators are contractions, to the dynamics of 
polynomial (rational) maps on the Riemann sphere, and to the random complex dynamics (section~\ref{Application}).  
There are many important examples of forward or backward self-similar systems (section~\ref{Preliminaries}).  
When we investigate the random complex dynamics, 
we have to see the minimal sets (which are forward invariant) and the Julia sets (which are backward invariant) 
of the associated semigroups (Lemma~\ref{minimalsetlem}, Remark~\ref{r:minsetimp}). 
We have many phenomena which can hold in rational semigroups and random complex dynamics, 
but cannot hold in the usual iteration dynamics of a single polynomial (rational) map (see \cite{SdpbpI, SdpbpIII, Srcd}).  
By using interaction cohomology, these interesting new phenomena can be systematically investigated (see \cite{S10}). 
\item[(8)]  
Given self-similar systems or iterated function systems, we have been investigating the (fractal) dimensions 
of the invariant sets by 
using analysis or ergodic theory so far. 
However, if the small copies overlap heavily, it is very difficult to analyze the fractal dimensions, 
and there has been no invariant to study or classify the self-similar systems which could be 
understood well. 
Interaction cohomology of the systems 
is a new interest, rather than fractal dimensions of the invariant sets, 
and can be a new strong research interest in the field of self-similar systems, 
iterated function systems,  
the dynamics of semigroups of holomorphic maps, random complex dynamics, and more general semigroup actions.  
Overlapping of the small copies of the invariant sets is the most difficult point in the 
study of self-similar systems. 
Nevertheless, by using the interaction cohomology, 
we positively study the overlapping of the small copies, 
rather than avoiding the difficulty. 
Interaction cohomology provides a new language to investigate self-similar systems.
For results when we have overlapping of small copies, 
see Theorems~\ref{precompojth}, \ref{ssprecompojth}, 
\ref{preicmainthm}, \ref{sspreicmainthm}, \ref{compojth}, 
\ref{icmainthm}, and Proposition~\ref{pbpch1inftyexprop}.   
     
\end{itemize} 

We remark that it is a new idea to use homological theory 
when we investigate self-similar systems (iterated function systems) and 
their invariant sets (fractal sets). Using homological theory, we can 
introduce many new topological invariants of self-similar systems. 
Those invariants are naturally and deeply related to the dynamical behavior of 
the systems and the topological properties of the invariant sets of the systems. 
Thus, developing the theory of ``interaction (co)homology,'' 
we can systematically investigate the dynamics of self-similar systems. 
The results are applicable to fractal geometry, the dynamics of rational semigroups, and 
random complex dynamics.

 In section~\ref{Preliminaries}, we give some basic notations and definitions on forward or backward 
 self-similar systems. In section~\ref{Main}, we present the main results of this paper. 
 We provide some fundamental tools to prove the main results in section~\ref{Tools} and  
present the proofs of the main results in section~\ref{Proofs}.   

\  
 
\noindent {\bf Acknowledgement:} 
The author thanks Rich Stankewitz for many valuable comments. 
The author thanks Eri Honda-Sumi for drawing Figure~\ref{fig:sierpn1n2-2}.     
\section{Preliminaries}
\label{Preliminaries}
In this section, we give some fundamental notations and definitions on forward or backward self-similar systems. 
We sometimes use the notation from \cite{Sp}.  
\begin{df}
If a semigroup $G$ is generated by a family 
$\{ h_{1},\ldots ,h_{m}\} $ of elements of $G$, 
then we write $G=\langle h_{1},\ldots ,h_{m}\rangle .$
\end{df} 
\begin{df}
\label{d:fss}
Let $(L,d)$ be a non-empty compact metric space. 
Let $h_{j}:L\rightarrow L\ (j=1,\ldots ,m)$ be a continuous map. 
We say that ${\frak L}=(L,(h_{1},\ldots ,h_{m}))$ is a 
{\bf forward self-similar system} if $L=\bigcup _{j=1}^{m}h_{j}(L).$  
The set $L$ is called the invariant set of ${\frak L}.$ Each $h_{j}$ is called a 
generator of ${\frak L}.$ 
\end{df}
\begin{df}
\label{d:fsscat}
Let 
${\frak L}_{1}=(L_{1},(h_{1},\ldots, h_{m}))$ and 
${\frak L}_{2}=(L_{2},(g_{1},\ldots ,g_{n}))$ be any two 
forward self-similar systems.  
A pair $\Lambda =(\alpha ,\beta )$, 
where $\alpha : L_{1}\rightarrow L_{2}$ is a continuous map and $\beta :\{ 1,\ldots ,m\} \rightarrow 
\{ 1,\ldots ,n\} $ is a map, is called a morphism of ${\frak L}_{1}$ to 
${\frak L}_{2}$ if $g_{\beta (j)}\circ \alpha =\alpha \circ h_{j}$ on $L$  for each $j=1,\ldots ,m.$ 
If $\Lambda $ is a morphism of ${\frak L}_{1}$ to ${\frak L}_{2}$,  we write $\Lambda :{\frak L}_{1}\rightarrow {\frak L}_{2}$. 
If $\Lambda _{1}=(\alpha _{1},\beta _{1}):{\frak L}_{1}\rightarrow {\frak L}_{2}$ 
and $\Lambda _{2}=(\alpha _{2},\beta _{2}):{\frak L}_{2}\rightarrow {\frak L}_{3}$ are 
two morphisms, then $\Lambda _{2}\circ \Lambda _{1}:=(\alpha _{2}\circ \alpha _{1}, \beta _{2}\circ \beta _{1})$ 
is a morphism of ${\frak L}_{1}$ to ${\frak L}_{3}.$ Moreover, for each system ${\frak L}=(L,(h_{1},\ldots ,h_{m}))$, 
the morphism $Id _{{\frak L}}=(Id _{L}, Id):{\frak L}\rightarrow {\frak L}$ is called the identity morphism.  
With these notations, we have a category. This is called the category of forward self-similar systems.  
\end{df}
\begin{df}
\label{d:bss} 
 Let $X$ be a metric space.
 Let $h_{j}:X\rightarrow X\ (j=1,\ldots ,m)$ be a 
 continuous map. 
 Let $L$ be a non-empty compact subset of $X.$ 
 We say that ${\frak L}= (L,(h_{1},\ldots ,h_{m}))$ is a 
  {\bf backward self-similar system} 
if  
(1) $ L=\bigcup _{j=1}^{m}h_{j}^{-1}(L)$ 
 and (2) for each 
 $z\in L$ and each $j\in \{ 1,\ldots ,m\} $, 
 $h_{j}^{-1}(\{ z\} )\neq \emptyset .$
The set $L$ is called the invariant set of ${\frak L}.$  
Each $h_{j}$ is called a generator of ${\frak L}.$ 
\end{df}
\begin{df}
\label{d:bsscat}
Let 
${\frak L}_{1}=(L_{1},(h_{1},\ldots, h_{m}))$ and 
${\frak L}_{2}=(L_{2},(g_{1},\ldots ,g_{n}))$ be any two 
backward self-similar systems.  
A pair $\Lambda =(\alpha ,\beta )$, 
where $\alpha : L_{1}\rightarrow L_{2}$ is a continuous map and $\beta :\{ 1,\ldots ,m\} \rightarrow 
\{ 1,\ldots ,n\} $ is a map, is called a morphism of ${\frak L}_{1}$ to 
${\frak L}_{2}$ if $\alpha (h_{j}^{-1}(L_{1}))\subset g_{\beta (j)}^{-1}(L_{2})$ and 
$g_{\beta (j)}\circ \alpha =\alpha \circ h_{j}$ on $h_{j}^{-1}(L_{1})$  for each $j=1,\ldots ,m.$ 
If $\Lambda $ is a morphism of ${\frak L}_{1}$ to ${\frak L}_{2}$,  we write $\Lambda :{\frak L}_{1}\rightarrow {\frak L}_{2}$. 
If $\Lambda _{1}=(\alpha _{1},\beta _{1}):{\frak L}_{1}\rightarrow {\frak L}_{2}$ 
and $\Lambda _{2}=(\alpha _{2},\beta _{2}):{\frak L}_{2}\rightarrow {\frak L}_{3}$ are 
two morphisms, then $\Lambda _{2}\circ \Lambda _{1}:=(\alpha _{2}\circ \alpha _{1}, \beta _{2}\circ \beta _{1})$ 
is a morphism of ${\frak L}_{1}$ to ${\frak L}_{3}.$ Moreover, for each system ${\frak L}=(L,(h_{1},\ldots ,h_{m}))$, 
the morphism $Id _{{\frak L}}=(Id _{L}, Id):{\frak L}\rightarrow {\frak L}$ is called the identity morphism.  
With these notations, we have a category. This is called the category of backward self-similar systems.  
\end{df}
\begin{rem}
\label{r:ssiso}
Let ${\frak L}_{1}=(L_{1},(h_{1},\ldots ,h_{m}))$ and 
${\frak L}_{2}=(L_{2},(g_{1},\ldots ,g_{n})$ be two forward (resp. backward) 
self-similar systems. 
By Definition~\ref{d:fsscat} and \ref{d:bsscat}, 
 ${\frak L}_{1}$ is isomorphic to 
${\frak L}_{2}$ (indicated by ${\frak L}_{1}\cong {\frak L}_{2}$) 
if and only if $m=n$ and there exists a homeomorphism $\alpha : L_{1}\rightarrow 
L_{2}$ and a bijection $\tau :\{ 1,\ldots ,m\} \rightarrow \{ 1,\ldots ,m\} $ 
such that for each $j=1,\ldots, m$, 
$\alpha h_{j}=g_{\tau (j)}\alpha $ on $L_{1}$ 
(resp. $\alpha (h_{j}^{-1}(L_{1}))\subset g_{\tau (j)}^{-1}(L_{2})$ and 
$\alpha h_{j}=g_{\tau (j)}\alpha $ on $h_{j}^{-1}(L_{1})$).  
\end{rem}

We give several examples of forward or backward self-similar systems.
\begin{df}
Let $(X,d)$ be a metric space. 
We say that a map $f:X\rightarrow X$ is contractive (with respect to $d$) if 
there exists a number $0<s<1$ such that 
for each $x,y,\in X$, $d(f(x),f(y))\leq sd(x,y).$ A contractive map $f:X\rightarrow X$ is called 
a contraction.  
\end{df}
\begin{df}
Let $(X,d)$ be a complete metric space. 
For each $i=1,\ldots ,m$, 
let $h_{i}:X\rightarrow X$ be a contraction with respect to $d.$ 
By \cite[Theorem 1.1.4]{K}, there exists a unique 
non-empty compact subset $M$ of $X$ such that 
$(M,(h_{1},\ldots ,h_{m}))$ is a forward self-similar system. 
We denote this set $M$ by $M_{X}(h_{1},\ldots ,h_{m}).$  
The set $M_{X}(h_{1},\ldots ,h_{m})$ is called the 
{\bf attractor} or {\bf invariant set} of the 
iterated function system $\{ h_{1},\ldots ,h_{m}\} $ on $X.$  
\end{df}
\begin{df}
Let $X$ be a compact metric space. 
Let $G$ be a semigroup of continuous maps on $X.$ 
We set 
$F(G):=\{ z\in X\mid G \mbox{ is equicontinuous on a neighborhood of }z\} .$ 
For the definition of equicontinuity, see \cite{B}. 
The set $F(G)$ is called the {\bf Fatou set} of $G.$ 
Moreover, we set $J(G):= X\setminus F(G).$ The set $J(G)$ is called the 
{\bf Julia set} of $G.$ 
Furthermore, for a continuous map $g:X\rightarrow X$, 
we set $F(g):= F(\langle g\rangle )$ and $J(g):= J(\langle g\rangle ).$    
\end{df}
\begin{rem}
By the definition above, we have that $F(G)$ is open and $J(G)$ is compact.  
\end{rem}
By the definition above, it is easy to prove that the following 
Lemmas~\ref{openjlem} and \ref{openjfinlem} hold. 
\begin{lem}
\label{openjlem}
Let $X$ be a compact metric space. 
Let $G$ be a semigroup of continuous maps on $X.$ 
Suppose that for each $h\in G$, $h:X\rightarrow X$ is an open map. 
Then, for each $h\in G$, 
$h(F(G))\subset F(G)$ and $h^{-1}(J(G))\subset J(G).$ 
\end{lem}
\begin{lem}
\label{openjfinlem}
Let $X$ be a compact metric space. 
Let $G$ be a semigroup of continuous maps on $X.$ 
Suppose that $G$ is generated by a finite family $\{ h_{1},\ldots ,h_{m}\} $ 
of continuous maps on $X.$ 
Suppose that for each $j=1,\ldots,m$, $h_{j}:X\rightarrow X$ is an open and surjective map. 
Moreover, suppose $J(G)\neq \emptyset .$ 
Then, ${\frak L}:= (J(G),(h_{1},\ldots ,h_{m}))$ is a backward self-similar system.
\end{lem}
\begin{df}[\cite{HM,GR}]
We denote by $\CCI $ the Riemann sphere $\CC \cup \{ \infty \} .$  
A {\bf rational semigroup} is a semigroup 
generated by a family of non-constant rational maps on 
$\CCI $ with the semigroup operation being the 
functional composition. A 
{\bf polynomial semigroup } is a 
semigroup generated by a family of non-constant 
polynomial maps on $\CCI $.
\end{df}
\begin{rem}
\label{r:jisbss}
If a rational semigroup $G$ is generated by 
$\{ h_{1},\ldots ,h_{m}\} $ and if 
$J(G)\neq \emptyset$, then by Lemma~\ref{openjfinlem}, 
$(J(G),(h_{1},\ldots ,h_{m}))$ is a backward self-similar system. 

\end{rem}
\begin{rem}
\label{juliamcrem}
For each $j=1,\ldots ,m$, let 
$a_{j}\in \CC $ with $|a_{j}|>1$ and let $p_{j}\in \CC .$  
Moreover, let $h_{j}:\CCI \rightarrow \CCI $ be the map 
defined by $h_{j}(z)=a_{j}(z-p_{j})+p_{j}$ for each $z\in \CC .$ 
Let $G=\langle h_{1},\ldots ,h_{m}\rangle. $ 
Then, it is easy to see $\infty \in F(G).$ Hence $\emptyset \neq J(G)\subset \CC .$ 
From Lemma~\ref{openjfinlem} and \cite[Theorem 1.1.4]{K}, 
it follows that 
$J(G)=M_{\CC }(h_{1}^{-1},\ldots ,h_{m}^{-1}).$ 
\end{rem}
\begin{df}[\cite{HM}]
Let $G$ be a polynomial semigroup. 
We denote by $K_{1}(G)$ the set of points $z\in \CC $ satisfying that 
there exists a sequence $\{ g_{j}\} _{j\in \NN }$ of mutually distinct elements of 
$G$ such that $\{ g_{j}(z)\} _{j\in \NN }$ is bounded in $\CC .$ 
Moreover, we set $K(G):= \overline{K_{1}(G)}$, where the closure is taken in $\CCI .$ 
The set $K(G)$ is called the {\bf filled-in Julia set} of $G.$ 
Furthermore, for a polynomial $g$, we set 
$K(g):= K(\langle g\rangle ).$ 
\end{df}
\begin{rem}
It is easy to see that for each $g\in G$, $g^{-1}(K(G))\subset K(G).$ 
Moreover, if a polynomial semigroup $G$ is generated by a finite family 
$\{ h_{1},\ldots ,h_{m}\} $ and if $K(G)\neq \emptyset $, then 
${\frak L}=(K(G),(h_{1},\ldots ,h_{m}))$ is a backward self-similar 
system (\cite[Remark 3]{S01}). 
Furthermore, 
it is easy to see that if $G$ is generated by finitely many elements $h_{j},j=1,\ldots ,m$ 
such that $\deg (h_{j})\geq 2$ for each $j$, then $\emptyset \neq K(G)\subset \CC .$ 
\end{rem}
\begin{df}
For each $m\in \NN $, we set $\Sigma _{m}:=
\{ 1,\ldots ,m\} ^{\NN }$ endowed with the 
product topology. Note that $\Sigma _{m}$ is a compact metric space.   
Moreover, we set 
$\Sigma _{m}^{\ast }:= 
\bigcup _{j=1}^{\infty }\{1,\ldots ,m\} ^{j}$ (disjoint union). 
Let $X$ be a space and for each $j=1,\ldots ,m$, let $h_{j}:X\rightarrow X$ be a map.   
For a finite word $w=(w_{1},\ldots ,w_{k})\in 
\{ 1,\ldots ,m\} ^{k}$, 
we set 
$h_{w}:=h_{w_{k}}\circ \cdots \circ h_{w_{1}}$, 
 $\overline{w}=(w_{k},w_{k-1},\ldots ,w_{1})$, and   
$| w| :=k.$ For an element $w\in \Sigma _{m}$, we set $|w|=\infty .$ 
For an element $w\in \Sigma _{m}\cup \Sigma _{m}^{\ast }$, $|w|$ is called the word length of $w.$  
 Moreover, for any $w=(w_{1},w_{2},\ldots )\in \Sigma _{m}\cup \Sigma _{m}^{\ast }$ and 
 any $l\in \NN $ with $l\leq |w|$, we set 
 $w|l:=(w_{1},w_{2},\ldots ,w_{l})\in \{ 1,\ldots ,m\} ^{l}.$ 
Furthermore, for any $w=(w_{1},\ldots,  w_{k}\in \Sigma _{m}^{\ast }$ and 
any $\tau =(\tau _{1},\tau _{2},\ldots )\in \Sigma _{m}^{\ast }\cup \Sigma _{m}$, 
we set $w\tau =(w_{1},w_{2},\ldots ,w_{k},\tau _{1},\tau _{2},\ldots )\in \Sigma _{m}^{\ast }\cup \Sigma _{m}.$ 

\end{df}

\begin{df}
\label{Rfssdf}
Let $K$ be a non-empty compact metric space and 
let $h_{j}:K\rightarrow K$ be a continuous map for each 
$j=1,\ldots ,m.$ 
We set $$R_{K,f}(h_{1},\ldots ,h_{m})
:=\bigcap _{n=1}^{\infty }\bigcup _{w\in \Sigma _{m}^{\ast }: |w|=n}
h_{\overline{w}}(K).$$ 
\end{df}
\begin{lem}
\label{Rfsslem}
Under Definition~\ref{Rfssdf}, 
we have that $R_{K,f}(h_{1},\ldots ,h_{m})$ is non-empty and compact, 
$R_{K,f}(h_{1},\ldots ,h_{m})=\bigcup _{w\in \Sigma _{m}}\bigcap _{k=1}^{\infty }
h_{\overline{w|k}}(K)$, and  
${\frak L}:= (R_{K,f}(h_{1},\ldots ,h_{m}),(h_{1},\ldots ,h_{m}))$ is a 
forward self-similar system. 
\end{lem}
\begin{proof}
It is easy to see that $R_{K,f}(h_{1},\ldots ,h_{m})$ is non-empty and compact. 
Moreover, it is easy to see that 
$R_{K,f}(h_{1},\ldots ,h_{m})\supset \bigcup _{w\in \Sigma _{m}}\bigcap _{k=1}^{\infty }
h_{\overline{w|k}}(K).$ 
To show the opposite inclusion, let 
$x\in R_{K,f}(h_{1},\ldots ,h_{m}).$ 
Then for each $n\in \NN $ there exists a word $w^{n}\in \Sigma _{m}^{\ast }$ 
with $|w^{n}|=n$ and a point $y_{n}\in K$ such that 
$x=h_{w_{1}^{n}}\cdots h_{w_{n}^{n}}(y_{n}).$ 
Then, there exists an infinite word $w^{\infty }\in \Sigma _{m}$ and  
a sequence $\{ n_{k}\} _{k\in \NN }$ of positive integers with $n_{k}>k$ 
such that for each $k\in \NN $, $w^{n_{k}}|k=w^{\infty }|k.$ 
Hence, for each $k\in \NN $, 
$x=h_{w_{1}^{\infty }}\cdots h_{w_{k}^{\infty }}h_{w_{k+1}^{n_{k}}}
\cdots h_{w_{n_{k}}^{n_{k}}}(y_{n_{k}}).$ Therefore, 
$x\in \bigcap _{k=1}^{\infty }h_{w_{1}^{\infty }}\cdots h_{w_{k}^{\infty }}(K).$ 
Thus, we have shown $R_{K,f}(h_{1},\ldots ,h_{m})=\bigcup _{w\in \Sigma _{m}}\bigcap _{k=1}^{\infty }
h_{\overline{w|k}}(K).$ 
 From this formula, it is easy to see that 
$R_{K,f}(h_{1},\ldots ,h_{m})\supset \bigcup _{j=1}^{m}h_{j}
(R_{K,f}(h_{1},\ldots , h_{m})).$ 
In order to show the opposite inclusion, 
let $x\in R_{K,f}(h_{1},\ldots ,h_{m})=\bigcup _{w\in \Sigma _{m}}\bigcap _{k=1}^{\infty }
h_{\overline{w|k}}(K)$ be a point. 
Let $w\in \Sigma _{m}$ be an element such that 
$x\in  \bigcap _{k=1}^{\infty }h_{w_{1}}\cdots h_{w_{k}}(K).$ 
Then for each $k\in \NN $ with $k\geq 2$, there exists a point 
$y_{k}\in h_{w_{2}}\cdots h_{w_{k}}(K)$ such that 
$x=h_{w_{1}}(y_{k}).$ Since $K$ is a compact metric space, 
there exists a subsequence $\{ y_{k_{l}}\} _{l\in \NN }$ of 
$\{ y_{k}\} _{k\in\ \NN }$ and a point $y_{k_{\infty }}\in K$ 
such that $y_{k_{l}}\rightarrow y_{k_{\infty }}$ as $l\rightarrow \infty .$ 
Then, it is easy to see that $y_{k_{\infty }}\in  
\bigcap _{i=2}^{\infty }h_{w_{2}}\cdots h_{w_{i}}(K).$ 
Hence, $x=h_{w_{1}}(y_{k_{\infty }})\in 
h_{w_{1}}(\bigcup _{\tau \in \Sigma _{m}}\bigcap _{k=1}^{\infty }
h_{\overline{\tau |k}}(K)).$ Thus, we have proved Lemma~\ref{Rfsslem}. 
\end{proof}
\begin{df}
\label{Rbssdf}
Let $X$ be a metric space and let 
$h_{j}:X\rightarrow X$ be a continuous map for each 
$j=1,\ldots ,m.$ Let $K$ be a compact subset of $X$ 
and suppose that for each $z\in K$ and $j=1,\ldots ,m$, 
we have $h_{j}^{-1}(\{ z\} )\neq \emptyset $. Moreover, 
suppose that  
$\bigcup _{j=1}^{m}h_{j}^{-1}(K)\subset K.$ 
Then we set $R_{K,b}(h_{1},\ldots ,h_{m}):= 
\bigcap _{n=1}^{\infty }\bigcup _{w\in \Sigma _{m}^{\ast }:|w|=n} 
h_{w}^{-1}(K).$  
\end{df}
Using the argument similar to that in the proof of Lemma~\ref{Rfsslem}, 
we can easily prove the following lemma.  
\begin{lem}
\label{Rbsslem}
Under Definition~\ref{Rbssdf}, 
we have that $R_{K,b}(h_{1},\ldots ,h_{m})$ is non-empty and compact, 
$R_{K,b}(h_{1},\ldots ,h_{m})=\bigcup _{w\in \Sigma _{m}}
\bigcap _{k=1}^{\infty }
h_{w|k}^{-1}(K)$, and  ${\frak L}:= 
(R_{K,b}(h_{1},\ldots ,h_{m}),(h_{1},\ldots ,h_{m}))$ is a 
backward self-similar system. 
\end{lem}
\begin{df}
Let $X$ be a compact metric space and 
let $G$ be a semigroup of continuous maps on $X$. 
A non-empty compact subset $M$ of $X$ is said to be 
minimal for $(G,X)$ if 
$M$ is minimal with respect to the inclusion 
in the space of all non-empty compact subsets 
$K$ of $X$ 
satisfying 
that for each $g\in G$, $g(K)\subset K.$  
\end{df}
\begin{lem}
\label{minimalsetlem}
Let $X$ be a compact metric space and 
let $G$ be a semigroup of continuous maps on $X$. 
Then, we have the following. 
\begin{enumerate}
\item Let $K$ be a non-empty compact subset of $X$ such that 
for each $g\in G$, $g(K)\subset K.$ Then, 
there exists a minimal set $L$ for $(G,X)$ such that $L\subset K.$  
\item If, in addition to the assumptions of our lemma, 
$G$ is generated by a finite family $\{ h_{1},\ldots ,h_{m}\} $ 
of continuous maps on $X$, then 
for any minimal set $M$ for $(G,X)$, 
$(M,(h_{1},\ldots ,h_{m}))$ is a forward self-similar system.  
\end{enumerate}
\end{lem}
\begin{proof}
Statement 1 easily follows from Zorn's lemma. 
In order to show statement 2, 
suppose that $G=\langle h_{1},\ldots ,h_{m}\rangle $ and 
$M$ is a minimal set for $(G,X).$  
Since $M$ satisfies that $g(M)\subset M$ for each 
$g\in G$, we have $\bigcup _{j=1}^{m}h_{j}(M)\subset M.$ 
Let $K:=\bigcup _{j=1}^{m}h_{j}(M).$ 
Since $G=\langle h_{1},\ldots ,h_{m}\rangle $, 
we have that for each $g\in G$, $g(K)\subset K.$ 
Thus, by statement 1, there exists a minimal set 
$L$ for $(G,X)$ such that $L\subset K.$ 
By the minimality of $M$, it must hold that $L=M.$ 
Hence, $K=M.$  Therefore, we have proved statement 2. 
Thus, we have proved Lemma~\ref{minimalsetlem}. 
\end{proof}
\begin{rem}
\label{r:minsetimp}
It is very important to consider the minimal sets for rational semigroups when we investigate 
random complex dynamics as well as rational semigroups. 
In \cite{Srcd}, it is shown that for any Borel probability measure $\tau $ on 
the space Rat of non-constant rational maps, if supp$\,\tau $ is compact,   
$\bigcap _{g\in G_{\tau }}g^{-1}(J(G_{\tau }))=\emptyset $ and 
$J(G_{\tau })\neq \emptyset $, where $G_{\tau }$ denotes the rational semigroup generated by supp$\,\tau $, 
then there exist at most finitely many minimal sets $K_{1},\ldots ,K_{r}$ for $(G_{\tau },\CCI )$, and 
for each $z\in \CCI $, for $(\otimes _{j=1}^{\infty }\tau )$-a.e. $\gamma =(\gamma _{1},\gamma _{2},\ldots )
\in (\mbox{Rat})^{\NN }$, $d(\gamma _{n}\cdots \gamma _{1}(z), \bigcup _{j=1}^{r}K_{j})\rightarrow 0$ 
as $n\rightarrow \infty .$ Note that $K_{j}$ may meet $J(G_{\tau })$, and 
for a $g\in G_{\tau }$, $g|_{K_{j}}$ is neither contractive nor injective in general. 
Thus it is important to investigate forward self-similar systems whose generators may be neither contractive nor injective.    
\end{rem}
The above examples give us a natural and strong motivation to 
investigate forward or backward self-similar systems. 

 We now give some definitions which we need later. 
\begin{df}
Let ${\frak L}_{1}=(L_{1},(h_{1},\ldots ,h_{m}))$ and 
${\frak L}_{2}=(L_{2},(g_{1},\ldots ,g_{n}))$ be two forward or backward 
self-similar systems. 
We say that ${\frak L}_{1}$ is a subsystem of ${\frak L}_{2}$ if 
$L_{1}\subset L_{2}$ and there exists an injection $\tau :\{ 1,\ldots ,m\} 
\rightarrow \{ 1,\ldots ,n\} $ such that 
for each $j=1,\ldots ,m$, $h_{j}=g_{\tau (j)} .$ 
\end{df}
\begin{df}
Let ${\frak L} =(L, (h_{1},\ldots ,h_{m}))$ be a 
forward (resp. backward) self-similar system. 
A forward (resp. backward) self-similar system 
${\frak M}=(L,(g_{1},\ldots ,g_{m^{n}}))$ is said to be 
an $n$-th iterate of ${\frak L}$ if 
there exists a bijection $\tau :\{ 1,\ldots ,m^{n}\} 
\rightarrow \{ w\in \Sigma _{m}^{\ast } \mid |w|=n\} $ such that 
for each $j=1,\ldots ,m^{n}$, 
$g_{j}=h_{\tau (j)} .$  
\end{df}
\begin{df}
For a topological space $X$, we denote by 
 Con$(X)$ the set of all connected components of $X.$ 
\end{df}
\begin{df}
Let $X$ be a space. 
For any covering ${\cal U}=\{ U_{\lambda }\} _{\lambda \in \Lambda }$ of 
$X$, 
we denote by $N({\cal U})$ the nerve of ${\cal U}.$  
By definition, 
the vertex set of $N({\cal U})$ is equal to $\Lambda $. 
\end{df}
\begin{df}
Let ${\cal S}$ be an abstract simplicial complex. 
Moreover, we denote by $|{\cal S}|$ the realization (see \cite[p.110]{Sp}). 
As in \cite{Sp}, we embed the vertex set of ${\cal S}$ into $|{\cal S}|.$ 
\end{df}
We now define a new kind of cohomology theory for forward or backward 
self-similar systems.
\begin{df}
\label{bssicdf}
Let ${\frak L}=(L,(h_{1},\ldots ,h_{m}))$ be a backward self-similar 
system. 
\begin{enumerate} 

\item 
For each $x=(x_{1},x_{2},\ldots ,)\in \Sigma _{m} $,  
we set 
$ L_{x}:=\bigcap _{j=1}^{\infty }
h_{x_{1}}^{-1}\cdots h_{x_{j}}^{-1}L \ (\neq \emptyset ).$ 

\item 
For any $k\in \NN $, 
let ${\cal U}_{k}={\cal U}_{k}({\frak L})$ be the finite covering of $L$ 
defined as: 
$ {\cal U}_{k}:=\{ h_{w}^{-1}(L)\} _{w\in \Sigma _{m}^{\ast }: |w|=k} .$ 
We denote by $N_{k}$ or $N_{k}({\frak L})$ the nerve  $N({\cal U}_{k})$ of 
${\cal U}_{k}.$ Thus $N_{k}$ is a simplicial complex such that 
the vertex set is equal to $\{ w\in \Sigma _{m}^{\ast }\mid |w|=k\} $ and 
mutually distinct $r$ elements $w^{1},\ldots ,w^{r}\in \Sigma _{m}^{\ast }$ with 
$|w^{1}|=\cdots =|w^{r}|=k$ make an $(r-1)$-simplex of $N_{k}$ if and only if 
$\bigcap _{j=1}^{r}h_{w^{j}}^{-1}(L)\neq \emptyset .$    
 Let $\varphi _{k}:N_{k+1}\rightarrow 
N_{k}$ be the simplicial map 
defined as:
$(w_{1},\ldots ,w_{k+1})\mapsto (w_{1},\ldots ,w_{k})$   
for each 
$ (w_{1},\ldots ,w_{k+1})\in \{ 1,\ldots ,m\} ^{k+1}.$ 
Moreover, for each $k,l\in \NN $ with $l>k$, we denote by 
$\varphi _{l,k}:N_{l}\rightarrow N_{k}$ the 
composition $\varphi _{l-1}\circ \cdots \circ \varphi _{k}.$ 
Then, $\{ N_{k}, \varphi _{l,k}\} _{k,l\in \NN, l>k}$ forms an 
inverse system of 
simplicial complexes.  
\item 
Let $ \{ (\varphi _{k})_{\ast }: \mbox{Con}(|N_{k+1}|)
\rightarrow \mbox{Con}(|N_{k}|) \} _{k\in \NN }$ be 
the inverse system induced by 
$\{ (\varphi _{k})_{\ast }\} _{k}.$      
\item 
\label{bssicdf2}
Let $R$ be a $\ZZ $ module and let $p\in \NN \cup \{ 0\} .$ 
We set 
$\check{H}^{p}(L,(h_{1},\ldots ,h_{m});R)_{k}:=
H^{p}(N_{k}; R).$ This is called 
{\bf the $p$-th interaction cohomology group 
of backward self-similar system ${\frak L}=(L, (h_{1},\ldots ,h_{m}) )$ 
at $k$-th stage 
with coefficients $R$ }. 
We sometimes use the notation 
$\check{H}^{p}({\frak L};R)_{k}$ to denote  
the above $\check{H}^{p}(L,(h_{1},\ldots ,h_{m});R)_{k}.$ 
Similarly, we set 
$\check{H}_{p}(L,(h_{1},\ldots ,h_{m});R)_{k}$ $ :=
H_{p}(N_{k}; R).$ This is called 
{\bf the $p$-th interaction homology group 
of backward self-similar system ${\frak L}=(L, (h_{1},\ldots ,h_{m}) )$ 
at $k$-th stage 
with coefficients $R$ }. 
We sometimes use the notation 
$\check{H}_{p}({\frak L};R)_{k}$ to denote  
the above $\check{H}_{p}(L,(h_{1},\ldots ,h_{m});R)_{k}.$  
\item 
Let $R$ be a $\ZZ $ module 
and let $p\in \NN \cup \{ 0\} .$ 
We denote by $\check{H}^{p}(L,(h_{1},\ldots ,h_{m});R)$ the direct limit of 
the direct system $\{ H^{p}(N_{k};R), \varphi _{l,k}^{\ast }\} _{l,k\in \NN , l>k} $ of $\ZZ $ modules.  
This is called the {\bf $p$-th interaction cohomology group of ${\frak L}=(L,(h_{1},\ldots ,h_{m}))$ 
with coefficients $R.$ } 
We sometimes use the notation 
$\check{H}^{p}({\frak L};R)$ in order to denote  
the above cohomology group $\check{H}^{p}(L,(h_{1},\ldots ,h_{m});R).$ 
\item 
We denote by $\mu _{k,p}: \check{H}^{p}(\frak{L};R)_{k}\rightarrow \check{H}^{p}(\frak{L};R)$ the 
canonical projection. 
\item Similarly, for any $\ZZ $ module $R$, 
we denote by $\check{H}_{p}(L,(h_{1},\ldots ,h_{m});R)$ the 
inverse limit of the inverse system 
$\{ H_{p}(N_{k};R), (\varphi _{l,k})_{\ast }\} _{l,k\in \NN ,l>k}$ of $\ZZ $ modules. 
This is called the 
{\bf $p$-th interaction homology group of   
${\frak L}=(L,(h_{1},\ldots ,h_{m}))$ with coefficients $R.$}
We sometimes use the notation 
$\check{H}_{p}({\frak L};R)$ in order to denote  
the above homology group $\check{H}_{p}(L,(h_{1},\ldots ,h_{m});R).$
\item For each $p\in \NN $, $k\in \NN $ and $w\in \Sigma _{m}$, we 
set $\check{\pi }_{p}({\frak L},w)_{k}:= 
\pi _{p}(|N_{k}|, w|k)$ and 
$\check{\pi }_{p}({\frak L},w):= \varprojlim _{k} \pi _{p}(|N_{k}|,w|k).$ 
We call $\check{\pi }_{p}({\frak L},w)_{k}$ the $p$-th {\bf interaction homotopy group 
of ${\frak L}$ at $k$-th stage} with base $w$ and 
we call $\check{\pi }_{p}({\frak L},w)$ the $p$-th {\bf interaction homotopy group 
of ${\frak L}$ with base $w.$} 
\end{enumerate}
\end{df}
\begin{df}
\label{fssicdf}
Let ${\frak L}=(L,(h_{1},\ldots ,h_{m}))$ be a forward self-similar system. 
For each $x\in \Sigma _{m}$, we set 
$L_{x}:= \bigcap _{j=1}^{\infty }h_{x_{1}}\cdots h_{x_{j}}(L).$ 
For any $k\in \NN $, let ${\cal U}_{k}={\cal U}_{k}({\frak L})$ be the finite covering of $L$ 
defined as:  
${\cal U}_{k}:= \{ h_{\overline{w}}(L)\} _{w\in \Sigma _{m}^{\ast }: |w|=k} .$ 
We denote by $N_{k}$ or $N_{k}({\frak L})$ the nerve $N({\cal U}_{k})$ of ${\cal U}_{k}.$ 
Thus $N_{k}$ is a simplicial complex such that 
the vertex set is equal to $\{ w\in \Sigma _{m}^{\ast }\mid |w|=k\} $ and 
mutually distinct $r$ elements $w^{1},\ldots ,w^{r}\in \Sigma _{m}^{\ast }$ with 
$|w^{1}|=\cdots =|w^{r}|=k$ make an $(r-1)$-simplex of $N_{k}$ if and only if 
$\bigcap _{j=1}^{r}h_{\overline{w^{j}}}(L)\neq \emptyset .$    
Let $\varphi _{k}:N_{k+1}\rightarrow N_{k}$ be the simplicial map defined as: 
$(w_{1},\ldots ,w_{k+1})\mapsto (w_{1},\ldots ,w_{k})$ for each $(w_{1},\ldots ,w_{k+1})\in 
\{ 1,\ldots ,m\} ^{k+1}.$ Moreover, we set $\varphi _{l,k}:= \varphi _{l-1}\circ \cdots \circ \varphi _{k}.$ 
We define the $p$-th interaction cohomology group 
$\check{H}^{p}({\frak L};R)$ and  the $p$-th interaction homology group $\check{H}_{p}({\frak L};R)$ 
as in Definition~\ref{bssicdf}. 
Moreover, we define the $p$-th interaction homotopy group $\check{\pi }_{p}({\frak L},w)_{k}$ of ${\frak L}$ 
at $k$-th stage with base $w\in \Sigma _{m}$, 
the $p$-th interaction homotopy group $\check{\pi }_{p}({\frak L},w)$ of ${\frak L} $ with base 
$w$,  
the $p$-th interaction cohomology group $\check{H}^{p}({\frak L};R)_{k}$ of ${\frak L}$ at $k$-th stage, 
and $p$-th interaction homology group $\check{H}_{p}({\frak L};R)_{k}$ of ${\frak L}$ at $k$-th stage,     
as in Definition~\ref{bssicdf}. 
Furthermore, 
we denote by $\mu _{k,p}: \check{H}^{p}(\frak{L};R)_{k}\rightarrow \check{H}^{p}(\frak{L};R)$ the 
canonical projection. 
\end{df}
\begin{df}
Let ${\cal A}=\{ A_{\lambda }\} _{\lambda \in \Lambda _{1}}$ 
and ${\cal B}=\{ B_{\mu }\} _{\mu \in \Lambda _{2}}$ be two 
coverings of a topological space $L.$ 
We say that ${\cal B}$ is a refinement of ${\cal A}$ if 
there exists a map $r_{{\cal A},{\cal B}}:\Lambda _{2}\rightarrow \Lambda _{1}$ 
such that $B_{\mu }\subset A_{r_{{\cal A},{\cal B}}(\mu )}$ for each 
$\mu \in \Lambda _{2}.$ The $r_{{\cal A},{\cal B}}$ is called the 
refining map. 
If ${\cal B}$ is a refinement of ${\cal A}$, we write 
${\cal A}\preceq {\cal B}.$  
\end{df}
\begin{rem}
\label{r:refvp}
Let ${\frak L}=(L,(h_{1},\ldots ,h_{m}))$ be a forward (resp. backward) self-similar system. 
Let $r_{{\cal U}_{k},{\cal U}_{k+1}}: \{ w\in \Sigma _{m}^{\ast } \mid |w|=k+1\} \rightarrow 
\{ w\in \Sigma _{m}^{\ast } \mid |w|=k\} $ be the map defined by 
$(w_{1},\ldots ,w_{k+1})\mapsto (w_{1},\ldots ,w_{k}).$ 
If $w=(w_{1},\ldots, w_{k+1})\in \{ 1,\ldots, m\} ^{k+1}$, 
then $h_{\overline{w}}(L)\subset h_{\overline{w|k}}(L)$ (resp. 
$h_{w}^{-1}(L)\subset h_{w|k}^{-1}(L)$).  Thus 
for each $k$, 
${\cal U}_{k}\preceq {\cal U}_{k+1}$ with refining map 
$r_{{\cal U}_{k},{\cal U}_{k+1}}.$  
Moreover, 
the simplicial map $(r_{{\cal U}_{k},{\cal U}_{k+1}})_{\ast }: N_{k+1}\rightarrow N_{k}$ induced by 
$r_{{\cal U}_{k},{\cal U}_{k+1}}$ is equal to 
the simplicial map $\varphi _{k}: N_{k+1}\rightarrow N_{k}.$ 
\end{rem}

From the  definition of interaction (co)homology groups  
and the continuity theorem for \v{C}ech (co)homology (\cite{W}), it is easy to prove the following lemma. 
\begin{lem}
\label{bssiclem1}
Let ${\frak L}$ be a forward or backward self-similar system and let $R$ be a $\ZZ $ module. Then 
$ \check{H}^{p}({\frak L};R) 
\cong \check{H}(\varprojlim _{k}|N_{k}|;R)$ and   
$ \check{H}_{p}({\frak L};R) 
\cong   \check{H}_{p}(\varprojlim _{k}|N_{k}|;R)$.  
\end{lem}
\begin{ex}[{\bf Sierpi\'{n}ski gasket}]
\label{ex:SG1}
Let $p_{1},p_{2},p_{3}\in \CC $ be mutually distinct three points 
such that $p_{1}p_{2}p_{3}$ makes an equilateral  
triangle. Let $h_{i}(z):=\frac{1}{2}(z-p_{i})+p_{i}$, 
for each $i=1,2,3.$ 
Let $L=M_{\CC }(h_{1},h_{2},h_{3}).$ 
Then, 
$L$ is equal to the Sierpi\'{n}ski gasket (\cite{K}, see Figure~\ref{fig:sgasket}). 
We consider the forward self-similar system 
${\frak L}=(L,(h_{1},h_{2} ,h_{3}))$. 
We see that the set of one-dimensional simplexes of $N_{1}$ is equal to 
$\{ \{ 1,2\} ,\{ 2,3\} , \{ 3,1\} \} $ and for each $r\geq 2$, there exists no 
$r$-dimensional simplexes of $N_{1}.$ 
Moreover, it is easy to see that the set of one-dimensional simplexes of $N_{2}$ is 
equal to 
$\{ \{ (1,1),(1,2)\} , \{ (1,2), (1,3)\} , \{ (1,3), (1,1)\} ,$
$ \{ (2,1),(2,2)\} , \{ (2,2), (2,3)\} , \{ (2,3), (2,1)\} ,$ 
$\{ (3,1),(3,2)\} , \{ (3,2), (3,3)\} , \{ (3,3), (3,1)\} ,$  
$\{ (1,2), (2,1)\} , \{ (2,3), (3,2)\} , \{ (3,1), (1,3)\} \} $ and 
for each $r\geq 2$ there exists no $r$-dimensional simplexes of $N_{2}$ 
(see Figure~\ref{fig:sierpn1n2-2}).  
Thus for each $\ZZ $ module $R$, 
$\check{H}^{1}({\frak L};R)_{1}=R,$ 
$\check{H}^{r}({\frak L};R)_{1}=0\ (\forall r\geq 2), $ 
$\check{H}^{1}({\frak L};R)_{2}=R^{4}, $ and 
$\check{H}^{r}({\frak L};R)_{2}=0\ (\forall r\geq 2).$  
\begin{figure}[htbp]
\caption{The Sierpi\'{n}ski gasket}
\ \ \ \ \ \ \ \ \ \ \ \ \ \ \ \ \ \ \ \ \ \ \ \ \ \ \ \ \ \ \ 
\ \ \ \ \ \ \ \ \ \ \ \ \ \ \ \ \ \ \ \ \ \ \ 
\includegraphics[width=2.3cm,width=2.3cm]{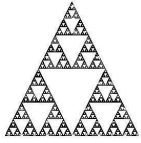}
\label{fig:sgasket}
\end{figure}
\begin{figure}[htbp]
\caption{The figure of $\varphi _{1}: N_{2}\rightarrow N_{1}$ for the system of the Sierpi\'{n}ski gasket}
\ \ \ \ \ \ \ \ \ \ \ \ \ \ \ \ \ \ \ \ \ \ \ \ \ \ \ \ \ \ \ \ \ \ \ \ \ \ \ \ \ \ 
\includegraphics[width=5.0cm,width=5.0cm]{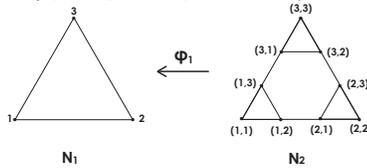}
\label{fig:sierpn1n2-2}
\end{figure}
\end{ex}
\begin{rem}
\label{r:functor}
${\frak L}\mapsto \{ N_{k}({\frak L}),\varphi _{l,k}\} _{l,k\in \NN ,l>k}$ 
is a covariant functor from the category of backward self-similar systems to the category 
of inverse systems of simplicial complexes. 
For any $\ZZ $ module $R$ and any $p\in \NN \cup \{ 0\} $, 
${\frak L}\mapsto \{ \check{H}_{p}({\frak L};R)_{k}, (\varphi _{l,k})_{\ast }\} _{l,k\in \NN , l>k}$ 
is a covariant functor from the category of backward self-similar systems to the category 
of inverse systems of $\ZZ $ modules, 
${\frak L}\mapsto \check{H}_{p}({\frak L};R)$ is a covariant functor from the category of 
backward self-similar systems to the category of $\ZZ $ modules, 
 ${\frak L}\mapsto \{ \check{H}^{p}({\frak L};R)_{k}, (\varphi _{l,k})^{\ast }\} _{l,k\in \NN , l>k}$ 
is a contravariant functor from the category of backward self-similar systems to the category 
of direct systems of $\ZZ $ modules, and 
${\frak L}\mapsto \check{H}^{p}({\frak L};R)$ is a contravariant functor from the category of 
backward self-similar systems to the category of $\ZZ $ modules. 
Thus the isomorphism classes of $\{ N_{k}({\frak L}),\varphi _{l,k}\} _{l,k\in \NN ,l>k} $, 
$\{ \check{H}_{p}({\frak L};R)_{k}, (\varphi _{l,k})_{\ast }\} _{l,k\in \NN , l>k}$, 
$\check{H}_{p}({\frak L};R)$, 
$\{ \check{H}^{p}({\frak L};R)_{k}, (\varphi _{l,k})^{\ast }\} _{l,k\in \NN , l>k}$, 
and  $\check{H}^{p}({\frak L};R)$ are invariant under the isomorphisms of backward self-similar 
systems. The same statements as above hold for forward self-similar systems. 
\end{rem}
\begin{rem}
\label{r:semigrinv}
Let ${\frak L}_{1}=(L,(h_{1},\ldots ,h_{m}))$ and ${\frak L}_{2}=(L,(g_{1},\ldots ,g_{n}))$  
be two backward self-similar systems such that $\langle h_{1},\ldots ,h_{m}\rangle =\langle g_{1},\ldots ,g_{n}\rangle .$ 
Then, by the definition of the interaction (co)homology, it is easy to see that 
there exist isomorphisms $\check{H}^{\ast }({\frak L}_{1};R)\cong \check{H}^{\ast }({\frak L}_{2};R)$ 
and $\check{H}_{\ast }({\frak L}_{1};R)\cong \check{H}_{\ast }({\frak L}_{2};R).$ 
Similar statement holds for two forward self-similar systems. 
\end{rem}
\noindent {\bf Notation:} 
Let $(X,d)$ be a metric space. Let $A$ be a non-empty subset of $X.$ 
Let $\delta >0.$ We set 
$B(A,\delta ):=\{ x\in X\mid d(y,A)<\delta \} .$  
\begin{df}
Let $(X,d)$ be a metric space. Let 
${\cal A}:=\{ L_{\lambda }\} _{\lambda \in \Lambda }$ 
be a covering of $X.$ 
For each $\delta >0$, we set 
${\cal A}^{\delta } := \{ B(L_{\lambda },\delta )\} _{\lambda \in \Lambda }$ 
and we denote by $\psi _{{\cal A},\delta } :
N({\cal A})\rightarrow N({\cal A}^{\delta })$ the simplicial map 
induced by the refinement $L_{\lambda }\subset B(L_{\lambda },\delta ), \lambda 
\in \Lambda .$
\end{df}
\begin{lem}
\label{naturalhomlem}
Let $(L,d)$ be a compact metric space. 
Let ${\cal A}=\{ L_{i}\} _{i=1}^{r}$ be a 
finite covering of $L$ such that for each $i=1,\ldots ,r$, 
$L_{i}$ is a non-empty compact subset of $L.$ Then, 
we have the following.
\begin{enumerate}
\item \label{naturalhomlem1}
There exists a number $\delta ({\cal A})>0$ such that for each 
$0<\delta <\delta ({\cal A})$,  
$\psi _{{\cal A},\delta }: N({\cal A})\rightarrow  N({\cal A}^{\delta })$ 
is a simplicial isomorphism.  

\item 
\label{naturalhomlem2} 
Let ${\cal B}=\{ M_{j}\} _{j=1}^{l}$ be another finite covering 
of $L$ such that 
for each $j=1,\ldots,l$, $M_{j}$ is a non-empty compact subset of 
$L.$ Assume that there exists a map 
$\beta _{{\cal A},{\cal B}}: \{ 1,\ldots ,r\} \rightarrow 
\{ 1,\ldots ,l\} $ 
such that $M_{j}\subset L_{\beta _{{\cal A},{\cal B}}(j)}$ for each 
$j=1,\ldots ,l.$ Then, there exists a $\delta _{0}>0$ such that 
for each $0<\delta <\delta _{0}$, we have the following 
(i),(ii), and (iii):   
(i) $B(M_{j},\delta )\subset B(L_{\beta _{{\cal A},{\cal B}}(j)},\delta )\ 
(j=1,\ldots ,l)$,  
(ii) the diagram 
$$
\begin{CD}
N({\cal B}) @>{\psi _{{\cal B},\delta }}>>N({\cal B}^{\delta })\\ 
@V{(\beta _{{\cal A},{\cal B}})_{\ast }}VV 
@VV{(\beta _{{\cal A},{\cal B}})_{\ast }}V\\ 
N({\cal A}) @>{\psi _{{\cal A},\delta }}>>N({\cal A}^{\delta })
\end{CD}
$$ 
commutes where 
$(\beta _{{\cal A},{\cal B}})_{\ast }: N({\cal B})\rightarrow N({\cal A})$ 
and 
$(\beta _{{\cal A},{\cal B}})_{\ast }: 
N({\cal B}^{\delta })\rightarrow N({\cal A}^{\delta })$ 
are simplicial maps induced by 
$\beta _{{\cal A},{\cal B}}:\{ 1,\ldots ,r\} \rightarrow 
\{ 1,\ldots ,l\} $,
 and (iii) the 
simplicial maps 
$\psi _{{\cal B},\delta }: N({\cal B})\rightarrow N({\cal B}^{\delta })$ 
and 
$\psi _{{\cal A},\delta }: N({\cal A})\rightarrow N({\cal A}^{\delta })$ 
are isomorphisms. 
\end{enumerate} 
\end{lem}
\begin{proof}
First, we will show statement \ref{naturalhomlem1}. 
If $\bigcap _{i=1}^{r}L_{i}\neq \emptyset $, 
then for any $\delta >0$, $\psi _{{\cal A},\delta }$ is an 
isomorphism. Hence we may assume that 
$\bigcap _{i=1}^{r}L_{i}=\emptyset .$ 
Let $(i_{1},\ldots ,i_{r})\in \{ 1,\ldots ,r\} ^{r}$ be any element 
with $\bigcap _{t=1}^{r}L_{i_{t}}=\emptyset .$  
Then there exists a $\delta =\delta (i_{1}, \ldots ,i_{r})>0 $ 
such that $\bigcap _{t=1}^{r}B(L_{i_{t}},\delta )=\emptyset .$ 
Let 
$$\delta ({\cal A}):= \min 
\{ \delta (i_{1},\ldots, i_{r})\mid (i_{1},\ldots i_{r})\in 
\{ 1,\ldots ,r\} ^{r},\ \bigcap _{t=1}^{r}L_{i_{t}}=\emptyset \} .$$ 
Then, $\delta ({\cal A})>0.$ Hence, 
for each $0<\delta <\delta ({\cal A})$, 
if $\bigcap _{t=1}^{r}L_{i_{t}}=\emptyset $, then 
$\bigcap _{t=1}^{r}B(L_{i_{t}},\delta )=\emptyset .$ 
Therefore, statement \ref{naturalhomlem1} holds. 
 Statement \ref{naturalhomlem2} follows easily 
 from statement \ref{naturalhomlem1}. 
\end{proof}
\begin{rem}
\label{nathomrem} 
Let ${\frak L}=(L,(h_{1},\ldots ,h_{m}))$ be a forward or backward self-similar system. 
Let $R$ be a $\ZZ $ module and let 
$\check{H}^{p}(L;R)$ be the $p$-th \v{C}ech cohomology group of $L$ with 
coefficients $R.$ Since $\check{H}^{p}(L;R)=
\varinjlim _{{\cal W}} H^{p}(N({\cal W});R)$, where 
${\cal W}$ runs over all open coverings of $L$, 
Lemma~\ref{naturalhomlem} implies that 
for each $k\in \NN $, there exists a 
homomorphism 
$\Psi _{{\cal U}_{k}}:H^{p}(N_{k};R)\rightarrow 
\check{H}^{p}(L;R)$ induced 
by $\psi _{{\cal U}_{k},\delta }.$ Using Lemma~\ref{naturalhomlem} again, 
$\{ \Psi _{{\cal U}_{k}}\} _{k\in \NN }$ induces a natural 
homomorphism 
\begin{equation}
\label{Psidf}
\Psi : \check{H}^{p}({\frak L}; R)\rightarrow \check{H}^{p}(L;R).
\end{equation}  
\end{rem} 
\begin{rem}
\label{r:Psiiso}
Suppose that either (a) ${\frak L}=(L,(h_{1},\ldots ,h_{m}))$ is a forward 
self-similar system such that for each $j=1,\ldots ,m$, 
$h_{j}:L\rightarrow L$ is a contraction, 
or (b) ${\frak L}=(L,(h_{1},\ldots ,h_{m}))$ is a 
backward self-similar system such that for each $j=1,\ldots ,m$, 
$h_{j}^{-1}:L\rightarrow L$ is well-defined and 
$h_{j}^{-1}: L\rightarrow L$ is a contraction. 
Then, for any $p$ and any $\ZZ $ module $R$, 
$\Psi :\check{H}^{p}({\frak L};R)\rightarrow \check{H}^{p}(L;R)$ 
is an 
isomorphism. For, given an open covering ${\cal W}$, there exists 
a $k\in \NN $ and a $\delta >0$ such that 
${\cal W}\preceq N_{k}^{\delta }.$ It means that 
$\{ H^{p}(N_{k}^{\delta };R) \} _{k,\delta }$ is cofinal in 
$\{ H^{p}(N({\cal W});R) \} _{{\cal W}}.$ From Lemma~\ref{naturalhomlem}, it follows that  
$\Psi :\check{H}^{p}({\frak L};R)\rightarrow \check{H}^{p}(L;R)$ is an 
isomorphism. Similarly, $\check{H}_{p}({\frak L};R)$ and $\check{H}_{p}(L;R)$ 
are naturally isomorphic. Moreover, $\check{\pi }_{r}({\frak L},w)$ and $\check{\pi }_{r}(L,x_{0})$ 
($x_{0}\in h_{w_{1}}\cdots h_{w_{k}}(L)$ for each $k\in \NN $)  
are naturally isomorphic, where $\check{\pi }_{r}(L,x_{0})$ denotes the $r$-th 
shape group of $L$ with base point $x_{0}$. (For the definition of shape groups, see \cite{MS}.)  

However, $\Psi $ is not an isomorphism in general. 
In fact, $\Psi $ may not be a monomorphism (see Proposition~\ref{rgeq2exprop}). 
Similarly, $\check{H}_{p}({\frak L};R)$ and $\check{H}_{p}(L;R)$ 
may not be isomorphic in general (Example~\ref{ex:finsetic}). Moreover, $\check{\pi }_{r}({\frak L},w)$ and $\check{\pi }_{r}(L,x_{0})$ 
($x_{0}\in h_{w_{1}}\cdots h_{w_{k}}(L)$ for each $k\in \NN $)  
may not be isomorphic in general (Example~\ref{ex:finsetic}).  
\end{rem}
%
\section{Main results}
\label{Main}
In this section, we present the main results of this paper. 
The proofs of the results are given in section~\ref{Proofs}.  
\subsection{General results}
\label{General}
In this subsection, we present some general results on 
the $0$-th and the first interaction (co)homology groups of 
forward or backward self-similar systems. The proofs are given in 
section~\ref{ProofsGeneral}.  

We investigate the space of all connected components of an invariant set 
of a forward or backward self-similar system. This is related to 
the $0$-th interaction (co)homology groups of forward or backward self-similar systems. 
Note that it is a new point of view to study the above space. 
As an application,     
we generalize and further develop the essence of 
the well-known result (Theorem~\ref{weakconnthm}) on the 
necessary and sufficient condition for the invariant sets of the  
forward self-similar systems to be connected.   
\begin{rem}
\label{r:n1con}
Let ${\frak L}=(L,(h_{1},\ldots ,h_{m}))$ be a forward (resp. backward) self-similar 
system. Then $|N_{1}|$ is connected if and only if 
for each $i,j\in \{ 1,\ldots ,m\} $, there exists a sequence $\{ i_{t}\} _{t=1}^{r}$ 
in $\{ 1,\ldots ,m\} $ such that $i_{1}=i,i_{t}=j$, and 
$h_{i_{t}}(L)\cap h_{i_{t+1}}(L)\neq \emptyset $ 
(resp. $h_{i_{t}}^{-1}(L)\cap h_{i_{t+1}}^{-1}(L)\neq \emptyset $) for each 
$t=1,\ldots ,s-1.$ 
\end{rem}
\begin{thm}
\label{precompojth}
Let ${\frak L}=(L,(h_{1},\ldots ,h_{m}))$ be a backward self-similar system  
such that $L_{x}$ is connected for each $x\in \Sigma _{m}.$  
Let $R$ be a field. Then, we have the following.
\begin{enumerate}
\item \label{precompojth1}
There exists a bijection:
$ \Phi : \varprojlim _{k}\mbox{{\em Con}}(|N_{k}|)
\cong  \mbox{{\em Con}}(L)$, where,  
the map $\Phi $ is defined as follows: 
let $B=(B_{k})_{k}\in  \varprojlim _{k}\mbox{{\em Con}}(|N_{k}|)$ 
where $B_{k}\in \mbox{{\em Con}}(|N_{k}|)$ 
and $(\varphi _{k})_{\ast }(B_{k+1})=B_{k}$ for each $k.$ 
Take a point $x\in \Sigma _{m}$ 
 such that $(x_{1},\ldots ,x_{k})\in 
 B_{k}$ for each $k.$ Take an element 
 $C\in \mbox{{\em Con}}(L)$ such that 
 $L_{x}\subset C.$ Let 
 $\Phi (B)=C.$ 
 
\item \label{precompojth2}
$L$ is connected if and only if 
$|N_{1}|$ is connected. (See Remark~\ref{r:n1con}.)

\item \label{precompojth3}
$\sharp \mbox{{\em Con}}(|N_{k}|)
\leq \sharp \mbox{{\em Con}}(|N_{k+1}|)$, for each 
$k\in \NN .$ Furthermore,  
$\{ \sharp \mbox{{\em Con}}(|N_{k}|)\} _{k\in \NN }$
is bounded if and only if $\sharp \mbox{{\em Con}}(L)<\infty .$ 
If $\sharp \mbox{{\em Con}}(L)<\infty $,  
 then 
$\lim\limits _{k\rightarrow \infty }
\sharp \mbox{{\em Con}}(|N_{k}|)=
\sharp \mbox{{\em Con}}(L).$
\item  \label{precompojth6}
$\dim _{R}\check{H}^{0}({\frak L}; R)<\infty $ 
if and only if 
$\sharp \mbox{{\em Con}}(L)<\infty .$ 
\item \label{precompojth7}
If $\dim _{R}\check{H}^{0}({\frak L};R)<\infty $, 
then $\dim _{R}\check{H}^{0}({\frak L};R)=
\sharp \mbox{{\em Con}}(L)$ and 
$\Psi :\check{H}^{0}({\frak L};R)\rightarrow 
\check{H}^{0}(L;R)$ is an isomorphism. 
\item \label{precompojth4}
Suppose that $m=2$ and $L$ is disconnected.   
 Then, $h_{1}^{-1}(L)\cap h_{2}^{-1}(L)=\emptyset $,  
there exists a bijection 
 $ \mbox{{\em Con}}(L)\cong \Sigma _{2}$, 
 and $\sharp \mbox{{\em Con}}(L)>\aleph _{0}.$ 
\item \label{precompojth5}
Suppose that $m=3$ and $L$ is disconnected. Then,  
$\sharp \mbox{{\em Con}}(L)\geq \aleph _{0}$ and 
there exists a $j\in \{ 1,2,3\} $ such that 
$L_{(j)^{\infty }}$ is a connected component of $L$, where 
$(j)^{\infty }:= (j,j,j,\ldots )\in \Sigma _{3}.$ 

\end{enumerate}

\end{thm}
\begin{thm}
\label{ssprecompojth}
Let ${\frak L}=(L,(h_{1},\ldots ,h_{m}))$ be a forward self-similar system  
such that $L_{x}$ is connected for each $x\in \Sigma _{m}.$  
Let $R$ be a field. Then, we have the following.
\begin{enumerate}
\item \label{ssprecompojth1}
There exists a bijection:
$ \Phi : \varprojlim \mbox{{\em Con}}(|N_{k}|)
\cong  \mbox{{\em Con}}(L).$ 
 
\item \label{ssprecompojth2}
$L$ is connected if and only if 
$|N_{1}|$ is connected. (See Remark~\ref{r:n1con}.)

\item \label{ssprecompojth3}
$\sharp \mbox{{\em Con}}(|N_{k}|)
\leq \sharp \mbox{{\em Con}}(|N_{k+1}|)$, for each 
$k\in \NN .$ Furthermore,  
$\{ \sharp \mbox{{\em Con}}(|N_{k}|)\} _{k\in \NN }$
is bounded if and only if $\sharp \mbox{{\em Con}}(L)<\infty .$ 
If $\sharp \mbox{{\em Con}}(L)<\infty $,  
 then 
$\lim\limits _{k\rightarrow \infty }
\sharp \mbox{{\em Con}}(|N_{k}|)=
\sharp \mbox{{\em Con}}(L).$ 
\item  \label{ssprecompojth6}
$\dim _{R}\check{H}^{0}({\frak L}; R)<\infty $ 
if and only if 
$\sharp \mbox{{\em Con}}(L)<\infty .$ 
\item \label{ssprecompojth7}
If $\dim _{R}\check{H}^{0}({\frak L};R)<\infty $, 
then $\dim _{R}\check{H}^{0}({\frak L};R)=
\sharp \mbox{{\em Con}}(L) $ and 
$\Psi :\check{H}^{0}({\frak L};R)\rightarrow \check{H}^{0}(L;R)$ is an 
isomorphism. 

\item \label{ssprecompojth4}
If $m=2$ and $L$ is disconnected,  
 then $h_{1}(L)\cap h_{2}(L)=\emptyset .$ 
\item \label{ssprecompojth4-1}
If $m=2$, $h_{j}:L\rightarrow L$ is injective for each $j=1,2,$ 
and $L$ is disconnected, then 
there exists a bijection 
 $ \mbox{{\em Con}}(L)\cong \Sigma _{2}$
and $\sharp ${\em Con}$(L)>\aleph _{0}.$ 
\item \label{ssprecompojth5}
If $m=3$, $h_{j}:L\rightarrow L$ is injective for each $j=1,2,3$, and $L$ is disconnected, then 
$\sharp \mbox{{\em Con}}(L)\geq \aleph _{0}$ and there exists a 
$j\in \{ 1,2,3\} $ such that 
$L_{(j)^{\infty }}$ is a connected component of $L$, 
where $(j)^{\infty }:= (j,j,j,\ldots )\in \Sigma _{3}.$ 
\end{enumerate}
\end{thm}
\begin{rem}
\label{r:slx1}
Let ${\frak L}=(L,(h_{1},\ldots ,h_{m}))$ be a forward self-similar system. 
If each $h_{j}:L\rightarrow L $ is a contraction, then for each $x\in \Sigma _{m}$, 
$\sharp L_{x}=1$ and $L_{x}$ is connected. 
\end{rem}
We now consider the first interaction cohomology groups of forward or 
backward self-similar systems. 
\begin{rem}
\label{intgroup0rem}
Let ${\frak L}:= (L, (h_{1},\ldots ,h_{m}))$ be a forward (resp. backward) 
self-similar system. Let 
$G=\langle h_{1},\ldots ,h_{m}\rangle $ and let $R$ be a $\ZZ $ module. 
If $\bigcap _{g\in G}g(L)\neq \emptyset $ (resp. if 
$\bigcap _{g\in G}g^{-1}(L)\neq \emptyset $),  
then, $\check{H}^{0}({\frak L};R)=R$ and $\check{H}^{p}({\frak L};R)=0$ for each 
$p\geq 1.$ In particular,  
if there exists a point $z\in L$ such that for each $j=1,\ldots ,m$, $h_{j}(z)=z$, 
then,  $\check{H}^{0}({\frak L};R)=R$ and $\check{H}^{p}({\frak L};R)=0$ for each 
$p\geq 1.$   
\end{rem}
By Remark~\ref{intgroup0rem}, 
we can find many examples of ${\frak L}$ such that 
$\check{H}^{p}({\frak L};R)=0$ for each $p\in \NN $ and each $\ZZ $ module $R$. 
\begin{rem}
\label{intgroupnon0rem}
For any $n\in \NN \cup \{ 0\} $, we also have many examples of forward or backward self-similar 
systems ${\frak L}=(L,(h_{1},\ldots ,h_{m}))$ such that 
for each field $R$, 
$0<\dim _{R}\check{H}^{n}({\frak L};R)<\infty .$ 
For example, let $M_{0}$ and $M_{1}$ be two cubes in $\RR ^{n+1}$ such that 
$M_{1}\subset \mbox{int}(M_{0}).$ Let $L:=M_{0}\setminus \mbox{int}(M_{1}).$ 
Then, we easily see that 
there exists a forward self-similar system ${\frak L}=(L,(h_{1},\ldots ,h_{m}))$ 
such that for each $j=1,\ldots ,m$, $h_{j}:L\rightarrow L$ is an injective contraction.  
For this ${\frak L}$, we have 
$\check{H}^{n}({\frak L};R)\cong \check{H}^{n}(L;R)=R.$ 
\end{rem}
We give a sufficient condition for the rank of the 
first interaction cohomology group of a system to be infinite.   
\begin{thm}
\label{preicmainthm}
Let $\frak{L}=(L,(h_{1},\ldots ,h_{m}))$ be a backward self-similar
 system. 
 Let $R$ be a field. 
We assume all of the following:
\begin{enumerate}
\item \label{preicmainthmc1}
$|N_{1}|$ is connected.  (See Remark~\ref{r:n1con}.) 

\item \label{preicmainthmc2}  
$(h_{1}^{2})^{-1}(L)\cap (\bigcup _{i: i\neq 1}h_{i}^{-1}(L))
=\emptyset . $
\item \label{preicmainthmc3}
There exist mutually distinct elements 
$j_{1}, j_{2}, j_{3}\in \{ 1,\ldots ,m\} $ such that 
$j_{1}=1$ and  such that for each $k=1,2,3,$ 
$h_{j_{k}}^{-1}(L)\cap h_{j_{k+1}}^{-1}(L)
\neq \emptyset ,$ where $j_{4}:=j_{1}.$ 
\item \label{preicmainthmc4}
For each $s,t\in \{ 1,\ldots ,m\} $, 
we have the following: if $s,t,1$ are mutually distinct, 
then $h_{1}^{-1}(L)\cap h_{s}^{-1}(L)\cap h_{t}^{-1}(L)=\emptyset .$  
\end{enumerate}  
Then, $\dim _{R}\check{H}^{1}({\frak L};R)=
\infty .$

\end{thm}
\begin{thm}
\label{sspreicmainthm}
Let $\frak{L}=(L,(h_{1},\ldots ,h_{m}))$ be a forward self-similar
 system such that for each $j=1,\ldots ,m$, 
$h_{j}:L\rightarrow L$ is injective.  
 Let $R$ be a field. 
We assume all of the following:
\begin{enumerate}
\item \label{sspreicmainthmc1}
$|N_{1}|$ is connected.  (See Remark~\ref{r:n1con}.) 

\item \label{sspreicmainthmc2}  
$h_{1}^{2}(L)\cap (\bigcup _{i: i\neq 1}h_{i}(L))
=\emptyset . $
\item \label{sspreicmainthmc3}
There exist mutually distinct elements 
$j_{1}, j_{2}, j_{3}\in \{ 1,\ldots ,m\} $ such that 
$j_{1}=1$ and  such that for each $k=1,2,3,$ 
$h_{j_{k}}(L)\cap h_{j_{k+1}}(L)
\neq \emptyset ,$ where $j_{4}:=j_{1}.$ 
\item \label{sspreicmainthmc4}
For each $s,t\in \{ 1,\ldots ,m\} $, 
we have the following: if $s,t,1$ are mutually distinct, 
then $h_{1}(L)\cap h_{s}(L)\cap h_{t}(L)=\emptyset .$  
\end{enumerate}  
Then, $\dim _{R}\check{H}^{1}({\frak L};R)=
\infty .$
\end{thm}
\begin{cor}
\label{c:sspreic}
Let ${\frak L}=(L,(h_{1},\ldots ,h_{m}))$ be a forward self-similar system
 such that each $h_{j}:L\rightarrow L$ is injective and such that 
 for each $x\in \Sigma _{m}$, $L_{x}$ is connected. 
Let $R$ be a field. Suppose that the conditions 1,2,3,4 in the assumptions 
of Theorem~\ref{sspreicmainthm} hold. Then, 
$\dim _{R}\check{H}^{1}({\frak L};R)=\dim _{R}\Psi (\check{H}^{1}({\frak L};R))=\infty $, 
$\dim _{R}\check{H}^{1}(L;R)=\infty $,  
and $\Psi :\check{H}^{1}({\frak L};R)\rightarrow \check{H}^{1}(L;R)$ is a monomorphism.     
\end{cor}
\begin{rem}
\label{r:rkfl} 
Let $K$ be a non-empty connected compact metric space and let 
$h_{j}: K\rightarrow K$ be a continuous map for each $j=1,\ldots , m.$ 
Let $L=R_{K,f}(h_{1},\ldots ,h_{m})$ and let ${\frak L}=(L,(h_{1},\ldots ,h_{m})).$ 
Regarding the forward self-similar system ${\frak L}$ (cf. Lemma~\ref{Rfsslem}), 
suppose that $|N_{1}|$ is connected. Then, $L$ is connected and $L_{x}$ is connected for each $x\in \Sigma _{m}.$ 
For, by Lemma~\ref{1conall}, which will be proved later, $|N_{k}|$ is connected for all $k\in \NN $, 
therefore $\bigcup _{|w|=k}h_{\overline{w}}(K)$ is connected for each $k\in \NN $. 
It implies that $L=\bigcap _{k=1}^{\infty }\bigcup _{|w|=k}h_{\overline{w}}(K)$ is connected.  
\end{rem}
\begin{ex}
\label{ex:csspreic}
Let ${\frak L}=(L,(h_{1},h_{2},h_{3}))$ be the forward self-similar system 
in Example~\ref{ex:SG1}. Then ${\frak L}$ satisfies the assumptions of Corollary~\ref{c:sspreic}. 
\end{ex}
\subsection{Application to the dynamics of polynomial semigroups}
\label{Application}
\begin{figure}[htbp]
\caption{The Julia set of a $3$-generator 
polynomial semigroup $G\in {\cal G}_{dis}$ with   
$\sharp (\mbox{Con}(J(G)))=\aleph _{0}.$}    
\ \ \ \ \ \ \ \ \ \ \ \ \ \ \ \ \ \ \ \ \ \ \ \ \ \ \ \ \ \ \ \ 
\ \ \ \  \ \ \ \ \ \ \ \ \ \ \ \ \ \ \ \ \ \ 
\includegraphics[width=2.5cm,width=2.5cm]{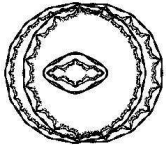}
\label{fig:3mapcountjulia2}
\end{figure}

In this subsection, we present some results on the Julia sets of postcritically bounded  polynomial semigroups $G$,  
which are obtained by applying the results in section~\ref{General}.  
The proofs of the results are given in section~\ref{Proofsappl}. 
\begin{df}
For each polynomial map $g:\CCI \rightarrow \CCI $, 
we denote by $CV(g)$ the set of all critical values of the holomorphic map 
$g:\CCI \rightarrow \CCI .$ 
Moreover, 
for a polynomial semigroup $G$, 
we set 
$ P(G)=
\overline{\bigcup _{g\in G}CV(g) } \ (\subset \CCI ).$   
The set $P(G)$ is called the {\bf postcritical set} of $G.$
Moreover, we set $P^{\ast }(G):= P(G)\setminus \{ \infty \} .$ 
The set $P^{\ast }(G)$  is called the {\bf planar postcritical set} of $G.$ 
We say that a polynomial semigroup $G$ is {\bf postcritically bounded} if 
$P^{\ast }(G)$ is bounded in $\CC .$  
\end{df}
\begin{df}
We denote by ${\cal G} $ the set of all postcritically bounded polynomial semigroups 
$G$ such that 
 for each $g\in G$,  $\deg (g)\geq 2.$   
Moreover, we set 
${\cal G}_{con}:=\{ G\in {\cal G}
\mid J(G)\mbox{ is connected}\} $ 
and ${\cal G}_{dis}:=
\{ G\in {\cal G}\mid J(G)
\mbox{ is disconnected}\} .$
\end{df}
\begin{rem}
\label{pgrem}
Let $G=\langle h_{1},\ldots ,h_{m}\rangle $ be a finitely generated 
polynomial semigroup. Then, 
$P(G)=\overline{\bigcup _{g\in G\cup \{ Id\} }g(\bigcup _{j=1}^{m}CV(h_{j}))}$ and 
$g(P(G))\subset P(G)$ for each $g\in G.$ From the above formula, one may use a computer 
to see if $G\in {\cal G}$ much in the same way as one verifies the boundedness 
of the critical orbit for the maps $f_{c}(z)=z^{2}+c.$  
\end{rem} 
\begin{df}
We set 
Rat$:=\{ g:\CCI \rightarrow \CCI \mid 
g \mbox{ is a non-constant rational map}\} $ endowed with the 
topology induced by the uniform convergence on $\CCI .$ Moreover, we set  
${\cal Y}:= \{ g: \CCI \rightarrow \CCI \mid g \mbox{ is a polynomial}, 
\deg (g)\geq 2\} $ endowed with the relative topology from Rat.  
Moreover, for each $m\in \NN $ we set 
${\cal Y}^{m}_{b}:= \{ (h_{1},\ldots ,h_{m})\in {\cal Y}^{m}\mid 
\langle h_{1},\ldots ,h_{m}\rangle \in {\cal G}\} .$ 
Furthermore, we set 
${\cal Y}^{m}_{b,con}:= \{ 
(h_{1},\ldots ,h_{m})\in {\cal Y}^{m}\mid 
\langle h_{1},\ldots ,h_{m}\rangle \in {\cal G}_{con}\} $ 
and ${\cal Y}^{m}_{b,dis}:=\{ (h_{1},\ldots ,h_{m})\in {\cal Y}^{m}
\mid \langle h_1,\ldots ,h_{m}\rangle \in {\cal G}_{dis}\} .$ 
\end{df}
\begin{rem}
It is well-known that for a polynomial $g\in {\cal Y}$, 
the semigroup $\langle g\rangle $ belongs to ${\cal G}$ if and only if 
$J(g)$ is connected (\cite{Mi}). However, for a general polynomial semigroup $G$, 
it is not true. For example, $\langle z^{3},z^{2}/4\rangle $ belongs to ${\cal G}_{dis}.$ 
There are many new phenomena about the dynamics of $G\in {\cal G}_{dis}$ which cannot hold 
in the dynamics of a single polynomial map. For the dynamics of $G\in {\cal G}_{dis}$, 
see \cite{SdpbpI, SdpbpII, SdpbpIII, S11, S10}. 
\end{rem}
We now present the first main result of this subsection. 
\begin{thm}
\label{compojth}
Let $G=\langle h_{1},\ldots ,h_{m}\rangle  
\in {\cal G}.$ 
 Then, for the backward self-similar 
 system ${\frak L}=(J(G),\ (h_{1},\ldots ,h_{m}))$,  
all of the statements 1,...,7 
in Theorem~\ref{precompojth} hold. 
\end{thm}
\begin{rem}
It is well known that if $G$ is a semigroup generated by a single $h\in $ Rat with $\deg (h)\geq 2$ 
or if $G$ is a non-elementary Kleinian group, then either $J(G)$ is connected or $J(G)$ has 
uncountably many connected components (\cite{B, Mi}). However, even for a finitely generated 
polynomial semigroup in ${\cal G}$, this is not true any more. In fact, 
in \cite{SdpbpI}, it was shown that for any positive integer $n$, there exists an  
element $(h_{1},\ldots ,h_{2n})\in {\cal Y}_{b}^{2n}$ such that 
$\sharp \mbox{Con}(J(\langle h_{1},\ldots ,h_{2n}\rangle ))=n.$ 
Moreover, in \cite{SdpbpI}, it was shown that 
there exists an element $(h_{1},h_{2},h_{3})\in {\cal Y}_{b}^{3}$ such that 
$\sharp \mbox{Con}(J(\langle h_{1},h_{2},h_{3}\rangle ))=\aleph _{0}$ (see Figure~\ref{fig:3mapcountjulia2}).    
\end{rem}
 \begin{figure}[htbp]
\caption{The Julia set of $G=\langle g_{1}^{2},g_{2}^{2}\rangle $, 
where $g_{1}(z):=z^{2}-1,\ g_{2}(z):=\frac{z^{2}}{4}$. 
$G\in {\cal G}_{dis}$ and 
$\sharp (\mbox{Con}(J(G)))>\aleph _{0}.$ }
\ \ \ \ \ \ \ \ \ \ \ \ \ \ \ \ \ \ \ \ \ \ \ \ \ \ \ \ \ \ \ \ \ \ \ \ \ \ \ \ \ \ \ 
\ \ \ \ \ \ \ \ \ \ \ \ \ 
\includegraphics[width=2.2cm,width=2.2cm]{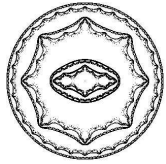}
\label{fig:dcgraph}
\end{figure}
By Remark~\ref{intgroup0rem}, for each $m\in \NN $, 
 there exists an element $(h_{1},\ldots ,h_{m})\in 
 {\cal Y}^{m}_{b}$ such that 
 setting $G=\langle h_{1},\ldots ,h_{m}\rangle$, 
we have $\check{H}^{1}(J(G),(h_{1},\ldots ,h_{m});R)=0.$ 
We will show that 
there exists an element $(h_{1},\ldots ,h_{4})\in {\cal Y}^{4}_{b}$
 such that setting $G=\langle h_{1},\ldots ,h_{4}\rangle $, 
 $\check{H}^{1}(J(G),(h_{1},\ldots ,h_{4});R)$ has 
 infinite rank.
\begin{thm}
\label{icmainthm}
Let $m\in \NN $ and 
let $(h_{1},\ldots ,h_{m})\in {\cal Y}^{m}_{b}.$ 
Let $G=\langle h_{1},\ldots ,h_{m}\rangle .$
 Let $R$ be a field. 
For the backward self-similar system  
${\frak L}=(J(G),\ (h_{1},\ldots ,h_{m}))$, 
suppose that all of the conditions 1, 2, 3, 4  
in the assumptions of Theorem~\ref{preicmainthm} hold.  
Then, we have the following statements 1, 2, 3.
\begin{enumerate}
\item $\dim _{R}\check{H}^{1}({\frak L}; R)=
\dim _{R}\Psi (\check{H}^{1}({\frak L};R))=\infty $ and 
$\dim _{R}\check{H}^{1}(J(G);R)=\infty .$ 
\item 
$\Psi : \check{H}^{1}({\frak L}; R)
\rightarrow \check{H}^{1}(J(G);R)$ 
is a monomorphism. 
\item $F(G)$ has infinitely many connected components.
\end{enumerate}

\end{thm}
\begin{prop}
\label{pbpch1inftyexprop}
There exists an element  
 $h=(h_{1},h_{2},h_{3},h_{4})\in {\cal Y}^{4}_{b}$ 
 which satisfies  the assumptions 
 of Theorem~\ref{icmainthm}. In particular,   
 for this $h$, setting $G=\langle h_{1},\ldots ,h_{4}\rangle$, 
 for any field $R$,  we have that  
 $\dim _{R}\check{H}^{1}(J(G),(h_{1},\ldots ,h_{4});R)=
 \dim _{R}\Psi (\check{H}^{1}(J(G),(h_{1},\ldots ,h_{4});R))=\infty $ and 
 $F(G)$ has infinitely many connected components (see Figure~\ref{fig:h1infty}).  
\end{prop}
\begin{prob}[Open]
Let $m\in \NN $ with $m\geq 2.$ 
Are there any $(h_{1},\ldots ,h_{m})\in {\cal Y}^{m}_{b}$ such that 
$$0<\dim _{R}\check{H}^{1}(J(\langle h_{1},\ldots ,h_{m}\rangle ), (h_{1},\ldots ,h_{m});R)<\infty \ ?$$
\end{prob}
\subsection{Postunbranched systems}
\label{Postunbranched}
In this subsection, we introduce ``postunbranched systems,'' and we present 
some results on the interaction (co)homology groups of such systems. 
The proofs of the main results are given in 
section~\ref{ProofsPost}. 
\begin{df}
\label{d:pu}
Let ${\frak L}=(L,(h_{1},\ldots ,h_{m}))$ be a forward (resp. backward)  
self-similar system. 
For each $(i,j)\in \{ 1,\ldots ,m\} ^{2}$ with 
$i\neq j$,  
we set $C_{i,j}=C_{i,j}({\frak L}):=h_{i}(L)\cap h_{j}(L)$ (resp. 
$C_{i,j}=C_{i,j}({\frak L}):=h_{i}^{-1}(L)\cap h_{j}^{-1}(L)$).  
We say that ${\frak L}$ 
is 
{\bf postunbranched }if for any 
$(i,j)\in \{ 1,\ldots ,m\} ^{2}$ such that $i\neq j$ and $C_{i,j}\neq \emptyset$, there exists a 
unique $x=x(i,j)\in \Sigma _{m}$ 
such that 
\begin{itemize}
\item 
$h_{i}^{-1}(C_{i,j})\subset L_{x}$ (resp. 
$h_{i}(C_{i,j})\subset L_{x}$) and 
\item 
for each $x'\in \Sigma _{m}$ with $x'\neq x$, we have 
$h_{i}^{-1}(C_{i,j})\cap L_{x'}=\emptyset $ (resp.  
$h_{i}(C_{i,j})\cap L_{x'}=\emptyset $ ).  
\end{itemize}
\end{df}
The following Lemmas~\ref{pssubsyslem}, \ref{puiteratelem}, \ref{psiterateblem}, \ref{psiterateflem}  
are easy to show from the definition above. 
\begin{lem}
\label{pssubsyslem}
Let ${\frak L}=(L,(h_{1},\ldots ,h_{m}))$ be a forward or backward self-similar 
system. Suppose that ${\frak L}$ is postunbranched. 
Then, any subsystem ${\frak M}$ of ${\frak L}$ is postunbranched. 
\end{lem}
\begin{lem}
\label{puiteratelem}
Let ${\frak L}=(L,(h_{1},\ldots ,h_{m}))$ be a 
forward or backward self-similar system. 
Suppose that ${\frak L}$ is postunbranched. 
When ${\frak L}$ is a forward self-similar system, 
we assume further that for each $j=1,\ldots ,m$, 
$h_{j}:L\rightarrow L$ is injective. 
Then, for each $n\in \NN $, an $n$-th iterate of 
${\frak L}$ is postunbranched.  
\end{lem}
\noindent {\bf Notation:} 
Let $m\in \NN .$  
For each $j=1,\ldots ,m$, we set 
$(j)^{\infty }:= (j,j,\ldots )\in \Sigma _{m}.$  
\begin{lem}
\label{psiterateblem}
Let ${\frak L}=(L,(h_{1},\ldots ,h_{m}))$ be a backward 
self-similar system. 
Suppose that for each $(i,j)\in \{ 1,\ldots ,m\} ^{2}$ such that $i\neq j$ and $C_{i,j}\neq \emptyset$, there exists an $r\in \{ 1,\ldots ,m\} $ 
such that 
$h_{i}(C_{i,j})\subset L_{(r)^{\infty }}$ 
and $L_{(r)^{\infty }}\subset (L\setminus  \bigcup _{k:k\neq r}h_{k}^{-1}(L)).$ 
Then, for any $n\in \NN $, an $n$-th iterate of ${\frak L}$ is postunbranched.  
\end{lem}
\begin{lem}
\label{psiterateflem}
Let ${\frak L}=(L,(h_{1},\ldots ,h_{m}))$ be a forward 
self-similar system such that for each $j=1,\ldots ,m$, 
$h_{j}:L\rightarrow L$ is injective.  
Suppose that for each $(i,j)\in \{ 1,\ldots ,m\} ^{2}$ such that $i\neq j$ and $C_{i,j}\neq \emptyset$, there exists an 
$r\in \{ 1,\ldots, m\} $ such that 
$h_{i}^{-1}(C_{i,j})\subset L_{(r)^{\infty }}$ and 
$L_{(r)^{\infty }}\subset (L\setminus \bigcup _{k:k\neq r}h_{k}(L)).$ 
Then, for any $n\in \NN $, an $n$-th iterate of ${\frak L}$ is postunbranched.  
\end{lem}
From Lemmas~\ref{pssubsyslem}, \ref{puiteratelem}, \ref{psiterateblem}, \ref{psiterateflem}, 
we can easily obtain many examples of postunbranched systems.\\ 
\noindent {\bf Notation:} 
We denote by Fix$(f)$ the set of all fixed points of $f.$ 
\begin{ex}[{\bf Sierpi\'{n}ski gasket}]
\label{SGex}
Let ${\frak L}=(L,(h_{1},h_{2} ,h_{3}))$ be the forward self-similar system in Example~\ref{ex:SG1}.  
Thus $L$ is the Sierpi\'{n}ski gasket. 
(See Figure~\ref{fig:sgasket}.)    
From Figure~\ref{fig:sgasket}, we see that for each $(i,j)\in \{ 1,2,3\}^{2}$ such that $i\neq j$ and 
$C_{i,j}\neq \emptyset $, 
$h_{i}^{-1}(C_{i,j})=\mbox{Fix}(h_{j})\cap L=L_{(j)^{\infty }}
\subset (L\setminus 
\bigcup _{k:k\neq j}h_{k}(L)).$ 
From Lemmas~\ref{psiterateflem} and \ref{pssubsyslem}, it follows that 
for any $n\in \NN $, if ${\frak M}=(M,(g_{1},\ldots ,g_{t}))$ is a subsystem of an 
$n$-th iterate of ${\frak L}$, then ${\frak M}$ is postunbranched.  
\end{ex}
\begin{ex}[{\bf Pentakun, Snowflake}]
\label{PSex}
Let ${\frak L}=(L,(h_{1},\ldots ,h_{m}))$ be a forward self-similar system 
in \cite[Example 3.8.11 (Pentakun)]{K} or \cite[Example 3.8.12 (Snowflake)]{K}. 
Hence $L$ is one of the snowflake, the pentakun, the heptakun, the octakun, and so on. 
(The definition of the snowflake is as follows: let $p_{k}=\exp (2k\pi \sqrt{-1}/6)$ for each $k=1,\ldots, 6$ and 
let $p_{7}=0.$ We define $h_{k}:\CC \rightarrow \CC $ by 
$h_{k}(z)=(z-p_{k})/3+p_{k}$ for each $k=1,\ldots ,7.$ The snowflake is 
$M_{\CC }(h_{1},\ldots ,h_{7}).$ The definition of the pentakun is as follows: 
for each $k=1,\ldots ,5$, let $q_{k}=\exp (2k\pi \sqrt{-1}/5).$ We define $g_{k}:\CC \rightarrow \CC $ 
by $g_{k}(z)=\frac{3-\sqrt{5}}{2}(z-p_{k})+p_{k}$ for each $k=1, \ldots ,5.$ 
The pentakun is $M_{\CC }(g_{1},\ldots ,g_{5}).$)  
Then, looking at Figure~\ref{fig:snowflake2}, 
it is easy to see that for each $(i,j)\in \{ 1,\ldots ,m\}^{2}$ such that $i\neq j$ and 
$C_{i,j}\neq \emptyset $, there exists an $r\in \{ 1,\ldots ,m\} $ such that 
$h_{i}^{-1}(C_{i,j})=\mbox{Fix}(h_{r})\cap L=L_{(r)^{\infty }}
\subset (L\setminus 
\bigcup _{k:k\neq r}h_{k}(L)).$ 
From Lemmas~\ref{psiterateflem} and \ref{pssubsyslem}, it follows that 
for any $n\in \NN $, if ${\frak M}=(M,(g_{1},\ldots ,g_{t}))$ is a subsystem of an 
$n$-th iterate of ${\frak L}$, then ${\frak M}$ is postunbranched. 
\begin{figure}[htbp]
\caption{(From left to right) Snowflake, Pentakun}    
\ \ \ \ \ \ \ \ \ \ \ \ \ \ \ \ \ \ \ \ \ \ \ \ \ \ \ \ \ \ \ \ 
\ \ \ \  \ \ \ \ \ \ \ \ \ \ \ \ \
\includegraphics[width=1.8cm,width=1.8cm]{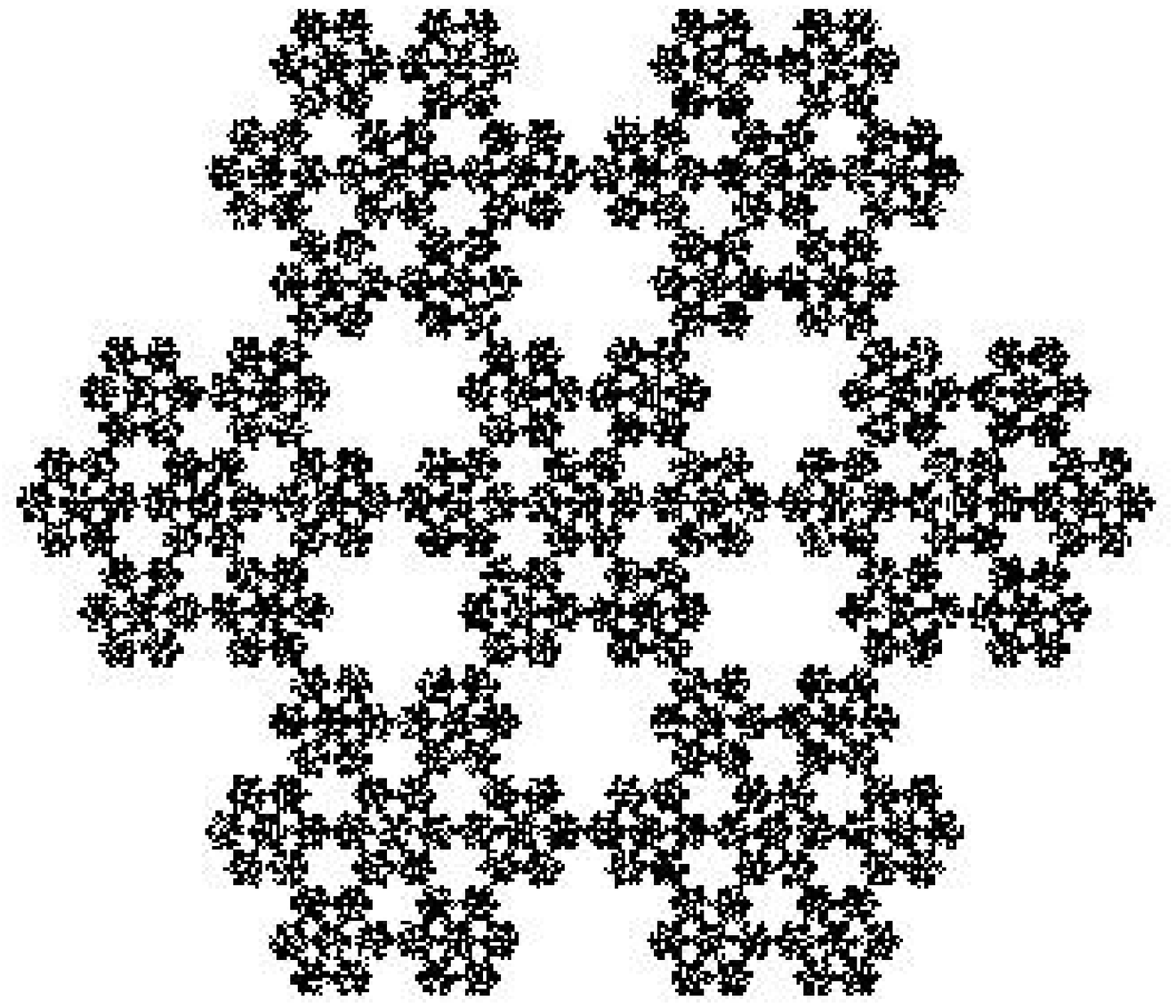}
\includegraphics[width=2.3cm, width=2.3cm,angle=17.5,totalheight=2.7cm]{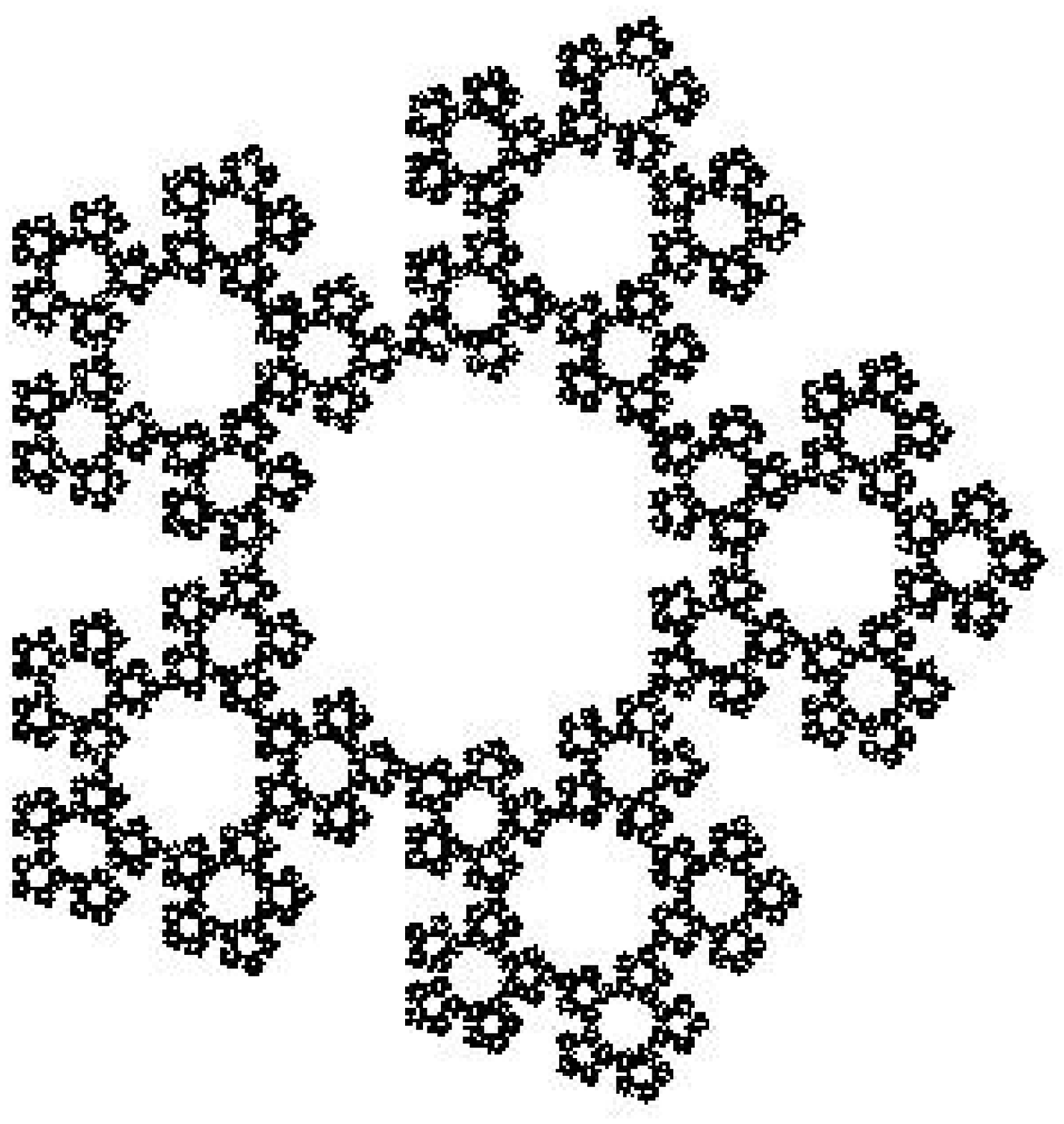}
\label{fig:snowflake2}
\end{figure}

\end{ex} 
In order to state the main results, we need some definitions. 
\begin{df}
Let ${\frak L}=(L,(h_{1},\ldots ,h_{m}))$ be a forward or backward self-similar system and let $R$ be a $\ZZ $ module. 
Let $w\in \Sigma _{m}^{\ast }$ with $|w|=l.$ 
Let $k\in \NN $ with $k> l.$ 
We denote by $N_{k,w}$ (or $N_{k,w}({\frak L})$) the 
unique full subcomplex of $N_{k}$ whose vertex set is equal to 
$\{ wx\mid x\in \{ 1,\ldots ,m\} ^{k-l}\} .$ 
Moreover, for each $j=1,\ldots ,m$, we set 
$N_{1,j}:= \{ j\}\ (\subset N_{1}).$ 
We denote by $w_{\ast }: N_{k}\rightarrow N_{k+l}$ the simplicial map 
assigning to each vertex $x=(x_{1},\ldots ,x_{k})\in \{ 1,\ldots ,m\} ^{k}$ 
the vertex $wx\in \{ 1,\ldots ,m\} ^{k+l}.$ 
We denote by $w_{\ast }: H_{p}(N_{k};R)\rightarrow H_{p}(N_{k+l,w};R)$
the homomorphism induced by the above simplicial map $ w_{\ast }: N_{k}\rightarrow N_{k+l,w}.$ 
Moreover, we denote by $w^{\ast }: H^{p}(N_{k+l,w};R)\rightarrow H^{p}(N_{k};R)$ the 
homomorphism induced by $w_{\ast }: N_{k}\rightarrow N_{k+l,w}.$ 
Moreover, we denote by  $\bigoplus _{j=1}^{m}(j)_{\ast }: \bigoplus _{j=1}^{m}H_{p}(N_{k};R)\rightarrow  \bigoplus _{j=1}^{m}H_{p}(N_{k+1,j};R)\cong 
H_{p}(\bigcup _{j=1}^{m}N_{k+1,j};R)$ the homomorphism $(\alpha _{j})_{j=1}^{m}\mapsto (j_{\ast }(\alpha _{j}))_{j=1}^{m} $. 
Moreover, let $\iota : \bigcup _{j=1}^{m}N_{k+1,j}\rightarrow N_{k+1}$ be the canonical embedding and 
let $(\eta _{k})_{\ast }: \bigoplus _{j=1}^{m}H_{p}(N_{k};R)\rightarrow H_{p}(N_{k+1};R)$ be the composition 
$\iota _{\ast }\circ (\bigoplus _{j=1}^{m}(j)_{\ast }).$ 
Similarly, we denote by $\bigoplus _{j=1}^{m}(j)^{\ast }: H^{p}(\bigcup _{j=1}^{m}N_{k+1,j};R)\cong 
\bigoplus _{j=1}^{m}H^{p}(N_{k+1,j};R)\rightarrow \bigoplus _{j=1}^{m}H^{p}(N_{k};R)$ the 
homomorphism $(\beta _{j})_{j=1}^{m}\mapsto (j^{\ast }\beta _{j})_{j=1}^{m}.$ 
Let $\eta _{k}^{\ast }: H^{p}(N_{k+1};R)\rightarrow \bigoplus _{j=1}^{m}H^{p}(N_{k};R)$ be the 
composition $(\bigoplus _{j=1}^{m}(j)^{\ast })\circ \iota ^{\ast }.$ 
\end{df}
From this definition, it is easy to see that the following lemma holds. 
\begin{lem}
\label{nkwlem}
Let ${\frak L}=(L,(h_{1},\ldots ,h_{m}))$ be a forward or backward self-similar system. 
When ${\frak L}$ is a forward self-similar system, we assume further that 
$h_{j}: L\rightarrow L$ is injective for each $j.$  
Let $w\in \Sigma _{m}^{\ast }$ with $|w|=l.$ 
Then, for each $k\in \NN $, 
the simplicial map 
$w_{\ast }: N_{k}\rightarrow N_{k+l,w}$ is isomorphic. 
\end{lem}
\begin{df}
Let ${\frak L}:=(L,(h_{1},\ldots ,h_{m}))$ be a forward or backward self-similar system and 
 let $R$ be a $\ZZ $ module.   
Let $w\in \Sigma _{m}^{\ast }$ with $|w|=l$ and let $k\in \NN .$  
We denote by $w_{\ast }: N_{k}\rightarrow N_{k+l}$ 
the simplicial map assigning to each vertex  
$x=(x_{1},\ldots ,x_{k})\in \{ 1,\ldots ,m\} ^{k}$ the 
vertex $wx\in \{ 1,\ldots ,m\} ^{k+l}.$ 
We denote by $w_{\ast }: H_{\ast }({\frak L};R)_{k}
\rightarrow H_{\ast }({\frak L};R)_{k+l}$ the homomorphism induced by 
the above simplicial map $w _{\ast }: N_{k}\rightarrow 
N_{k+l}.$ Moreover, 
we denote by $w^{\ast }:H^{\ast }({\frak L};R)_{k+l}\rightarrow 
H^{\ast }({\frak L};R)_{k}$ the homomorphism induced by 
$w_{\ast }: N_{k}\rightarrow 
N_{k+l}.$ Moreover,  
we denote by $q_{w}:N_{1}\rightarrow N_{l}$ the 
constant simplicial map assigning to each vertex 
$x\in \{ 1,\ldots ,m\} $ the vertex $w.$  
\end{df}
From the above definition, it is easy to see that the following lemma holds.
\begin{lem}
\label{wmaplem}
Let ${\frak L}=(L,(h_{1},\ldots ,h_{m}))$ be a forward or backward self-similar system. 
Then, for each $k\in \NN $ with $k\geq 2$, we have 
$\varphi _{k}j_{\ast }(x)=j_{\ast }\varphi _{k-1}(x)$ for each 
$x\in N_{k}$, and 
$\varphi _{1}j_{\ast }(x)=q_{j}(x)$ for each $x\in N_{1}.$ 
More generally, 
let $w\in \Sigma _{m}^{\ast }$ with $|w|=l. $
Then, for each $k\in \NN $ with $k\geq 2$, we have 
$\varphi _{l+k-1}w_{\ast }(x)=w_{\ast }\varphi _{k-1}(x)$ for each 
$x\in N_{k}$, and 
$\varphi _{l}w_{\ast }(x)=q_{w}(x)$ for each $x\in N_{1}.$ 

\end{lem}
\begin{df}
Let ${\frak L}=(L,(h_{1},\ldots ,h_{m}))$ be a forward or backward self-similar system and 
let $R$ be a $\ZZ $ module.  
Let $w=(w_{1},\ldots ,w_{l})\in \Sigma _{m}^{\ast }$ with $|w|=l. $ 
We define a homomorphism $w_{\ast }: \check{H}_{p}({\frak L};R)\rightarrow 
\check{H}_{p}({\frak L};R)$ as follows. 
Let $a=(a_{k})\in \check{H}_{p}({\frak L};R)= 
\varprojlim _{k}H_{p}(N_{k};R)$ be an element, where 
for each $k\in \NN $, $a_{k}\in H_{p}(N_{k};R)$ and 
$(\varphi _{k})_{\ast }(a_{k+1})=a_{k}.$ 
For each $k\in \NN $, we set $b_{k+l}:= w_{\ast }(a_{k})\in H_{p}(N_{k+l};R).$ 
Moreover, for each $s\in \NN $ with $1\leq s\leq l$, we set 
$b_{s}:= (q_{w|s})_{\ast }(a_{1})\in H_{p}(N_{s};R).$ 
Then, by Lemma~\ref{wmaplem}, 
$b=(b_{t})_{t=1}^{\infty }$ determines an element in 
$\check{H}_{p}({\frak L};R)= \varprojlim _{k}\check{H}_{p}(N_{k};R).$ 
We set $w_{\ast }(a):= b.$ 

 Similarly, we define 
 a homomorphism 
 $w^{\ast }:\check{H}^{p}({\frak L};R)\rightarrow \check{H}^{p}({\frak L};R)$ 
 as follows. 
Let $a\in \check{H}^{p}({\frak L};R)= \varinjlim _{k} H^{p}(N_{k};R)$ 
be an element. When $a$ is represented by an element 
$c\in H^{p}(N_{k};R)$ with $k\geq l+1$,  
 we set $c_{1}:= w^{\ast }(c)\in H^{p}(N_{k-l};R)$ and 
let $w^{\ast }(a):= \mu _{k-l,p}(c_{1})\in \check{H}^{p}({\frak L};R).$ 
When $a$ is represented by an element 
$c\in H^{p}(N_{k};R)$ with $k\leq l$, 
we set $c_{1}:=q_{w|k}^{\ast }(c)\in H^{p}(N_{1};R)$ and 
let $w^{\ast }(a)=\mu _{1,p}(c_{1})\in \check{H}^{p}({\frak L};R).$ 
By Lemma~\ref{wmaplem}, $w^{\ast }(a)\in \check{H}^{p}({\frak L};R)$ is 
well defined and independent of the choice of $c.$ 

 Furthermore, we define a homomorphism 
$\theta : \check{H}^{p}(\frak{L};R)\rightarrow \bigoplus _{j=1}^{m}\check{H}^{p}(\frak{L};R)$ 
by $ \theta (c):= (j^{\ast }(c))_{j=1}^{m}.$   
\end{df}
\begin{df}
\label{d:cohcomp}
Let ${\frak L}=(L,(h_{1},\ldots ,h_{m}))$ be a forward or backward self-similar
 system.  
 Let $R$ be a field and let $T$ be a $\ZZ $ module. 
 Let $a_{r,k}=a_{r,k}({\frak L}; R):=\dim _{R}\check{H}^{r}({\frak L}; R)_{k}$ 
 for each $r,k\in \ZZ $ with $r\geq 0,k\geq 1.$ 
Moreover, we set 
$u^{r}({\frak L};R):= \limsup _{k\rightarrow \infty }\frac{1}{k}\log a_{r,k}\in \{ -\infty \} \cup [0 ,+\infty)$ and 
$l^{r}({\frak L}; R):= \liminf _{k\rightarrow \infty }\frac{1}{k}\log a_{r,k}\in \{ -\infty \} \cup [0,+\infty ).$ 
The quantity $u^{r}({\frak L};R)$ is called the $r$-th {\bf upper cohomological complexity} of 
${\frak L}$ with coefficients $R$, and $l^{r}({\frak L};R)$ is called the $r$-th 
{\bf lower cohomological complexity} of ${\frak L}$ with coefficients $R.$  
 Moreover, let 
 $a_{r,\infty }=a_{r,\infty }({\frak L}; R):=\dim _{R}\check{H}^{r}({\frak L};R)$ and 
 $b_{1,\infty }=b_{1,\infty }({\frak L}; R):= \dim _{R}\mbox{Im}\mu _{1,1}.$ 
Moreover, let $S_{1}=S_{1}({\frak L})$ be the CW complex defined by $S_{1}:=|N_{1}|/\{ 1,\ldots ,m\} .$
Moreover, for each $k\in \NN $ with $k>1$, we set 
$A_{k}=A_{k}({\frak L}; T):= 
\mbox{Im} ((\varphi _{k,1})_{\ast }: \check{H}_{1}({\frak L};T)_{k}
\rightarrow \check{H}_{1}({\frak L};T)_{1})$, 
$B_{k}=B_{k}({\frak L}; R):= \mbox{Im} (\varphi _{k,1}^{\ast }: \check{H}^{1}({\frak L};R)_{1}
\rightarrow \check{H}^{1}({\frak L};R)_{k})$, and 
$\lambda _{k}=\lambda _{k}({\frak L};R):= \dim _{R}B_{k}.$  
\end{df}
\begin{rem}
\label{r:luinv}
From the above notation, we have $0\leq a_{r,k}\leq \frac{m^{k}(m^{k}-1)\cdots (m^{k}-r)}{(r+1)!}$ and  
$-\infty \leq l^{r}({\frak L};R)\leq u^{r}({\frak L};R)\leq (r+1)\log m.$ 
Moreover, by Remark~\ref{r:functor}, 
it follows that if ${\frak L}_{1}\cong {\frak L}_{2}$, then 
$a_{r,k}({\frak L}_{1};R)=a_{r,k}({\frak L}_{2};R)$, 
$a_{r,\infty }({\frak L}_{1};R)=a_{r,\infty }({\frak L}_{2};R)$,
$b_{1,\infty }({\frak L}_{1};R)=b_{1,\infty }({\frak L}_{2};R)$, 
$l^{r}({\frak L}_{1};R)=l^{r}({\frak L}_{2};R)$, 
$u^{r}({\frak L}_{1};R)=u^{r}({\frak L}_{2};R)$, 
$A_{k}({\frak L}_{1};T)\cong A_{k}({\frak L}_{2};T)$, 
$B_{k}({\frak L}_{1};R)\cong B_{k}({\frak L}_{2};R)$, 
and $\lambda _{k}({\frak L}_{1};R) =\lambda _{k}({\frak L}_{2};R).$ 
\end{rem}

We now state one of the main results on the interaction (co)homology groups of postunbranched systems.
\begin{thm}
\label{puthm}
Let ${\frak L}=(L,(h_{1},\ldots ,h_{m}))$ be a forward or backward self-similar
 system.  
When ${\frak L}=(L,(h_{1},\ldots ,h_{m}))$ is a  
 forward self-similar system, we assume further that   
$h_{j}:L\rightarrow L$ is 
 injective for each $j=1,\ldots ,m$. 
Furthermore, let $R$ be a field and let $T$ be a $\ZZ $ module.  
Suppose that ${\frak L}$ 
 is postunbranched. Then,
  we have all of the following statements \ref{puthm1},...,\ref{puthm7}.
\begin{enumerate}
\item \label{puthm1}
Let $r\geq 2 $ and $k\geq 1$. Then, 
$a_{r,k+1}=ma_{r,k}+a_{r,1}$ and 
there exists an exact sequence: 
\begin{equation}\label{puthm1eq}
0\longrightarrow \bigoplus _{j=1}^{m}\check{H}_{r}(\frak{L};T)_{k}\overset{(\eta _{k})_{\ast }}{\longrightarrow }
\check{H}_{r}(\frak{L};T)_{k+1}\overset{(\varphi _{k+1,1})_{\ast }}\longrightarrow 
\check{H}_{r}(\frak{L};T)_{1}\longrightarrow 0.
\end{equation} 
\item \label{puthm2}
Let $r\geq 2$ and $k\geq 1.$ 
If $\check{H}_{r}(\frak{L};T)_{1}=0$, then 
$\check{H}_{r}(\frak{L};T)_{k}=\check{H}_{r}(\frak{L};T)=0.$ 
\item \label{puthmrgeq2}
Let $r\geq 2.$ Then, there exists an exact sequence of $R$ modules: 
\begin{equation}
\label{puthmrgeq2eq}
0\longrightarrow \check{H}^{r}(\frak{L};R)_{1}
\overset{\mu _{1,r} }{\longrightarrow }\check{H}^{r}(\frak{L};R)
\overset{\theta }{\longrightarrow }\bigoplus _{j=1}^{m}\check{H}^{r}(\frak{L};R)\longrightarrow 0. 
\end{equation}
\item \label{puthmmukrinj} 
Let $r\neq 1$ and $k\geq 1.$ Then, 
$\mu _{k,r}:\check{H}^{r}({\frak L};R)_{k}\rightarrow 
\check{H}({\frak L};R)$ and 
$\varphi _{k}^{\ast }: \check{H}^{r}({\frak L};R)_{k}\rightarrow 
\check{H}^{r}({\frak L};R)_{k+1}$ are monomorphisms. 
\item \label{puthmcohsupp}
Let $r\geq 2.$ 
\begin{enumerate}
\item If $\check{H}^{r}({\frak L};R)_{1}=0$, then for each $k\in \NN $, 
$\check{H}^{r}({\frak L};R)_{k}=0$ and $\check{H}^{r}({\frak L};R)=0.$
\item If $\check{H}^{r}({\frak L};R)_{1}\neq 0$, then $a_{r,\infty }=\infty .$ 
\end{enumerate}
\item \label{puthmr1fhom}
Let $k\in \NN $. Then we have the following exact sequences:
\begin{equation}
\label{puthmr1fhomeq1}
0\rightarrow \bigoplus _{j=1}^{m}\check{H}_{1}({\frak L};T)_{k}
\overset{(\eta _{k})_{\ast }}{\rightarrow }\check{H}_{1}({\frak L};T)_{k+1}
\overset{(\varphi _{k+1,1})_{\ast }}{\rightarrow }A_{k+1}\rightarrow 0
\end{equation}
and 
\begin{equation}
\label{puthmr1fhomeq2}
0\rightarrow A_{k+1}\rightarrow H_{1}(S_{1};T)\rightarrow 
\bigoplus _{j=1}^{m}\check{H}_{0}({\frak L};T)_{k}
\overset{(\eta _{k})_{\ast }}{\rightarrow } \check{H}_{0}({\frak L};T)_{k+1}\rightarrow 0.
\end{equation}
\item \label{puthmr1f}
Let $k\in \NN $. Then we have the following exact sequences of $R$ modules:
\begin{equation}
\label{puthmr1feq2}
0\rightarrow B_{k+1}\rightarrow \check{H}^{1}({\frak L};R)_{k+1}
\overset{\eta _{k}^{\ast }}{\rightarrow }\bigoplus _{j=1}^{m}\check{H}^{1}({\frak L};R)_{k}
\rightarrow 0
\end{equation}
and 
\begin{equation}
\label{puthmr1feq1}
0\rightarrow \check{H}^{0}({\frak L};R)_{k+1}\overset{\eta _{k}^{\ast }}{\rightarrow }
\bigoplus _{j=1}^{m}\check{H}^{0}({\frak L};R)_{k}\rightarrow H^{1}(S_{1};R)\rightarrow 
B_{k+1}\rightarrow 0.
\end{equation}

\item \label{puthmr1}
We have the following exact sequences of $R$ modules: 
\begin{equation}
\label{puthmr1spliteq2}
0\rightarrow \mbox{{\em Im}} \mu _{1,1}\rightarrow \check{H}^{1}({\frak L};R)
\overset{\theta }{\rightarrow } 
\bigoplus _{j=1}^{m}\check{H}^{1}({\frak L};R)\rightarrow 0
\end{equation}
and 
\begin{multline}
\label{puthmr1spliteq1}
0\longrightarrow \check{H}^{0}({\frak L};R)\overset{\theta }{\longrightarrow} \bigoplus _{j=1}^{m}\check{H}^{0}({\frak L};R)\longrightarrow  
H^{1}(S_{1} ;R)\longrightarrow \mbox{{\em Im}} \mu _{1,1}\longrightarrow 0.
\end{multline}

\item \label{puthmakvalues1} 
Let $k\in \NN .$ Then, we have that 
$a_{1,k+1}=ma_{1,k}+\lambda _{k+1}$
and 
$a_{0,k+1}=ma_{0,k}-m+a_{0,1}-a_{1,1}+\lambda _{k+1}.$  
\item \label{puthmlambda}
For each $k\in \NN $, $0\leq b_{1,\infty }\leq \lambda _{k+2}\leq \lambda _{k+1}
\leq \lambda _{2}\leq a_{1,1}.$ 
Moreover, there exists a positive integer $l$ such that 
for each $k\in \NN $ with $k\geq l$, $\lambda _{k}=b_{1,\infty }.$   
\item \label{puthm3}
For each $k\in \NN $, 
$a_{0,k+1}=ma_{0,k}-m+a_{0,1}-a_{1,1}-ma_{1,k}+a_{1,k+1}.$  
\item \label{puthmakvalues2}
For each $k\in \NN $, 
$ma_{1,k}\leq a_{1,k+1}\leq ma_{1,k}+a_{1,1}.$ 
\item \label{puthm4}
For each $k\in \NN $, $ma_{0,k}-m+a_{0,1}-a_{1,1}\leq 
a_{0,k+1}\leq ma_{0,k}-m+a_{0,1}.$ 
\item \label{puthmulrgeq1}
Let $r\geq 1 .$ Then, 
either (a) $l^{r}({\frak L};R)=u^{r}({\frak L};R)=-\infty $ or 
(b)$l^{r}({\frak L};R)=u^{r}({\frak L};R)=\log m.$
\item \label{puthmulr0}
Either (a) $l^{0}({\frak L};R)=u^{0}({\frak L};R)=0$ or 
(b) $l^{0}({\frak L};R)=u^{0}({\frak L};R)=\log m.$ 
\item \label{puthmavalues}
Let $r\geq 1 .$ Then, either $a_{r,\infty }=0$ or $a_{r,\infty }=\infty .$ 
\item \label{puthm8}
If $a_{0,\infty }<\infty  $, 
then $m-a_{0,1}+a_{1,1}=(m-1)a_{0,\infty }+b_{1,\infty }.$ 
\item \label{puthm9}
If  $m\geq 2$ and $\frac{m-a_{0,1}+a_{1,1}}{m-1}\not\in  
\NN $, then 
at least one of $a_{0,\infty }$ and $a_{1,\infty }$ is 
equal to $\infty .$
\item \label{puthm5}
If $m\geq 2$ and there exists an element 
$k_{0}\in \NN $ such that  
$a_{0,k_{0}}>\frac{1}{m-1}
(m-a_{0,1}+a_{1,1}), $ then 
$a_{0,k+1}>a_{0,k}$ for each $k\geq k_{0}.$ 
\item \label{puthm6}
If $m\geq 2$, then $a_{0,\infty }\in 
\{ x\in \NN \mid a_{0,1}\leq x\leq \frac{1}{m-1}(m-a_{0,1}+a_{1,1})\}  
\cup 
\{ \infty \} .$
\item \label{puthmless6}
If $2\leq m\leq 6$ and $|N_{1}|$ is disconnected, then 
$a_{0,\infty }=\infty $ and $L$ has infinitely many connected components. 
\item \label{puthmB2zero}
If $B_{2}=0$, then $\check{H}^{1}({\frak L};R)=0.$ 
\item \label{puthm7}
If $|N_{1}|$ is connected, then 
we have the following. 
 \begin{enumerate}
 \item 
 \label{puthm7a}
For each $k\in \NN $, we have the following exact sequence: 
\begin{equation}
\label{puthm7aeq}
0\longrightarrow \bigoplus _{j=1}^{m}\check{H}_{1}(\frak{L};T)_{k}\overset{(\eta _{k})_{\ast }}{\longrightarrow }
\check{H}_{1}(\frak{L};T)_{k+1}\overset{(\varphi _{k+1,1})_{\ast }}\longrightarrow 
\check{H}_{1}(\frak{L};T)_{1}\longrightarrow 0.
\end{equation} 

 \item \label{puthm7b}
 $a_{1,k+1}=ma_{1,k}+a_{1,1}.$
 \item \label{puthm7c}
 If $ a_{1,1}=0, $ then 
 $a_{1,\infty }=0.$
 If $ a_{1,1}\neq 0, $
 then $a_{1,\infty }=\infty .$
 \item \label{puthm7d}
If $\check{H}_{1}(\frak{L};\ZZ )_{1}=0$, 
then, for each $k\in \NN $, 
$\check{H}_{1}(\frak{L};T)_{k}=0$ and $\check{H}^{1}(\frak{L};T)_{k}=0$, and 
 $\check{H}_{1}(\frak{L};T)=0$ and $\check{H}^{1}(\frak{L};T)=0.$  
 \item \label{puthm7e}
 There exists an exact sequence of $R$ modules: 
 \begin{equation}\label{puthm7eeq}
 0\longrightarrow \check{H}^{1}(\frak{L};R)_{1}
\overset{\mu _{1,1}}{\longrightarrow }\check{H}^{1}(\frak{L};R)
\overset{\theta }{\longrightarrow }
 \bigoplus _{j=1}^{m}\check{H}^{1}(\frak{L};R)\longrightarrow 0. 
 \end{equation}
\end{enumerate}

\end{enumerate}  
\end{thm}
We now give some important examples of postunbranched systems.
\begin{prop} 
\label{rgeq2exprop}\ 
\begin{enumerate}
\item \label{rgeq2exprop1}
For each $n\in \NN \cup \{ 0\} $, there exists a postunbranched backward self-similar system 
${\frak L}=(L, (h_{1},\ldots ,h_{n+2}))$ such that 
$X=\CCI $, $L\subset \CC $, $h_{j}:X\rightarrow X$ is a topological branched covering for each $j=1,\ldots ,n+2$, 
and $\dim _{R}\check{H}^{n}({\frak L};R)=\infty $ for each field $R.$  
In particular, if $n\geq 2$, then the above ${\frak L}$ satisfies that 
$\Psi :\check{H}^{n}({\frak L};R)\rightarrow \check{H}^{n}(L;R)$ is not a 
monomorphism for each field $R.$ 
\item \label{rgeq2exprop2} 
For each $n\in \NN \cup \{ 0\} $, there exists a postunbranched forward self-similar system 
${\frak L}=$ \\ $(L, (h_{1},\ldots ,h_{n+2}))$ such that 
$L\subset \RR^{3}$, $h_{j}:L\rightarrow L$ is injective for each $j=1,\ldots, n+2$, 
and $\dim _{R}\check{H}^{n}({\frak L};R)=\infty $ for each field $R.$ 
In particular, if $n\geq 3$, then the above ${\frak L}$ satisfies that 
$\Psi :\check{H}^{n}({\frak L};R)\rightarrow \check{H}^{n}(L;R)$ is not a 
monomorphism for each field $R.$ 
\end{enumerate}
\end{prop}
Theorem~\ref{puthm}-\ref{puthmmukrinj}  
implies that under the assumptions of Theorem~\ref{puthm}, for each nonnegative 
integer $r$ with $r\neq 1$, $\mu _{1,r}: \check{H}^{r}({\frak L};R)_{1}
\rightarrow \check{H}^{r}({\frak L};R)$ is a monomorphism. 
However, as illustrated in the following Proposition~\ref{notinjmuexprop}, 
even under the assumptions of Theorem~\ref{puthm}, if $|N_{1}|$ is disconnected, 
$\mu _{1,r}: \check{H}^{r}({\frak L};R)_{1}
\rightarrow \check{H}^{r}({\frak L};R)$ is not a monomorphism in general. 
\begin{prop}
\label{notinjmuexprop}
There exists a postunbranched forward self-similar system 
${\frak L}=(L,(h_{1},\ldots ,h_{5}))$ such that 
$L\subset \CC $, such that 
$h_{j}$ is a contracting similitude on $\CC $ (hence 
$h_{j}:L\rightarrow L$ is injective) for each $j=1,\ldots ,5$, 
and such that  
for each field $R$, we have 
$a_{1,1}\neq 0$, $B_{2}=0$, 
$\check{H}^{1}({\frak L};R)\cong  \check{H}^{1}(L;R)=0$, 
$|N_{1}|$ is disconnected, $\CC \setminus L$ is connected, 
and $\mu _{1,1}$ is not injective. See Figure~\ref{fig:mu11notinj1}.  
\end{prop} 
\begin{rem}
\label{r:nonly1}
Proposition~\ref{notinjmuexprop} means  
that for a postunbranched system ${\frak L}=(L,(h_{1},\ldots ,h_{m}))$, 
if $|N_{1}|$ is disconnected, then we need information on 
not only $\check{H}^{1}({\frak L};R)_{1}$ but also 
$B_{2}$ (or $B_{k}\ (k\geq 3)$), to determine  
$\check{H}^{1}({\frak L};R)_{k}\ (k\geq 2)$ and 
$\check{H}^{1}({\frak L};R).$ This provides us a new 
problem: ``Investigate $B_{k}$ of postunbranched systems with disconnected 
$|N_{1}|$.''  
\end{rem}
\begin{ex}[{\bf Sierpi\'{n}ski gasket}]
\label{SGcalex}
Let ${\frak L}=(L,(h_{1},h_{2},h_{3}))$ be the 
postunbranched forward self-similar system 
in Example \ref{ex:SG1}. (Hence $L$ is the Sierpi\'{n}ski gasket (Figure~\ref{fig:sgasket}).) 
We easily see that 
$|N_{1}|$ is connected,  
the set of all $1$-simplexes of  
$N_{1}$ is  
$\{ \{ 1,2\} ,\ \{ 1,3\} ,\{ 2,3\} \} $,  
and there exists no $r$-simplex of  
$N_{1}$, for each $r\geq 2.$ 
Let $R$ be a field. Then  
we have $\dim _{R}H^{1}(N_{1}; R)=1.$
Hence, by Theorem~\ref{puthm}, 
we obtain that for each $k\in \NN $, 
$ a_{1,k+1}=3a_{1,k}+1$, and that 
$\dim _{R}\check{H}^{1}({\frak L};R)=\dim _{R}\check{H}^{1}(L;R)
=\infty .$ Combining it with the Alexander duality theorem (\cite{Sp}), 
we see that $\CCI \setminus L$ has infinitely many connected components.  
Note that $\CCI \setminus L=F(\langle h_{1}^{-1},h_{2}^{-1} ,h_{3}^{-1}\rangle ).$ 
\end{ex}
\begin{ex}[{\bf Pentakun}]
\label{Pencalex} 
Let ${\frak L}=(L,(h_{1},\ldots ,h_{5}))$ be the forward self-similar
 system in \cite[Example 3.8.11]{K}. Hence $L$ is the pentakun (Figure~\ref{fig:snowflake2}). 
By Example~\ref{PSex}, ${\frak L}$ is postunbranched. 
Let $R$ be a field. 
By \cite[Example 3.8.11 (Pentakun)]{K} or Figure~\ref{fig:snowflake2}, 
we get that $|N_{1}|$ is connected and  
$\dim _{R}H^{1}(N_{1};R)=1.$ 
Hence, by Theorem~\ref{puthm}, 
we obtain that for each $k\in \NN $, 
$ a_{1,k+1}=5a_{1,k}+1$, and that 
$\dim _{R}\check{H}^{1}({\frak L};R)=\dim _{R}\check{H}^{1}(L;R)
=\infty .$ Combining it with the Alexander duality theorem (\cite{Sp}), 
we see that $\CCI \setminus L$ has infinitely many connected components. 
Note that $\CCI \setminus L=F(\langle h_{1}^{-1},\ldots ,h_{5}^{-1}\rangle ).$
\end{ex} 
\begin{ex}[{\bf Snowflake}]
\label{Sfcalex}
Let ${\frak L}=(L,(h_{1},\ldots ,h_{7}))$ be the 
forward self-similar system in \cite[Example 3.8.12 (Snowflake)]{K}. 
(Hence $L$ is the snowflake (Figure~\ref{fig:snowflake2}).)  
By Example~\ref{PSex}, ${\frak L}$ is postunbranched. 
Let $R$ be a field. 
By \cite[Example 3.8.12 (Snowflake)]{K} or Figure~\ref{fig:snowflake2},  
we get that $|N_{1}|$ is connected and  
$\dim _{R}H^{1}(N_{1};R)=6.$ 
Hence, by Theorem~\ref{puthm}, 
we obtain that for each $k\in \NN $, 
$ a_{1,k+1}=7a_{1,k}+6$, and that 
$\dim _{R}\check{H}^{1}({\frak L};R)=\dim _{R}\check{H}^{1}(L;R)
=\infty .$ Combining it with the Alexander duality theorem (\cite{Sp}), 
we see that $\CCI \setminus L$ has infinitely many connected components. 
Note that $\CCI \setminus L=F(\langle h_{1}^{-1},\ldots ,h_{7}^{-1}\rangle ).$ 
\end{ex}
\begin{ex}
\label{SGsub7ex}
Let ${\frak L}=(L,(h_{1},h_{2},h_{3}))$ be the 
postunbranched forward self-similar system 
in Example \ref{ex:SG1}. (Hence $L$ is the Sierpi\'{n}ski gasket (Figure~\ref{fig:sgasket}).)  
Let $g_{1}:=h_{1}^{2}, g_{2}:=h_{1}\circ h_{3}$, 
$g_{3}:=h_{2}^{2}, $ $g_{4}:= h_{2}\circ h_{3},$ $g_{5}:= h_{3}\circ h_{1}$, 
$g_{6}:= h_{3}\circ h_{2}$, and $g_{7}:= h_{3}^{2}.$ 
Let $L':=M_{\CC }(g_{1},\ldots, g_{7})$ and 
let ${\frak L}':=(L', (g_{1},\ldots ,g_{7})).$ 
For the figure of  $L'$, see Figure~\ref{fig:sgsubcal3-2}. 
By Lemma~\ref{puiteratelem} and Lemma~\ref{pssubsyslem}, 
${\frak L}'$ is postunbranched. 
It is easy to see that the set of $1$-simplexes of $N_{1}({\frak L}')$ is equal to 
$\{ \{ 1,2\} ,\{ 3,4\} ,\{ 5,6\} ,$ $ \{ 6,7\} ,\{ 7,5\}, \{ 2,5\} ,\{ 4,6\} \} $ 
and there exists no $r$-simplex of $N_{1}({\frak L}')$ for each $r\geq 2$ 
(cf. Figure~\ref{fig:sierpn1n2-2} and \ref{fig:sgsubcal3-2}).  
Thus $|N_{1}({\frak L}')|$ is connected and 
for each field $R$, $\check{H}^{1}({\frak L}';R)_{1}=R.$ 
Hence, by Theorem~\ref{puthm}, 
for each $k\in \NN $, 
$a_{1,k+1}=7a_{1,k}+1$ and $\dim _{R}\check{H}^{1}({\frak L}';R)=\dim _{R}\check{H}^{1}(L';R)=\infty .$ 
Combining it with the Alexander duality theorem (\cite{Sp}), 
$\CCI \setminus L'$ has infinitely many connected components. Note that 
$\CCI \setminus L'=F(\langle g_{1}^{-1},\ldots ,g_{7}^{-1}\rangle ).$ 
\begin{figure}[htbp]
\caption{The invariant set $L'$ in Example~\ref{SGsub7ex}}
\ \ \ \ \ \ \ \ \ \ \ \ \ \ \ \ \ \ \ \ \ \ \ \ \ \ \ \ \ \ \ 
\ \ \ \ \ \ \ \ \ \ \ \ \ \ \ \ \ \ \ \ \ \ 
\includegraphics[width=2.5cm,width=2.5cm]{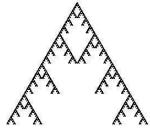}
\label{fig:sgsubcal3-2}
\end{figure}
\end{ex}
\begin{ex}
\label{SGsubcalex}
Let ${\frak L}=(L,(h_{1},h_{2},h_{3}))$ be the 
postunbranched forward self-similar system 
in Example \ref{ex:SG1}. (Hence $L$ is the Sierpi\'{n}ski gasket (Figure~\ref{fig:sgasket}).)  
Let $g_{1}:=h_{1}^{2}, g_{2}:=h_{1}\circ h_{2}$, 
$g_{3}:= h_{2}\circ h_{1}$, $g_{4}:= h_{2}^{2},$ 
$g_{5}:= h_{3}\circ h_{1}$, $g_{6}:= h_{3}\circ h_{2}$, and $g_{7}:= h_{3}^{2}.$ 
Let $L':= M_{\CC }(g_{1}\ldots ,g_{7})$ and let 
${\frak L}'=(L', (g_{1},\ldots ,g_{7})).$ For the figure of $L'$, 
see Figure~\ref{fig:sgsubcal2}. 
By Lemma~\ref{puiteratelem} and Lemma~\ref{pssubsyslem}, 
${\frak L}'$ is postunbranched. 
It is easy to see that the set of $1$-simplexes of $N_{1}({\frak L}')$ is 
equal to $\{ \{ 1,2\} ,\{ 2,3\} ,\{ 3,4\} ,\{ 5,6\} ,$ $\{ 6,7\} ,\{ 7,5\} \} $ and 
there exists no $r$-simplexes of $N_{1}({\frak L}')$ for each $r\geq 2$ 
(cf. Figures~\ref{fig:sierpn1n2-2} and \ref{fig:sgsubcal2}).   
Therefore $|N_{1}({\frak L}')|$ is disconnected and 
$\dim _{R}\check{H}^{1}({\frak L}';R)_{1}=1$ for each field $R.$ 
By Theorem~\ref{puthm}-\ref{puthm6} and Remark~\ref{r:Psiiso}, it follows that 
$\dim _{R}\check{H}^{0}({\frak L}';R)=\infty $ and $L'$ has infinitely many connected components. 
Note that $\CCI \setminus L'=F(\langle g_{1}^{-1},\ldots ,g_{7}^{-1}\rangle ).$ 
\begin{figure}[htbp]
\caption{The invariant set $L'$ in Example~\ref{SGsubcalex}}
\ \ \ \ \ \ \ \ \ \ \ \ \ \ \ \ \ \ \ \ \ \ \ \ \ \ \ \ \ \ \ 
\ \ \ \ \ \ \ \ \ \ \ \ \ \ \ \ \ \ \ \ \ \ \ 
\includegraphics[width=2.5cm,width=2.5cm]{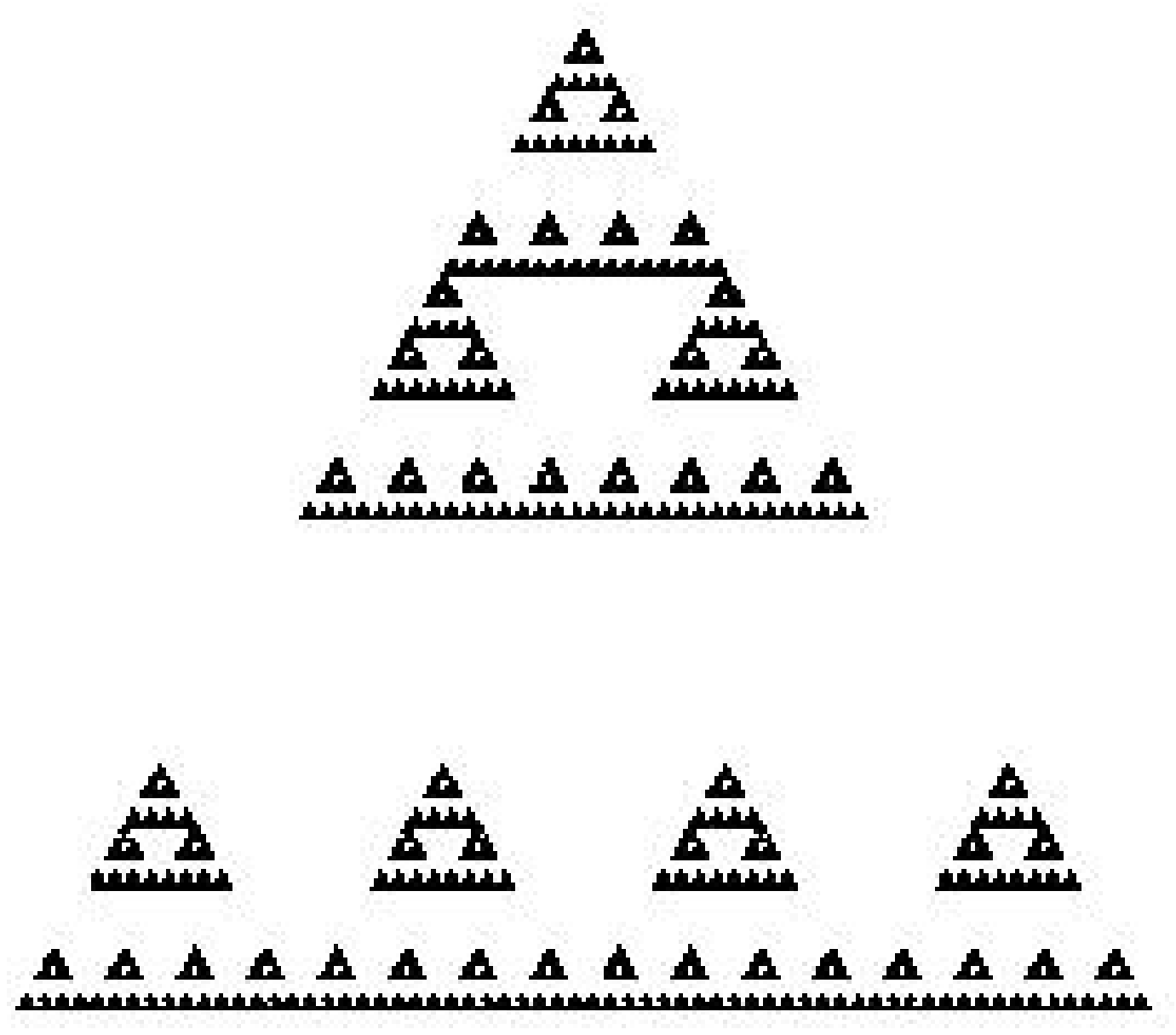}
\label{fig:sgsubcal2}
\end{figure}
\end{ex}
Regarding the postunbranched systems, we have the following lemma. 
\begin{lem}
\label{pucontcij1lem}
Let $\frak{L}=(L,(h_{1},\ldots ,h_{m}))$ be a postunbranched forward self-similar system such that 
for each $j=1,\ldots ,m$, 
$h_{j}: L\rightarrow L$ is a contraction. Then, for each $(i,j)\in \{ 1,\ldots ,m\} ^{2}$ with $i\neq j$, 
$\sharp C_{i,j}\leq 1.$ 
\end{lem}
\begin{proof}
Let $(i,j)\in \{ 1,\ldots m\} ^{2}$ be any element such that $i\neq j$ and $C_{i,j}\neq \emptyset .$ 
Since ${\frak L}$ is postunbranched, there exists an element 
$x\in \Sigma _{m}$ such that 
$h_{i}^{-1}(C_{i,j})\subset L_{x}.$ 
Since $h_{k}:L\rightarrow L$ is a contraction for each $k$, we have 
that $\sharp L_{x}=1.$ Hence $\sharp C_{i,j}\leq 1.$ 
\end{proof}
From Lemma~\ref{pucontcij1lem}, 
it is natural to consider the case $\sharp C_{i,j}\leq 1$ for each 
$(i,j)\in \{ 1,\ldots ,m\} ^{2}$ with $i\neq j.$ 
\begin{thm}\label{sscij1thm}
Let $\frak{L}=(L,(h_{1},\ldots ,h_{m}))$ be a forward self-similar system such that 
for each $j=1,\ldots ,m$, $h_{j}: L\rightarrow L$ is injective. 
 Let $T$ be a $\ZZ $ module and $R$ a field. 
Moreover, for each $r\in \NN \cup \{ 0\} $ and $k\in \NN $, 
let $a_{r,k} := \dim _{R}\check{H}^{r}({\frak L};R)_{k}.$ 
Furthermore, let $a_{1,\infty }:= \dim _{R}\check{H}^{1}({\frak L};R). $ 
Suppose that $\sharp C_{i,j}\leq 1 $ for each $(i,j)$ with $i\neq j$.  Then,  
we have the following. 
\begin{enumerate}
\item \label{sscij1thm-1}
 Let $k,r\in \NN $ with $r\geq 2$. Then,  
 $\check{H}_{r}({\frak L};T)_{k}=0$ and $\check{H}_{r}({\frak L};T)=0.$  
 \item \label{sscij1thm-2}
For each $k\in \NN ,$ $ma_{1,k}\leq a_{1,k+1}.$
 \item \label{sscij1thm-3}
If $|N_{1}|$ is 
 connected and $\check{H}^{1}({\frak L};R )\neq 0$,
  then 
  $a_{1,\infty }=\infty .$ 
\end{enumerate}  
\end{thm}
We present a result on the \v{C}ech cohomology groups of 
the invariant sets of the forward self-similar systems. 
This is also related to Lemma~\ref{pucontcij1lem}. 
\begin{prop}
\label{checkprop}
Let $X$ be a finite-dimensional topological manifold with a 
distance. Let $L$ be a non-empty compact subset of $X.$ Let $R$ be a field. Let $m\in \NN $ with $m\geq 2.$ 
 Let 
${\frak L}=(L,(h_{1},\ldots ,h_{m}))$ be a 
forward self-similar system.   
  Suppose that 
(a)for each $i=1,2$,  $h_{i}:L\rightarrow L$ is injective,  
  and 
(b) for each $(i,j)\in \{ 1,\ldots ,m\} ^{2}$ with 
$i\neq j$,  $\dim _{T}(C_{i,j})\leq n$ , where 
  $\dim _{T}$ denotes the topological dimension. 
  Then,\ $\dim _{R}\check{H}^{n+1}(L;R)$ 
  is either $0$ or $\infty .$ 
\end{prop}
\section{Tools}
\label{Tools}
In this section, we give some tools to show the main results.
\subsection{Fundamental properties of interaction cohomology}
In this subsection, we show some fundamental lemmas on the interaction (co)homology groups. 
We sometimes use the notation from \cite{Sp}. 
\begin{df}
Let ${\frak L}=(L,(h_{1},\ldots ,h_{m}))$ be a forward or backward self-similar system. 
For each $k\in \NN $, we denote by $\Gamma _{k}=\Gamma _{k}({\frak L})$ the 
$1$-dimensional skeleton of $N_{k}.$ 
\end{df}
\begin{lem}
\label{phionto} 
Let ${\frak L}=(L,(h_{1},\ldots ,h_{m}))$ be a forward or backward self-similar system.  
Then, 
for each $k\in \NN $, the simplicial map 
$\varphi _{k}: N_{k+1}\rightarrow N_{k}$ is 
surjective. That is, for any $r\in \NN $, 
if $x=\{ x^{1},\ldots ,x^{r}\} $ is an $r-1$ simplex of  
$N_{k}$, then there exists an 
$r-1$ simplex $y=\{ y^{1},\ldots ,y^{r}\} $ of $N_{k+1}$ such that 
$\varphi _{k}(y)=x.$ 
In particular, 
$ (\varphi _{k})_{\ast }: 
\mbox{{\em Con}}(|\Gamma _{k+1}|)\rightarrow 
\mbox{{\em Con}}(|\Gamma _{k}|)$ is surjective.
\end{lem} 
\begin{proof}
We will prove the statement of our lemma when ${\frak L}$ is a backward self-similar 
system (when ${\frak L}$ is 
a forward self-similar system, we can prove the statement by using an argument similar to the below). 
Let $x=\{ x^{1},\ldots ,x^{r}\} $ be an $r-1$ simplex of $N_{k}$, 
where for each $j=1,\ldots ,r$, $x^{j}=(x_{1}^{j},\ldots ,x_{k}^{j})\in \{ 1,\ldots ,m\} ^{k}.$ 
Then 
$\bigcap _{j=1}^{r}h_{x_{1}^{j}}^{-1}\cdots 
h_{x_{k}^{j}}^{-1}(L)\neq \emptyset .$
Let 
$z\in \bigcap _{j=1}^{r}h_{x_{1}^{j}}^{-1}\cdots 
h_{x_{k}^{j}}^{-1}(L). $ 
Then for each $j=1,\ldots ,r$, 
$h_{x_{k}^{j}}\cdots h_{x_{1}^{j}}(z)\in L=\bigcup _{i=1}^{m}h_{i}^{-1}(L).$ 
Hence, for each $j=1,\ldots, r$, there exists an 
$x_{k+1}^{j}\in \{ 1,\ldots ,m\} $ such that 
$h_{x_{k+1}^{j}}\cdots h_{x_{1}^{j}}(z)\in L.$ 
Therefore,  
$\bigcap _{j=1}^{r}h_{x_{1}^{j}}^{-1}\cdots h_{x_{k+1}^{j}}^{-1}(L)\neq \emptyset .$ 
Thus, setting $y^{j}:= (x_{1}^{j},\ldots ,x_{k+1}^{j})\in \{ 1,\ldots ,m\} ^{k+1}$ 
for each $j=1,\ldots ,r$, we have that $y=\{ y^{1},\ldots ,y^{r}\} $ is an $r-1$ simplex 
 of $N_{k+1}$  such that $\varphi _{k}(y)=x.$ 
\end{proof}
\begin{lem}
\label{1conall}
Let ${\frak L}=(L,(h_{1},\ldots ,h_{m}))$ 
be a forward or backward self-similar system. 
If 
$|\Gamma _{1}|$ is connected, then, 
for any $k\in \NN $, $|\Gamma _{k}|$ and $|N_{k}| $ are connected.

\end{lem}
\begin{proof}
We will prove the statement of our lemma when ${\frak L}$ is a backward self-similar system 
(when ${\frak L}$ is a forward self-similar system, we can prove the statement of our lemma by using 
an argument similar to the below). 
First, we show the following claim.\\ 
Claim:  Let 
 $w^{1}$ and $w^{2}$ be two elements in 
 $\{ 1,\ldots ,m\} ^{k}$ such that 
 $h_{w^{1}}^{-1}(L)\cap h_{w^{2}}^{-1}(L)\neq \emptyset .$ 
 Then, 
  for any $j_{1}$ and $j_{2}$ in 
 $\{ 1,\ldots ,m\} $, 
 there exists an edge path $\gamma $ 
 of $\Gamma _{k+1}$ 
 from $w^{1}j_{1}$ to 
 $w^{2}j_{2}.$ (For the definition of edge path, see \cite{Sp}.) 

 To show this claim, 
 since 
 $L=\bigcup _{j=1}^{m}h_{j}^{-1}(L)$, 
 we obtain that there exist 
 $w_{k+1}^{1}$ and $w_{k+1}^{2}$ in 
 $\{ 1,\ldots ,m\} $ such that 
 $h_{w^{1}}^{-1}h_{w_{k+1}^{1}}^{-1}(L)\cap 
 h_{w^{2}}^{-1}h_{w_{k+1}^{2}}^{-1}(L)\neq \emptyset .$ 
 Hence, there exists an edge 
 path $\alpha $ of $\Gamma _{k+1}$ from  
 $w^{1}w_{k+1}^{1}$ to 
 $w^{2}w_{k+1}^{2}.$ 
 Furthermore, since $|\Gamma _{1}|$ is 
 connected, we have that for each $i=1,2$, 
 there exists an edge path 
 $\tau _{i}$ of $\Gamma _{1}$ from  
 $j_{i}$ to $w_{k+1}^{i}.$ 
 Then, for each $i=1,2$, there exists an edge path 
 $\beta _{i}$ 
 of $\Gamma _{k+1}$ from  
 $w^{i}j_{i}$ to $w^{i}w_{k+1}^{i}.$ 
  Hence, there exists an edge path 
  of $\Gamma _{k+1}$ from  
  $w^{1}j_{1}$ to 
  $w^{2}j_{2}.$ Therefore, we have shown the above claim. 

 We now show the statement of our lemma by induction on $k.$ 
 Suppose that $|\Gamma _{k}|$ is connected. 
Let $x$ and $y$ be any elements in $\{ 1,\ldots ,m\} ^{k+1}.$ 
Then, there exists an edge 
path of $\Gamma _{k}$ from  
$x|k$ and $y|k.$ 
By the above claim, 
we easily obtain that there exists an edge path 
of $\Gamma _{k+1}$ from  
$x$ and $y.$ Hence, $|\Gamma _{k+1}|$ is connected. 
Thus, the induction is completed.  
\end{proof}

\begin{df}
Let $K$ be a simplicial complex and let $R$ be a $\ZZ $ module. 
We denote by $C_{\ast }(K)$ the oriented chain complex of $K$ 
(\cite[p.159]{Sp}). Moreover,  
we set $C_{\ast }(K;R):=C_{\ast }(K)\otimes R$ and 
$C^{\ast }(K;R):=\mbox{Hom}(C_{\ast }(K),R).$ 
Similarly, we denote by $\triangle _{\ast }(K)$ the 
ordered chain complex of $K$ (\cite[p.170]{Sp}) and we set 
$\triangle _{\ast }(K;R):=\triangle _{\ast }(K)\otimes R$ and 
$\triangle ^{\ast }(K;R):=\mbox{Hom}(\triangle _{\ast }(K),R).$  
Moreover, for a relative CW complex $(X,A)$, we denote by 
$C_{\ast }(X,A)$ the chain complex given in \cite[p. 475]{Sp}. 
Furthermore, we set $C_{\ast }(X,A;R):=C_{\ast }(X,A)\otimes R$ and 
$C^{\ast }(X,A;R):=\mbox{Hom}(C_{\ast }(X,A),R).$ 
\end{df}
\begin{df}
Let $X$ be a topological space and let $R$ be a $\ZZ $ module. 
We regard $R$ as a constant presheaf on $X$ (\cite[p. 323]{Sp}). 
Moreover, we denote by $\hat{R}$ the completion of the presheaf $R$ 
(\cite[p. 325]{Sp}).  
Thus $\hat{R}$ is a sheaf assigning to each non-empty open subset $U$ of $X$ 
the $\ZZ $ module of all locally constant 
functions $a: U\rightarrow R.$ Moreover, for an open covering ${\cal U}$ of $X$ and 
a presheaf $\Gamma $ on $X$, we denote by $C^{\ast }({\cal U};\Gamma )$ the 
cochain complex in \cite[p. 327]{Sp} and $H^{\ast }({\cal U};\Gamma )$ its cohomology group.
 Note that by definition, $\check{H}^{\ast }(X;\Gamma )=\varinjlim_{{\cal U}}H^{\ast }({\cal U};\Gamma )$ 
(\cite[p. 327]{Sp}).  
\end{df}
\begin{rem}
\label{r:rrhat}
There is a natural homomorphism $\alpha :R\rightarrow \hat{R}$ 
such that for each open subset $U$ of $X$, $\alpha $ assigns $\gamma \in R(U)$ to 
locally constant function $\hat{\gamma }:U\rightarrow R$ with $\hat{\gamma }(a)=\gamma $ for all $a\in U$. (See \cite[p. 325]{Sp}.) 
Thus $\alpha $ induces a natural homomorphism $\alpha _{\ast }: C^{\ast }({\cal U};R)\rightarrow 
C^{\ast }({\cal U};\hat{R})$ for any open covering ${\cal U}$ of $X.$  
\end{rem}
\begin{lem}
\label{injlem1}
Let $(L,d)$ be a compact metric space. 
Let ${\cal A}=\{ L_{i}\} _{i=1}^{r}$ be a 
finite covering of $L$ such that for each $i=1,\ldots ,r$, 
$L_{i}$ is a non-empty compact subset of $L.$
Let $\delta ({\cal A})$ be the number in Lemma~\ref{naturalhomlem}.
Let $0<\delta <\delta ({\cal A})$ and 
let $\psi _{0}:\triangle ^{\ast }(N({\cal A});R)\cong 
\triangle ^{\ast }(N({\cal A}^{\delta });R)\rightarrow 
C^{\ast }({\cal A}^{\delta }; R)$ be the natural homomorphism. 
Let $\psi :\triangle ^{\ast }(N({\cal A});R)\rightarrow C^{\ast }({\cal A}^{\delta }; \hat{R})$ 
be the composition $\alpha _{\ast }\circ \psi _{0}.$  
Moreover, let 
$\psi _{\ast }: H^{\ast }(N({\cal A});R)\cong H^{\ast }({\cal A}^{\delta };R)
\rightarrow H^{\ast }({\cal A}^{\delta };\hat{R})$ 
be the homomorphism induced by $\psi .$   
Then, we have the following.
\begin{enumerate}
\item \label{injlem1-1}
$\psi _{\ast }: H^{0}(N({\cal A});R)\cong H^{0}({\cal A}^{\delta };R)\rightarrow 
H^{0}({\cal A}^{\delta };\hat{R})$ is a monomorphism.
\item \label{injlem1-2} 
In addition to the assumptions of the lemma, suppose that 
for each $i=1,\ldots ,r$, $L_{i}$ is connected. Then,  
$\psi _{\ast }:H^{1}(N({\cal A});R)\cong H^{1}({\cal A}^{\delta };R)\rightarrow  
H^{1}({\cal A}^{\delta };\hat{R})$ is a monomorphism. 
Moreover, the natural homomorphism 
$\Psi _{{\cal A}}: H^{1}(N({\cal A});R)\cong H^{1}({\cal A}^{\delta };R)
\rightarrow \check{H}^{1}(L;R)$ is monomorphic.  
\end{enumerate}
\end{lem}
\begin{proof}
It is easy to see that statement \ref{injlem1-1} holds. 
We now prove statement \ref{injlem1-2}. 
Let $a=(a_{ij})_{(i,j): L_{i}\cap L_{j}\neq \emptyset }\in \triangle ^{1}(N({\cal A});R)$ be a cocycle, 
where $a_{ij}:L_{i}\cap L_{j}\rightarrow R$ is a constant function for each 
$(i,j)$ with $L_{i}\cap L_{j}\neq \emptyset .$ 
We write $\psi (a)$ as $(b_{ij})_{(i,j):L_{i}\cap L_{j}\neq \emptyset }$, 
where $b_{ij}:B(L_{i},\delta )\cap B(L_{j},\delta )\rightarrow R$ is a constant function 
which is an extension of $a_{ij}.$  
Suppose that $\psi (a)\in C^{1}({\cal A}^{\delta };\hat{R})$ is a coboundary. 
Then, there exists an element $(b_{i})_{i=1,\ldots ,r}\in C ^{0}({\cal A}^{\delta };\hat{R})$, 
where each $b_{i}:B(L_{i}:\delta )\rightarrow R$ is a locally constant function,  
such that $b_{ij}=b_{j}-b_{i}$ on $B(L_{i},\delta )\cap B(L_{j},\delta ).$ 
Hence 
\begin{equation}
\label{eq:injlem1eq1}
a_{ij}=\left((b_{j}|_{L_{j}})-(b_{i}|_{L_{i}})\right)|_{L_{i}\cap L_{j}} \mbox{ on } 
L_{i}\cap L_{j}.
\end{equation}
Moreover, for each $i$, since $L_{i}$ is connected and $b_{i}:B(L_{i},\delta )\rightarrow R$ is 
locally constant, we have that $b_{i}|_{L_{i}}:L_{i}\rightarrow R$ is constant. 
Combining it with (\ref{eq:injlem1eq1}), we obtain that 
$a$ is a coboundary. Thus, 
we have proved that $\psi _{\ast }:H^{1}(N({\cal A});R)\cong H^{1}({\cal A}^{\delta };R)\rightarrow  
H^{1}({\cal A}^{\delta };\hat{R})$ is a monomorphism. 
Moreover, by Leray's theorem (\cite[Theorem 5 in page 56 and Theorem 11 in page 61]{G}), 
the natural homomorphism $H^{1}({\cal A}^\delta ;\hat{R})\rightarrow \check{H}^{1}(L;\hat{R})$ is 
monomorphic. Furthermore, by \cite[p.329]{Sp}, 
the natural homomorphism $\check{H}^{1}(L;R)\rightarrow \check{H}^{1}(L,\hat{R})$ is isomorphic. 
Therefore, the natural homomorphism 
$\Psi _{{\cal A}}: H^{1}(N({\cal A});R)\rightarrow \check{H}^{1}(L;R)$ is monomorphic. 
Thus, we have proved statement \ref{injlem1-2}. 
\end{proof}

\begin{lem}
\label{icfundlem}
Let ${\frak L}=(L,(h_{1},\ldots ,h_{m}))$ be a forward or backward self-similar 
system. Let $G=\langle h_{1},\ldots ,h_{m}\rangle .$
Let $R$ be a $\ZZ $ module. 
Then, we have the following:
\begin{enumerate}
\item \label{icfundlem1}
For each $k\in \NN ,$ $\varphi _{k}^{\ast }:
\check{H}^{0}({\frak L};R)_{k}\rightarrow 
\check{H}^{0}({\frak L};R)_{k+1}$ is 
a monomorphism. In particular, 
for each $k\in \NN $, the projection map 
$\mu _{k,0}:\check{H}^{0}({\frak L};R)_{k}\rightarrow 
\check{H}^{0}({\frak L};R)$ is injective. 
\item \label{icfundlem2}
$\Psi :\check{H}^{0}({\frak L};R)\rightarrow 
\check{H}^{0}(L;R)$ is a monomorphism.
\item \label{icfundlem3-01}
Suppose that $|N_{1}|$ is connected. Then,  
for each $k$ and each $w\in \Sigma _{m}$, 
$(\varphi _{k})_{\ast }:
\check{\pi}_{1}({\frak L},w)_{k+1}\rightarrow 
\check{\pi}_{1}({\frak L},w)_{k}$ is an epimorphism and 
$(\varphi _{k})_{\ast }:
\check{H}_{1}({\frak L}; R )_{k+1}\rightarrow 
\check{H}_{1}({\frak L}; R )_{k}$ is an epimorphism. 
\item \label{icfundlem3}
Suppose that $|N_{1}|$ is connected. Then,  
for each $k$, 
$\varphi _{k}^{\ast }:
\check{H}^{1}({\frak L};R)_{k}\rightarrow 
\check{H}^{1}({\frak L};R)_{k+1}$ is 
a monomorphism and  
the projection map 
$\mu _{k,1}: \check{H}^{1}({\frak L}; R)_{k}\rightarrow 
\check{H}^{1}({\frak L};R)$ is a monomorphism. 
\item \label{h1inj}
Suppose that either (a) ${\frak L}$ is a forward self-similar system and $L$ is connected, 
or (b) ${\frak L}$ is a backward self-similar system such that $g^{-1}(L)$ is connected for each 
$g\in G$. Then,  
for each $k\in \NN $, 
the natural homomorphism 
$\Psi _{{\cal U}_{k}}:H^{p}(N_{k};R)\rightarrow \check{H}^{p}(L;R)$ in 
Remark~\ref{nathomrem} is monomorphic,   
$\Psi :\check{H}^{1}({\frak L};R)\rightarrow 
\check{H}^{1}(L;R)$ is a monomorphism, and 
for each $k\in \NN $, 
$\varphi _{k}^{\ast }:\check{H}^{1}({\frak L};R)_{k}\rightarrow \check{H}^{1}({\frak L};R)_{k+1}$ 
is a monomorphism. 
\end{enumerate}
\end{lem}
\begin{proof}
It is easy to see that statement \ref{icfundlem1} holds. 
Using Lemma~\ref{naturalhomlem}, it is easy to see that statement 
\ref{icfundlem2} holds. 

We now prove statements \ref{icfundlem3-01} and \ref{icfundlem3}. 
If $|N_{1}|$ is connected, then 
Lemma~\ref{1conall} implies that for each $k\in \NN $, 
$|N_{k}|$ is connected. 
Let $w\in \Sigma _{m}.$ 
Let $\zeta \in \pi _{1}(|N_{k}|,w)$ be an element. 
We use the notation in \cite{Sp}. 
By \cite{Sp}, there exists a closed edge path $\gamma =
\gamma _{1}\gamma _{2}\cdots \gamma _{r}$, 
where each $\gamma _{j}=(x^{j},x^{j+1})$ is an edge of $N_{k}$, 
such that $\gamma $ represents the element $\zeta .$ 
For each $j=1,\ldots ,r+1$, we write $x^{j}$ 
as $(x_{1}^{j},\ldots ,x_{k}^{j})\in \{ 1,\ldots ,m\} ^{k}.$  
By Lemma~\ref{phionto}, 
for each $j=1,\ldots ,r$ 
there exists an edge $\tau _{j}$ of $N_{k+1}$ such that 
$\varphi _{k}(\tau _{j})=\gamma _{j}.$ Then, there exists 
$y^{j}, z^{j}\in \{ 1,\ldots ,m\} $ such that 
the origin of $\tau _{j}$ is equal to $x^{j}y^{j}$ and the end of 
$\tau _{j}$ is  
 equal to $x^{j+1}z^{j}.$ Since we are assuming that 
$|N_{1}|$ is connected, 
for each $j=2,\ldots, r$, there exists an edge path 
$\beta _{j}=(v_{j}^{1},v_{j}^{2})(v_{j}^{2},v_{j}^{3})\cdots (v_{j}^{s_{j}-1},v_{j}^{s_{j}})$ of 
$N_{1}$, 
where each $v_{j}^{k}$ is a vertex of $N_{1}$, 
such that $v_{j}^{1}=z^{j-1}$ and $v_{j}^{s_{j}}=y^{j}.$ 
Similarly, there exists an edge path 
$\beta _{r+1}=(v_{r+1}^{1},v_{r+1}^{2})\cdots (v_{r+1}^{s_{r+1}-1},v_{r+1}^{s_{r+1}})$ of $N_{1}$ 
such that $v_{r+1}^{1}=z^{r}$ and $v_{r+1}^{s_{r+1}}=y^{1}.$ 
For each $j=2,\ldots, r+1$, let 
$\delta _{j}:= 
(x^{j}v_{j}^{1},x^{j}v_{j}^{2})\cdots (x^{j}v_{j}^{s_{j}-1},x^{j}v_{j}^{s_{j}}).$ 
Then, for each $j=2,\ldots ,r$, $\delta _{j}$ is an edge path of $N_{k+1}$ from  
$x^{j}z^{j-1}$ to  $x^{j}y^{j}.$ Moreover, 
$\delta _{r+1}$ is an edge path of $N_{k+1}$ from  
$x^{r+1}z^{r}$ to $x^{1}y^{1}.$ 
Let $\delta := \tau _{1}\delta _{2}\tau _{2}\delta _{3}\tau _{3}\delta _{4}\cdots 
\tau _{r}\delta _{r+1}.$ Then, $\delta $ is a closed edge path of $N_{k+1}$ 
such that $\varphi _{k}(\delta )\sim \gamma .$ Therefore, 
$(\varphi _{k})_{\ast }: \check{\pi }_{1}({\frak L},w)_{k+1}\rightarrow 
\check{\pi }_{1}({\frak L},w)_{k}$ is an epimorphism. 
Moreover, by Lemma~\ref{1conall} and \cite[p.394]{Sp}, 
for each $k\in \NN $, the natural homomorphism 
$\check{\pi }_{1}({\frak L},w)_{k}\rightarrow \check{H}_{1}({\frak L}; \ZZ )_{k}$ is 
an epimorphism.  Therefore, 
$(\varphi _{k})_{\ast }: \check{H }_{1}({\frak L}; \ZZ )_{k+1}\rightarrow 
\check{H }_{1}({\frak L}; \ZZ )_{k}$ is an epimorphism. 
From the universal-coefficient theorem for homology (\cite[p.222]{Sp}), 
it follows that for any $\ZZ $ module $R$,  
$(\varphi _{k})_{\ast }: \check{H }_{1}({\frak L}; R)_{k+1}\rightarrow 
\check{H }_{1}({\frak L}; R )_{k}$ is an epimorphism. 
Similarly, from the universal-coefficient theorem for cohomology (\cite[p.243]{Sp}), 
it follows that for any $\ZZ $ module $R$, 
$(\varphi _{k})_{\ast }: \check{H }^{1}({\frak L};R )_{k}\rightarrow 
\check{H }^{1}({\frak L};R )_{k+1}$ is a monomorphism. 
Thus we have proved statement \ref{icfundlem3-01} and \ref{icfundlem3}. 
 
 Statement \ref{h1inj} follows from  Lemma~\ref{injlem1}.

Hence, we have completed the proof of Lemma~\ref{icfundlem}.     
\end{proof}
\begin{ex}
\label{ex:finsetic}
Let $L=\{ a_{1},a_{2},a_{3}\} $ be a set where $a_{j}$ are mutually distinct points. 
Let $h_{1}:L\rightarrow L$ be the  map defined by 
$h_{1}(a_{1})=a_{1}, h_{1}(a_{2})=h_{1}(a_{3})=a_{2}.$ 
Similarly, let $h_{2}:L\rightarrow L$ be the map defined by 
$h_{2}(a_{2})=a_{2}, h_{2}(a_{3})=h_{2}(a_{1})=a_{3}.$ 
Finally, let $h_{3}:L\rightarrow L$ be the map defined by 
$h_{3}(a_{3})=a_{3}, h_{3}(a_{1})=h_{3}(a_{2})=a_{1}.$ 
Then ${\frak L}:= (L, (h_{1},h_{2},h_{3}))$ is a 
forward self-similar system. 
It is easy to see that 
the set of one-dimensional simplexes of $N_{1}$ is equal to 
$\{ \{ 1,2\} , \{ 2,3 \} , \{ 3,1\} \} $ and for each $r\geq 2$, 
there exists no $r$-dimensional simplex of $N_{1}.$  
Therefore $|N_{1}|$ is connected and 
for each $\ZZ $ module $R$, 
$\check{H}^{1}({\frak L};R)_{1}=R\neq 0.$ 
By Lemma~\ref{icfundlem}, 
it follows that $\check{H}_{1}({\frak L};R)\neq 0$, 
$\check{H}^{1}({\frak L};R)\neq 0$, and 
$\check{\pi }_{1}({\frak L},w)$ is not trivial for each 
$w\in \Sigma _{m}.$ However, 
since $L$ is a finite set, 
$\check{\pi }_{r}(L,x)$, 
$\check{H}^{r}(L;R)$ and $\check{H}_{r}(L;R)$ are trivial for each $r\geq 1$ and 
each $x\in L.$ This example means that the interaction cohomology groups  
of self-similar systems may have more  information than the 
\v{C}ech (co)homology groups and the shape groups of the invariant sets of the systems.   
\end{ex}
\begin{ex}
\label{ex:fintriv}
Let $L=\{ a_{1},a_{2},a_{3}\} $ be a set where $a_{j}$ are mutually distinct points. 
For each $j=1,2,3,$  
let $g_{j}:L\rightarrow L$ be the  map defined by
 $g_{j}(L)=\{ a_{j}\} .$ 
Then ${\frak L}'=(L, (g_{1},g_{2},g_{3}))$ is a forward self-similar system. 
It is easy to see that for each $r\geq 1$, there exists no $r$-dimensional simplexes of $N_{1}.$ 
Moreover, since each $g_{j}$ is a contraction and $L$ is a finite set, 
it follows that  
$\check{H}^{r}({\frak L}';R)_{k}$, $\check{H}^{r}({\frak L}';R)$, 
$\check{H}_{r}({\frak L}';R)_{k}$, $\check{H}_{r}({\frak L}';R)_{k}$, 
 and $\check{\pi }_{r}({\frak L};w)$ are trivial, for all $r\geq 1, $ $k\geq 1$, $w\in \Sigma _{3}$, 
and $\ZZ $ modules $R.$  
\end{ex}
\begin{rem}
\label{r:icdlh}
Example~\ref{ex:finsetic} and Example~\ref{ex:fintriv} mean that 
for any two self-similar systems ${\frak L}_{1}=(L_{1},(h_{1},\ldots ,h_{m}))$ 
and ${\frak L}_{2}=(L_{2},(g_{1},\ldots ,g_{n}))$, 
interaction cohomology groups $\check{H}^{r}({\frak L}_{1};R)$ and 
$\check{H}^{r}({\frak L}_{2};R)$ may not be isomorphic even when  
$L_{1}$ and $L_{2}$ are homeomorphic.  
\end{rem}
\subsection{Fundamental properties of rational semigroups}
We give some fundamental properties of rational semigroups. 
Let $G$ be a rational semigroup. We set 
$ E(G):= \{z\in \CCI \mid \sharp \bigcup _{g\in G}
g^{-1}(z) <\infty  \} .$ This is called the 
{\bf exceptional set} of $G.$ 
If \( z\in  \CCI \setminus E(G),\ \) then 
 \( J(G)\subset \overline{\bigcup _{g\in G}g^{-1}(\{ z\} )} .\) 
 In particular 
if  \( z\in J(G)\setminus E(G), \) \ then
$ \overline{\bigcup _{g\in G}g^{-1}(\{ z\} )}=J(G). $
If \( \sharp J(G) \geq 3\) ,\ then \( J(G) \) is a 
perfect set, $ \sharp E(G) \leq 2 $, 
\( J(G) \) is the smallest in 
$\{ K\subset \CCI \mid K:\mbox{compact},\ \sharp K\geq 3,\ 
 \mbox{and }g^{-1}(K)\subset K \mbox{ for each }g\in G\} ,$
%
and 
$$J(G)=\overline{\{ z\in \CCI \mid \exists g\in G \mbox{ s.t. } g(z)=z \mbox{ and }|g'(z)|>1\} } 
=\overline{\bigcup _{g\in G}J(g)}.$$ 
For the proofs of these results, see \cite{HM,GR} and \cite[Lemma 2.3]{S3}. 

\subsection{Fiberwise (Wordwise) dynamics}
In this subsection, we give some notations and fundamental 
properties of skew products related to finitely generated 
rational semigroups. 
\begin{df}[\cite{S1, S3}]
Let $G=\langle h_{1},\ldots ,h_{m}\rangle $ be a 
finitely generated rational semigroup. We define  
a map $\sigma :\Sigma _{m}\rightarrow \Sigma _{m}$ by:   
$\sigma (x_{1},x_{2},\ldots ):=(x_{2},x_{3},\ldots )
.$ This is called the {\bf shift map} on $\Sigma _{m}.$ Moreover,\ 
we define a map
$f:\Sigma _{m}\times \CCI \rightarrow 
\Sigma _{m}\times \CCI $ by:
$(x,y) \mapsto (\sigma (x),\ h_{x_{1}}(y)),\ $
where $x=(x_{1},x_{2},\ldots ).$
This is called 
{\bf the skew product associated with 
the multi-map} $(h_{1},\ldots ,h_{m})\in (\mbox{Rat})^{m}.$ 
Let $\pi :\Sigma _{m}\times \CCI \rightarrow 
\Sigma _{m}$ and  
$\pi _{\CCI }:\Sigma _{m}\times \CCI 
\rightarrow  \CCI $ be the projections. 
For each $x\in \Sigma _{m}$ and each $n\in \NN $, we set  
$f_{x}^{n}:=
f^{n}|_{\pi ^{-1}(\{ x\} )}:\pi ^{-1}(\{ x\} )\rightarrow 
\pi ^{-1}(\{ \sigma ^{n}(x)\} )\subset \Sigma _{m}\times \CCI $ 
and $f_{x,n}:=f_{x_{n}}\circ \cdots \circ f_{x_{1}}.$ 
Moreover, 
we denote by $F_{x}(f)$ 
the set of all points 
$y\in \CCI $ which has a neighborhood $U$ 
in $\CCI $ such that 
$\{ f_{x,n}:U\rightarrow 
\CCI \} _{n\in \NN }$
is normal on $U.$  
Moreover, we set  
$J_{x}(f):=\CCI \setminus 
F_{x}(f).$ Furthermore, we set 
$F^{x}(f):= \{ x\} \times F_{x}(f)$ and 
$J^{x}(f):=\{ x \} \times J_{x}(f).$  
We set $\tilde{J}(f):=
\overline {\bigcup _{x\in \Sigma _{m}}J^{x}
(f)}$, where the closure is taken in the product space $\Sigma _{m}\times \CCI .$
Moreover, for each $x\in \Sigma _{m}$, 
we set $\hat{J}^{x}(f):=\pi ^{-1}(\{ x\} )\cap \tilde{J}(f)$ 
and $\hat{J}_{x}(f):= \pi _{\CCI }(\hat{J}^{x}(f)).$ 
Furthermore, we set $\tilde{F}(f):=(\Sigma _{m}\times \CCI)\setminus 
\tilde{J}(f).$
\end{df}
\begin{rem} (See \cite[Lemma 2.4]{S1}.) 
$\tilde{J}(f)$, $J_{x}(f), J^{x}(f), \hat{J}_{x}(f), $ and $ \hat{J}^{x}(f)$ are compact. 
We have that $f^{-1}(\tilde{J}(f))=
\tilde{J}(f)=
f(\tilde{J}(f))$, 
$f^{-1}J^{\sigma (x)}(f)
=J^{x}(f)$,  
$f^{-1}\hat{J}^{\sigma (x)}(f)=
\hat{J}^{x}(f)$, and  
$\hat{J}_{x}(f)\supset 
J_{x}(f).$ However, 
the equality $\hat{J}_{x}(f)=J_{x}(f)$ does not hold in general. 
(This is one of the difficulties when we investigate the dynamics of rational semigroups or 
random complex dynamics.)  
\end{rem}
\begin{rem}
[\cite{J,S1}]
\label{lowersemicontirem}
({\bf Lower semicontinuity of 
$x\mapsto J_{x}(f)$})
Suppose that  $\deg (h_{j})\geq 2$ 
for each $j=1,\ldots ,m .$ 
Then, for each $x\in \Sigma _{m}$, 
$J_{x}(f)$ is a non-empty perfect set.
Furthermore, 
$x\mapsto J_{x}(f)$ is 
{\bf lower semicontinuous}, that is,  
for any point $y\in J_{x}(f)$  
and 
any sequence $\{ x^{n}\} _{n\in \NN }$ in $\Sigma _{m}$ 
with $x^{n}\rightarrow x,\ $  
there exists a sequence $\{ y^{n}\} _{n\in \NN }$ 
in $\CCI $ with 
$y^{n}\in J_{x^{n}}(f)\ (\forall n)$ such that 
$y^{n}\rightarrow y.$ 
The above result was shown by using the potential theory. 
We remark that  
$x\mapsto J_{x}(f)$ is not
 continuous with respect to the Hausdorff topology in 
general.
\end{rem}
\begin{lem}\label{fundfiblem}
Let $(h_{1},\ldots ,h_{m})\in (\mbox{{\em Rat}})^{m}$ and 
let $f:\Sigma _{m}\times \CCI \rightarrow \Sigma _{m}\times \CCI $ be 
the skew product associated with $(h_{1},\ldots ,h_{m}).$ 
Let $G=\langle h_{1},\ldots ,h_{m}\rangle .$ Suppose 
$\sharp J(G)\geq 3.$ Then, 
$\pi _{\CCI }(\tilde{J}(f))=J(G)$ and 
for each 
$x=(x_{1},x_{2},\ldots )\in \Sigma _{m}$, 
$\hat{J}_{x}(f)
=\bigcap _{j=1}^{\infty }h_{x_{1}}^{-1}
\cdots h_{x_{j}}^{-1}(J(G)).$  
\end{lem}
\begin{proof}
Since $J_{x }(f)\subset J(G)$ for each $x \in \Sigma _{m} $, 
we have $\pi _{\CCI }(\tilde{J}(f))\subset J(G).$ 
By \cite[Corollary 3.1]{HM} (see also \cite[Lemma 2.3 (g)]{S3}), we have 
$J(G)=\overline{\bigcup _{g\in G}J(g)}.$ 
Since $\bigcup _{g\in G}J(g)\subset \pi _{\CCI }(\tilde{J}(f))$, 
we obtain $J(G)\subset \pi _{\CCI }(\tilde{J}(f)).$ 
Therefore, we obtain $\pi _{\CCI }(\tilde{J}(f))=J(G).$ 

 We now show the latter statement. 
Let $x=(x_{1},x_{2}\ldots )\in \Sigma _{m} .$ 
By \cite[Lemma 2.4]{S1}, 
we see that for each $j\in \NN $, 
$h_{x_{j}}\cdots h_{x_{1}}(
\hat{J}_{x }(f))=
\hat{J}_{\sigma ^{j}(x )}(f)
\subset J(G).$ 
Hence, 
$\hat{J}_{x }(f)\subset 
\bigcap _{j=1}^{\infty }
h_{x _{1}}^{-1}\cdots 
h_{x _{j}}^{-1}(J(G)).$
 Suppose that there exists a point 
 $(x ,y)\in \Sigma _{m} \times \CCI $ such that  
 $y\in $ $\left( \bigcap _{j=1}^{\infty }
h_{x _{1}}^{-1}\cdots 
h_{x _{j}}^{-1}(J(G))\right) \setminus 
\hat{J}_{x }(f).$ 
Then, we have 
$(x ,y)\in (\Sigma _{m} \times \CCI ) 
\setminus \tilde{J}(f).$ 
Hence, there exists a 
neighborhood $U$
of $x $ in $\Sigma _{m} $ and a 
neighborhood $V$ of $y$ in $\CCI $ such that 
$U\times V\subset \tilde{F}(f).$ 
Then, there exists an $n\in \NN $ such that 
$\{ w\in \Sigma _{m}\mid w_{j}=x_{j}, j=1,\ldots ,n\} \subset U.$ 
Combining it with \cite[Lemma 2.4]{S1}, 
we obtain $\tilde{F}(f)\supset 
f^{n}(U\times V)\supset \Sigma _{m} \times 
\{ f_{x ,n}(y)\} .$ 
Moreover, 
since we have 
$f_{x ,n}(y)\in J(G)=
\pi _{\CCI }(\tilde{J}(f))$, 
we get that there exists an element 
$x '\in \Sigma _{m} $ such that 
$(x ', f_{x ,n}(y))\in 
\tilde{J}(f).$ However, it contradicts 
$(x ',f_{x ,n}(y))\in 
\Sigma _{m} \times \{ f_{x ,n}(y)\} 
\subset \tilde{F}(f).$ 
Hence, we obtain 
$\hat{J}_{x }(f)= 
\bigcap _{j=1}^{\infty }
h_{x_{1}}^{-1}\cdots 
h_{x_{j}}^{-1}(J(G)).$
\end{proof}
\begin{df}
Let $h_{1},\ldots ,h_{m}$ be polynomials and let 
$f: \Sigma _{m}\times \CCI \rightarrow \Sigma _{m}\times \CCI $ be 
the skew product associated with $(h_{1},\ldots ,h_{m}).$ 
For each $x\in \Sigma _{m}$, 
we set 
$K_{x}(f):=\{ y\in \CC \mid \{ f_{x,n}(y)\} _{n\in \NN } \mbox{ is bounded in } \CC \} $ 
and $A_{x}(f):= \{ y\in \CCI \mid f_{x,n}(y)\rightarrow \infty \mbox{ as } n\rightarrow \infty \} .$ 
\end{df} 
By using the method in \cite{B,Mi}, the following Lemma~\ref{jxkxlem} is easy to show and we omit the proof. 
\begin{lem}
\label{jxkxlem} 
Let $h_{1},\ldots ,h_{m}\in {\cal Y}$ and 
let $G:= \langle h_{1},\ldots ,h_{m}\rangle .$ 
Let $f:\Sigma _{m}\times \CCI 
\rightarrow \Sigma _{m}\times \CCI $ be the skew product associated 
with $(h_{1},\ldots ,h_{m}).$ 
Then, $\infty \in F(G)$ and for each $x\in 
\Sigma _{m},\ $ we have that 
$\infty \in F_{x}(f)$,  
$J_{x}(f)=\partial K_{x}(f) =\partial A_{x}(f)$,  
and $A_{x}(f)$ is the connected component of $F_{x}(f)$ containing $\infty .$ 
\end{lem}
\begin{lem}
\label{finfibconnlem}
Let $h_{1},\ldots ,h_{m}\in {\cal Y}$ 
and  
let $f:\Sigma _{m}\times \CCI \rightarrow 
\Sigma _{m}\times \CCI $ be the skew product 
associated with $(h_{1},\ldots ,h_{m}).$ 
Let $G=\langle h_{1},\ldots ,h_{m}\rangle .$ 
Then, 
the following {\em (1),(2),(3)} are equivalent.
{\em (1)} $G\in {\cal G}.$ 
{\em (2)} For each $x\in X$, $J_{x}(f)$ is 
connected.
{\em (3)} For each $x\in X$, $\hat{J}_{x}(f)$ is 
connected.
\end{lem}
\begin{proof} 
First, we show (1)$\Rightarrow $(2). 
Suppose that (1) holds.
Let 
$R>0$ be a number such that 
for each $x\in X$, 
$B:=\{ y\in \CCI 
\mid |y|>R\} \subset A_{x}(f)$ 
and $\overline{f_{x,1}(B)}\subset B.$ 
Then, for each $x\in X$, we have 
$A_{x}(f)=\bigcup _{n\in \NN }(f_{x,n})^{-1}(B)$ 
and $(f_{x,n})^{-1}(B)\subset 
(f_{x,n+1})^{-1}(B)$, 
for each $n\in \NN .$ Furthermore, since 
we assume (1), we see that 
for each $n\in \NN $, $(f_{x,n})^{-1}(B)$ 
is a simply connected domain, by the Riemann-Hurwitz formula (\cite{B,Mi}). 
Hence, for each $x\in X$, 
$A_{x}(f)$ is a simply connected domain.   
Since $\partial A_{x}(f)=J_{x}(f)$ for each $x\in X,$ 
we conclude that for each $x\in X$, 
$J_{x}(f)$ is connected. 
Hence, we have shown (1)$\Rightarrow $(2). 

 Next, we show 
 (2)$\Rightarrow $(3).  
 Suppose that (2) holds. 
 Let $z_{1}\in \hat{J}_{x}(f)$ and $z_{2}\in J_{x}(f)$ be two points. 
Let 
 $\{ x^{n}\} _{n\in \NN }$ be a sequence in $\Sigma _{m}$ 
 such that $x^{n}\rightarrow x$ as $n\rightarrow \infty $, and 
 such that 
 $d(z_{1},J_{x^{n}}(f))\rightarrow 0$ as 
 $n\rightarrow \infty .$ We may assume that 
 there exists a non-empty compact set $K$ in 
 $\CCI $ such that 
 $J_{x^{n}}(f)\rightarrow K$ as $n\rightarrow \infty $, 
 with respect to the Hausdorff topology in the 
 space of non-empty compact subsets of $\CCI .$ 
 Since we assume (2),  
 $K$ is connected. 
 By Remark~\ref{lowersemicontirem}, 
 we have $d(z_{2}, J_{x^{n}}(f))\rightarrow 0$ as 
 $n\rightarrow \infty .$ Hence, 
 $z_{i}\in K$ for each $i=1,2.$ 
 Therefore, denoting by $J$ the connected component of $\hat{J}_{x}(f)$ containing 
$K$, $z_{1}$ and $z_{2}$ belong to the same connected 
 component $J$ of $\hat{J}_{x}(f).$ Thus, 
 we have shown (2)$\Rightarrow $(3). 

 Next, we show (3)$\Rightarrow $(1).  
 Suppose that (3) holds.
 It is easy to see that 
 $A_{x}(f)\cap \hat{J}_{x}(f)=\emptyset $ for each 
 $x\in X.$ 
 Hence, $A_{x}(f)$ is a connected component of 
 $\CCI \setminus \hat{J}_{x}(f).$ 
 Since we assume (3), we have that 
 for each $x\in X$, 
 $A_{x}(f)$ is a simply connected domain. 
 Since $(f_{x,1})^{-1}(A_{g(x)}(f))=
 A_{x}(f)$ for each $x\in \Sigma _{m}$, the Riemann-Hurwitz formula implies that 
 for each $x\in X$, 
 there exists no critical point of 
 $f_{x,1}$ in $A_{x}(f)\cap \CC .$ 
 Therefore, we obtain (1).  
 Thus, we have shown (3)$\Rightarrow $(1). 
\end{proof}   
 
\begin{cor}\label{Lxconn} 
Let $G=\langle h_{1},\ldots ,h_{m}\rangle 
\in {\cal G}.$ Let $f:\Sigma _{m}\times \CCI 
\rightarrow \Sigma _{m}\times \CCI $ be the 
skew product associated with $(h_{1},\ldots ,h_{m}).$ 
Then, 
for each $x\in \Sigma _{m},$ the following sets 
$J_{x}(f), 
\hat{J}_{x}(f)$, and   
$\bigcap _{j=1}^{\infty }
h_{x_{1}}^{-1}\cdots h_{x_{j}}^{-1}(J(G))$ are connected. 
\end{cor}
\begin{proof}
From Lemma~\ref{fundfiblem} and Lemma~\ref{finfibconnlem}, 
the statement of the corollary easily follows. 
\end{proof}
\subsection{Dynamics of postcritically 
bounded polynomial semigroups}
We show a lemma on the dynamics of polynomial semigroups in ${\cal G}.$ 
\begin{lem}\label{biconn}
Let $G\in {\cal G}.$ Suppose that $J(G)$ is connected. 
Then, for any element $g\in G$, 
$g^{-1}(J(G))$ is connected.
\end{lem}
\begin{proof}
Let $g\in G.$ Since $G\in {\cal G}$, we have 
that $J(g)$ is a non-empty connected subset of $J(G).$ 
Let $J\in \mbox{Con}(g^{-1}(J(G)))$ be any element.  
By \cite{N} or \cite[Lemma 5.7.2]{B}, 
we have that $g(J)=J(G).$  Since $g^{-1}(J(g))=J(g)$, 
it follows that $J\cap J(g)\neq \emptyset .$ Hence $J(g)\subset J.$ 
Since this holds for any $J\in \mbox{Con}(g^{-1}(J(G)))$, 
$g^{-1}(J(G))$ is connected. 
\end{proof}
\begin{rem}
For further results on the dynamics of $G\in {\cal G}$, 
see \cite{SdpbpI,SdpbpII, SdpbpIII, S11, S10}. 
\end{rem}
\section{Proofs of results}
\label{Proofs}
In this section, we give the proofs of the main results in section~\ref{Main}. 
\subsection{Proofs of results in section~\ref{General}}
\label{ProofsGeneral}
In this subsection, 
we give the proofs of the results in section~\ref{General}.
We need some lemmas. 
\begin{df}
For each $j\in \{ 1,\ldots ,m\} $ and 
each $k\in \NN $, 
we set $(j)^{k}:=(j,j,\ldots ,j)\in \{ 1,\ldots ,m\} ^{k}.$  
\end{df}
\begin{lem}
\label{1disconlem}
Let $m\geq 2$ and 
let ${\frak L}=(L,(h_{1},\ldots ,h_{m}))$ be a backward self-similar 
system. 
Suppose 
that for each $j$ with 
$j\neq 1$, $h_{1}^{-1}(L)\cap 
h_{j}^{-1}(L)=\emptyset .$
For each $k$, let $C_{k}\in \mbox{{\em Con}}(|\Gamma _{k}|)$ 
 be the element containing 
 $(1)^{k}\in \{ 1,\ldots ,m\} ^{k}.$ 
 Then, we have the following.
 \begin{enumerate}
 \item \label{1disconlem1}
 For each $k\in \NN $, 
 $C_{k}=\{ (1)^{k}\} .$
  
\item \label{1disconlem2} For each $k\in \NN $, 
$\sharp (\mbox{{\em Con}}(|\Gamma _{k}|))
<\sharp (\mbox{{\em Con}}(|\Gamma _{k+1}|)).$
\item \label{1disconlem3}
$L$ has infinitely many connected components. 
\item \label{1disconlem4}
Let $x:=(1)^{\infty }\in \Sigma _{m}$ and 
let $x'\in \Sigma _{m}$ be an element with $x\neq x'.$ Then, 
for any $y\in L_{x}$ and $y'\in L_{x'}$, there exists no 
connected component $A$ of $L$ such that 
$y\in A$ and $y'\in A.$  
\end{enumerate}
\end{lem}
\begin{proof}
 We show statement \ref{1disconlem1} by induction on $k.$ 
 We have  $C_{1}=\{ 1\} .$ 
 Suppose $C_{k}=\{ (1)^{k}\} .$ 
 Let $w\in \{ 1,\ldots ,m\} ^{k+1}\cap C_{k+1}$ be 
 any element. 
 Since 
 $(\varphi _{k})_{\ast }(C_{k+1})=C_{k}$, 
 we have $\varphi _{k}(w)=(1)^{k}.$ 
 Hence, $w|k=(1)^{k}.$ 
 Since $
 h_{1}^{-1}(L)\cap h_{j}^{-1}(L)=\emptyset $ 
for each $j\neq 1$, 
we obtain $w=(1)^{k+1}.$ 
Hence, the induction is completed. Therefore, 
we have shown statement \ref{1disconlem1}. 

 Since both $(1)^{k+1}\in \{ 1,\ldots ,m\} ^{k+1}$ 
 and $(1)^{k}2\in \{ 1,\ldots ,m\} ^{k+1}$ are 
 mapped to $(1)^{k}$ by 
 $\varphi _{k}$, combining statement \ref{1disconlem1} and 
 Lemma~\ref{phionto}, 
 we obtain statement \ref{1disconlem2}. 
 For each $k\in \NN $, 
 we have 
\begin{equation}
\label{1disconlemeq1}
 L=\coprod _{C\in \mbox{Con}(|\Gamma _{k}|)}
 \ \bigcup _{w\in \{ 1,\ldots ,m\} ^{k}\cap C}
 h_{w}^{-1}(L).
\end{equation} 
 Hence, by statement \ref{1disconlem2}, we obtain 
 that $L$ has infinitely many connected components.

 We now show statement \ref{1disconlem4}. 
 Let $k_{0}:=\min \{ l\in \NN \mid x_{l}'\neq 1\} . $ 
 Then, by (\ref{1disconlemeq1}) and statement 
 \ref{1disconlem1}, we obtain that there exist compact sets 
 $B_{1}$ and $B_{2}$ such that 
 $B_{1}\cap B_{2}=\emptyset ,\ B_{1}\cup B_{2}=L,\ 
L_{x}\subset (h_{1}^{k_{0}})^{-1}(L)\subset  B_{1},$ and 
$L_{x'}\subset h_{x_{1}'}^{-1}\cdots h_{x_{k_{0}}'}^{-1}(L)\subset B_{2}.$ 
Hence, statement \ref{1disconlem4} holds.
\end{proof}
By an argument similar to that of the proof of Lemma~\ref{1disconlem}, 
we can prove the following.
\begin{lem}
\label{ss1disconlem}
Let $m\geq 2$ and 
let ${\frak L}=(L,(h_{1},\ldots ,h_{m}))$ be a forward self-similar 
system such that for each $j=1,\ldots ,m$, 
$h_{j}:L\rightarrow L$ is injective.  
Suppose 
that for each $j$ with 
$j\neq 1$, $h_{1}(L)\cap 
h_{j}(L)=\emptyset .$
For each $k$, let $C_{k}\in \mbox{{\em Con}}(|\Gamma _{k}|)$ 
 be the element containing 
 $(1)^{k}\in \{ 1,\ldots ,m\} ^{k}.$ 
 Then, all of the statements \ref{1disconlem1}--\ref{1disconlem4} in Lemma~\ref{1disconlem}
 hold. 
\end{lem}

To prove Theorem~\ref{precompojth}, we need the 
following lemma.

\begin{lem}
\label{ldisj}
Under the assumptions of Theorem~\ref{precompojth},
let $M_{1},\ldots ,M_{r}$ be 
mutually disjoint non-empty compact subsets of 
$L$ with $L=\bigcup _{i=1}^{r}M_{i}.$ 
Then there exists a number 
$l_{0}\in \NN $ such that 
for each $x\in \Sigma _{m}$ 
and each $l\in \NN $ with $l\geq l_{0}$,  
there exists a number $i=i(x,l)\in \{ 1,\ldots ,r\} $ 
with $h_{x|l}^{-1}(L)\subset M_{i}.$  
\end{lem}
\begin{proof}
Suppose that the statement is not true. 
Then for each $n\in \NN $, 
there exist an element $w^{n}\in \Sigma _{m}$,  
 an $l(n)> n$, and elements
$i_{1,n},i_{2,n}\in \{ 1,\ldots ,r\} $ 
with $M_{i_{1,n}}\neq M_{i_{2,n}}$,  
  such that 
$(h_{w^{n}|l(n)})^{-1}(L)\cap 
M_{i}\neq \emptyset $, 
for each $i=i_{1,n},i_{2,n}.$ 
Since $\Sigma _{m}$ is compact, we may 
assume that there exists an element 
$w\in \Sigma _{m}$ 
such that  
for each $n\in \NN $,  
$w^{n}|l(n)=(w|n) \alpha _{n}$ 
for some $\alpha _{n}\in \Sigma _{m}^{\ast }.$   
 
 Then, we have 
 $h_{w|n}^{-1}h_{\alpha _{n}}^{-1}(L)
 \cap M_{i}\neq \emptyset $, 
 for each $i=i_{1,n},i_{2,n}.$ 
 Hence, $h_{w|n}^{-1}(L)
 \cap M_{i}\neq \emptyset $, 
 for each $i=i_{1,n},i_{2,n}.$ 
Since $h_{w|n}^{-1}(L)\rightarrow 
L_{w}$ as $n\rightarrow \infty $ with respect 
to the Hausdorff topology  
and  
$L_{w}$ is connected (the assumption), 
we obtain a contradiction.   
\end{proof}
By an argument similar to that of the proof of Lemma~\ref{ssldisj}, 
we can prove the following. 
\begin{lem}
\label{ssldisj}
Under the assumptions of Theorem~\ref{ssprecompojth},
let $M_{1},\ldots ,M_{r}$ be 
mutually disjoint non-empty compact subsets of 
$L$ with $L=\bigcup _{i=1}^{r}M_{i}.$ 
Then, there exists a number 
$l_{0}\in \NN $ such that 
for each $x\in \Sigma _{m}$ and each $l\in \NN $  
with $l\geq l_{0}$,  
there exists a number $i=i(x,l)\in \{ 1,\ldots ,r\} $ 
with $h_{\overline{x|l}}(L)\subset M_{i}.$  
\end{lem}

We now demonstrate Theorem~\ref{precompojth}.

\noindent {\bf Proof of Theorem~\ref{precompojth}:}
Step 1: 
First, we show the following:\\ 
 Claim 1: Let $B=(B_{k})_{k}\in  \varprojlim \mbox{Con}(|\Gamma _{k}|)$ 
where $B_{k}\in \mbox{Con}(|\Gamma _{k}|)$ 
and $(\varphi _{k})_{\ast }(B_{k+1})=B_{k}$ for each $k.$ 
Take a point $x\in \Sigma _{m}$ 
 such that $x|k\in 
 B_{k}$ for each $k.$ Take an element 
 $C_{x}\in \mbox{Con}(L)$ such that 
 $L_{x}\subset C_{x}.$ Then, 
 $C_{x}$ does not depend on the 
 choice of $x\in \Sigma _{m}$ 
 such that $x|k\in 
 B_{k}$ for each $k.$ 
Hence, the map $\Phi : B\mapsto 
C_{x}$ is well-defined. 

 To show Claim 1, suppose 
 that there exist $x\in 
 \Sigma _{m}$ and $y\in \Sigma _{m}$ 
 such that $x|k,y|k\in B_{k}$ for each $k\in \NN $ 
 and such that there exist mutually different 
 connected components $J_{1}$ and $J_{2}$ of 
 $L$ with $L_{x}\subset J_{1}$ and 
 $L_{y}\subset J_{2}.$ 
 By the ``Cut Wire Theorem'' in \cite{N}, 
 there exist mutually disjoint 
 compact subsets $M_{1}$ and 
 $M_{2}$ of $L$ such that 
 $J_{i}\subset M_{i}$ for each $i=1,2.$ 
We apply Lemma~\ref{ldisj} to 
the disjoint union $L=M_{1}\cup M_{2}$ and 
let $l_{0}$ be the number in the lemma. 
Then, we have 
$h_{x|l_{0}}^{-1}(L)\subset M_{1}$,  
$h_{y|l_{0}}^{-1}(L)\subset M_{2}$,  
and 
$L=\bigcup _{|w|=l_{0}}h_{w}^{-1}(L)
=\coprod _{i=1}^{2}
\bigcup _{h_{w}^{-1}(L)\subset M_{i}, |w|=l_{0}}
h_{w}^{-1}(L).$ This implies that 
$x|l_{0}$ and $y|l_{0}$ do not 
belong to the same connected component 
of $|\Gamma _{l_{0}}|.$ This is a contradiction. 
Hence, we have shown Claim 1. 

Step 2: Next, we show the following:\\ 
Claim 2: $\Phi : \varprojlim \mbox{Con}(|\Gamma _{k}|)
\rightarrow \mbox{Con}(L)$ 
is bijective.

 To show Claim 2, since 
$ L=\bigcup _{j=1}^{m}h_{j}^{-1}(L)  $, 
we have $L=\bigcup _{x\in \Sigma _{m}}L_{x}.$ 
Hence, 
$\Phi $ is surjective.
 To show that $\Phi $ is injective, let $B=(B_{k})$ and 
 $B'=(B_{k}')$ be distinct elements in $ \varprojlim \mbox{Con}(|\Gamma _{k}|),$ 
 let $x\in \Sigma _{m}$ be such that  
 $x|k\in B_{k}$ for each $k\in \NN $, and 
 let $y\in \Sigma _{m}$ be such that  
 $y|k\in B_{k}'$ for each $k\in \NN .$ 
 Then, there exists a $k\in \NN $ with 
 $B_{k}\neq B_{k}'.$ 
 Combining this with 
$L=\coprod _{C\in \mbox{Con}(|\Gamma _{k}|)}\ 
 \bigcup _{w\in \Sigma _{m}^{\ast }\cap C, |w|=k}
 h_{w}^{-1}(L),$
 which follows from $L=\bigcup _{w\in \Sigma _{m}^{\ast }, |w|=k}
 h_{w}^{-1}(L)$, 
 we obtain that there exist two 
 compact subsets $K_{1}$ and $K_{2}$ of $L$ such that 
 $L=K_{1}\coprod K_{2}$, 
 $L_{x}\subset h_{x|k}^{-1}(L)\subset K_{1}$,  and 
 $L_{y}\subset h_{y|k}^{-1}(L)\subset K_{2}.$ 
 Hence, $\Phi (B)\neq \Phi (B').$ 
 Therefore, $\Phi $ is injective.

 Step 3: We now show statement \ref{precompojth2}.
 Since $L=\bigcup _{j=1}^{m}h_{j}^{-1}(L)$, it is 
 easy to see that 
 if $L$ is connected, then 
 $|\Gamma _{1}|$ is connected. 
 Conversely, suppose that 
 $|\Gamma _{1}|$ is connected. 
 Then, by Lemma~\ref{1conall}, 
 we obtain that for each $k\in \NN $,\ 
 $|\Gamma _{k}|$ is connected.
From statement \ref{precompojth1}, it follows that 
$L$ is connected.  
Hence, we have shown statement \ref{precompojth2}.  

Step 4: 
 Statement \ref{precompojth3} follows from 
 statement \ref{precompojth1} and Lemma~\ref{phionto}.
Statement \ref{precompojth6} and \ref{precompojth7} easily follow 
from statement \ref{precompojth3}. 

 Step 5: 
 We now show statement \ref{precompojth4}. 
 If $m=2$ and $L$ is disconnected, then 
 by statement \ref{precompojth2}, 
 we have $h_{1}^{-1}(L)\cap h_{2}^{-1}(L)=\emptyset .$ 
 Combining this with statement \ref{precompojth1}, 
 we obtain 
 $\mbox{Con}(L)\cong \{ 1,2\} ^{\NN }.$ 

 Step 6: 
 We now show statement \ref{precompojth5}.
 Suppose that $m=3$ and $L $ is disconnected. 
 By statement \ref{precompojth2}, 
 we may assume  
 $
 h_{1}^{-1}(L)\cap 
 h_{2}^{-1}(L)=h_{1}^{-1}(L)\cap 
 h_{3}^{-1}(L)=\emptyset .
 $  
 By Lemma~\ref{1disconlem}, we obtain 
 that $L$ has infinitely many connected components and that 
 $L_{(1)^{\infty }}$ is a connected component of $L.$ 
 
Thus we have completed the proof of Theorem~\ref{precompojth}.  
 \qed 

\

We now prove Theorem~\ref{ssprecompojth}. 

\noindent {\bf Proof of Theorem~\ref{ssprecompojth}:} 
The statements of the theorem easily follow from the argument of 
the proof of Theorem~\ref{precompojth}, Lemma~\ref{ssldisj}, and 
Lemma~\ref{ss1disconlem}.
\qed  

\ 

In order to prove Theorem~\ref{preicmainthm}, we need the following 
notations and lemmas. 

\begin{lem}
\label{preicmainthmpflem1}
Under the assumptions of Theorem~\ref{preicmainthm} or Theorem~\ref{sspreicmainthm}, 
let $k\in \NN .$ Then, for any simplex $s$ of $N_{k}$ with $(1)^{k}\in s$, 
the dimension $\dim s$ of $s$ is less than or equal to $1.$  
\end{lem}
\begin{proof}
We will show the conclusion of our lemma for a backward self-similar 
system ${\frak L}=(L,(h_{1},\ldots ,h_{m}))$ satisfying the assumptions of 
Theorem~\ref{preicmainthm} (using an argument similar to the below, 
we can show the conclusion of our lemma for a forward self-similar 
system ${\frak L}$ satisfying the assumptions of Theorem~\ref{sspreicmainthm}).  
We will show the conclusion of our lemma by induction on $k\in \NN .$ 
If $k=1$, then, assumption \ref{preicmainthmc4} of Theorem~\ref{preicmainthm} 
implies that for any simplex $s$ of $N_{1}$ with $1\in s$, we have $\dim s\leq 1.$ 
Let $l\in \NN $ and we now suppose that for any 
simplex $s$ of $N_{l}$ with $(1)^{l}\in s$, we have $\dim s\leq 1.$ 
Then, Lemma~\ref{nkwlem} implies that 
for any simplex $s$ of $N_{l+1,1}$ with $(1)^{l+1}\in s$, we have $\dim s\leq 1.$ 
Moreover, by assumption \ref{preicmainthmc2} of Theorem~\ref{preicmainthm}, 
we have $(h_{1}^{r})^{-1}(L)\cap (\bigcup _{i\neq 1}h_{i}^{-1}(L))=\emptyset $ for each 
$r\geq 2.$ Hence, it follows that for any $i\in \{ 1,\ldots ,m\} $ with $i\neq 1$ and 
any $w\in \Sigma _{m}^{\ast }$ with $|w|=l$,    
$\{ (1)^{l+1}, iw\} $ is not a simplex of $N_{l+1}.$ Therefore, 
for any simplex $s$ of $N_{l+1}$ with $(1)^{l+1}\in s$, we have $\dim s\leq 1.$ 
Thus, the induction is completed.   
\end{proof}
\begin{df}
Let $S$ be a simplicial complex and 
let $\tau =(v_{1},v_{2})(v_{2},v_{3})\cdots (v_{n-1},v_{n})$ be 
 an edge path of $S.$ We denote by $|\tau |$ the curve in 
 $|S|$ which is induced by $\tau $ in the way as in \cite[p.136]{Sp}. 

\end{df}
\begin{df}
Let ${\frak L}$ be a  forward or backward self-similar system, 
let $k\in \NN $, and let $w\in \Sigma _{m}^{\ast }$. Then for any 
edge path $\tau =(v_{1},v_{2})(v_{2},v_{3})\cdots (v_{n-1},v_{n})$   
of $N_{k}$, we denote by $w_{\ast }(\tau )$ the 
edge path $(wv_{1},wv_{2})(wv_{2},wv_{3})\cdots (wv_{n-1},wv_{n})$ of 
$N_{k+|w|}.$   
\end{df}
\begin{lem}
\label{preicmainthmpflem2}
Under the assumptions of Theorem~\ref{preicmainthm} or Theorem~\ref{sspreicmainthm}, 
let $\tau $ be the closed edge path $(1,j_{2})(j_{2},j_{3})(j_{3},1)$ of $N_{1}.$ 
Moreover, let $\gamma \in H_{1}(|N_{1}|;R)$ be the element induced by 
the closed curve $|\tau | $ in $|N_{1}|.$ Then, 
for each $k\in \NN $, the element $((1)^{k})_{\ast }(\gamma )\in 
H_{1}(|N_{k+1}|;R)$ is not zero.  
\end{lem}
\begin{proof}
For each $k\in \NN $, let $M_{k}$ be the 
unique full subcomplex of $N_{k}$ whose vertex set is equal to 
$\{ 1,\ldots ,m\} ^{k}\setminus \{ (1)^{k}\} .$ 
Moreover, let $P_{k}$ be the set of all 
$1$-simplexes $e$ of $N_{k+1}$ such that 
$(1)^{k+1}\in e$, $(1)^{k}j_{2}\not\in e$, and $(1)^{k}j_{3}\not\in e.$ Furthermore, 
let $Q_{k}=\bigcup _{e\in P_{k}}e.$ Note that $Q_{k}$ is a subcomplex of $N_{k+1}.$ 
Lemma~\ref{preicmainthmpflem1} implies that 
for each $k\in \NN $, 
$|N_{k+1}|=|((1)^{k})_{\ast }(\tau )|\cup |Q_{k}|\cup |M_{k+1}|.$ Moreover,  
$\left( |((1)^{k})_{\ast }(\tau )|\cup |Q_{k}|\right)\cap |M_{k+1}|=
|((1)^{k}j_{2},(1)^{k}j_{3})|\cup \bigcup _{e\in P_{k}}\{ e_{0}\} $, 
where for each $e\in P_{k}$, $e_{0}$ denotes the 
vertex of $e$ which is not equal to $(1)^{k+1}.$ In particular, 
each connected component of 
$ \left( |((1)^{k})_{\ast }(\tau )|\cup |Q_{k}|\right)\cap |M_{k+1}|$ is contractible. 
Using the Mayer-Vietoris sequence of $\{ |((1)^{k})_{\ast }\tau |\cup |Q_{k}|, |M_{k+1}|\}$, we obtain the following exact sequence: 
\begin{multline}
0 =  H_{1}(\left( |((1)^{k})_{\ast }(\tau )|\cup |Q_{k}|\right)\cap |M_{k+1}|;R)  \rightarrow \\  
   H_{1}(|((1)^{k})_{\ast }(\tau )|\cup |Q_{k}|;R)\oplus H_{1}(|M_{k+1}|;R)     \rightarrow 
H_{1}(|N_{k+1}|;R).
\end{multline}   
Let $u _{1}:|((1)^{k})_{\ast }(\tau )|\rightarrow |((1)^{k})_{\ast }(\tau )|\cup |Q_{k}|$, 
$u _{2}: |((1)^{k})_{\ast }(\tau )|\cup |Q_{k}|\rightarrow |N_{k+1}|$, 
and $u _{3}:|((1)^{k})_{\ast }(\tau )|\rightarrow |N_{k+1}|$ be the inclusion maps. 
Then, $u _{3}=u _{2}\circ u _{1}.$ 
Moreover, $(u _{1})_{\ast }: 
H_{1}(|((1)^{k})_{\ast }(\tau )|;R)\rightarrow H_{1}(|((1)^{k})_{\ast }(\tau )|\cup |Q_{k}|;R)$ 
is an isomorphism. Furthermore, 
$((1)^{k})_{\ast }(\gamma )=(u _{3})_{\ast }(a)$ in 
$H_{1}(|N_{k+1}|;R)$, where $a$ is a generator in $H_{1}(|((1)^{k})_{\ast }(\tau )|;R).$ 
From these arguments, it follows that the element $((1)^{k})_{\ast }(\gamma )\in H_{1}(|N_{k+1}|;R)$ is 
not zero. Thus, we have proved the lemma.  
\end{proof}
\begin{lem}
\label{preicmainthmpflem3}
Under the assumptions of Theorem~\ref{preicmainthm} or Theorem~\ref{sspreicmainthm}, 
we have that for each 
$k\in \NN $, 
$\dim _{R}\check{H}^{1}({\frak L};R)_{k}=\dim _{R}\check{H}_{1}({\frak L};R)_{k}<
\dim _{R}\check{H}^{1}({\frak L};R)_{k+1}=\dim _{R}\check{H}_{1}({\frak L};R)_{k+1}.$ 
\end{lem}
\begin{proof}
We use the notation in Lemma~\ref{preicmainthmpflem2}. 
By Lemma~\ref{preicmainthmpflem2}, 
we have that for each $k\in \NN $, 
$((1)^{k})_{\ast }(\gamma )\in H_{1}(|N_{k+1}|;R)$ is not zero. 
Moreover, by Lemma~\ref{wmaplem}, 
we have that for each $k\in \NN $, 
$(\varphi _{k})_{\ast }(((1)^{k})_{\ast }(\gamma ))=0.$ 
Hence $(\varphi _{k})_{\ast }: H_{1}(|N_{k+1}|;R)\rightarrow 
H_{1}(|N_{k}|;R)$ is not a monomorphism. 
Furthermore, by assumption \ref{preicmainthmc1} of Theorem~\ref{preicmainthm} and Theorem~\ref{sspreicmainthm}
and Lemma~\ref{icfundlem}-\ref{icfundlem3-01},    
we have that $(\varphi _{k})_{\ast }: H_{1}(|N_{k+1}|;R)\rightarrow 
H_{1}(|N_{k}|;R)$ is an epimorphism. 
It follows that for each $k\in \NN $, 
$\dim _{R}H_{1}(|N_{k}|;R)<\dim _{R}H_{1}(|N_{k+1}|;R).$ 
We are done. 
\end{proof}
We now prove Theorem~\ref{preicmainthm} and Theorem~\ref{sspreicmainthm}.

\noindent {\bf Proof of Theorem~\ref{preicmainthm} and Theorem~\ref{sspreicmainthm}:} 
By the assumption \ref{preicmainthmc1} of Theorem~\ref{preicmainthm} and Theorem~\ref{sspreicmainthm} and 
Lemma~\ref{icfundlem}-\ref{icfundlem3},  
the projection map $\mu _{k,1}:\check{H}^{1}({\frak L};R)_{k}\rightarrow \check{H}^{1}({\frak L};R)$ 
is injective for each $k\in \NN .$  
Combining it and Lemma~\ref{preicmainthmpflem3}, 
we obtain that $\dim _{R}\check{H}^{1}({\frak L};R)=\infty .$ Thus, 
we have proved Theorem~\ref{preicmainthm} and Theorem~\ref{sspreicmainthm}.
\qed 

\ 

We now prove Corollary~\ref{c:sspreic}.\\ 
\noindent {\bf Proof of Corollary~\ref{c:sspreic}:} 
Since $|N_{1}|$ is connected and $L_{x}$ is connected for each $x\in \Sigma _{m}$, 
Theorem~\ref{ssprecompojth} implies that $L$ is connected. Thus for each $w\in \Sigma _{m}^{\ast }$, 
$h_{\overline{w}}(L)$ is connected. Combining it with Lemma~\ref{icfundlem}-\ref{h1inj} and Theorem~\ref{sspreicmainthm}, 
we obtain that the statement our corollary holds. 
\qed 

\subsection{Proofs of results in section~\ref{Application}}
\label{Proofsappl}
In this subsection, we give the proofs of the results in subsection~\ref{Application}. 

We now prove Theorem~\ref{compojth}. 

\noindent {\bf Proof of Theorem~\ref{compojth}
:}
From Theorem~\ref{precompojth} and 
Corollary~\ref{Lxconn}, the statement of the theorem follows.
\qed \\ 
\ 

We now prove Theorem~\ref{icmainthm}.\\  
\noindent {\bf Proof of Theorem~\ref{icmainthm}:} 
By Theorem~\ref{compojth}, 
$J(G)$ is connected. 
Combining it with Lemma~\ref{biconn}, 
we obtain that for each $g\in G$, 
$g^{-1}(J(G))$ is connected. 
By Lemma~\ref{icfundlem}-\ref{h1inj}, 
it follows that 
$\Psi : \check{H}^{1}({\frak L};R)\rightarrow \check{H}^{1}(J(G),R)$ is a monomorphism. 
Moreover, by Theorem~\ref{preicmainthm}, 
we obtain 
$\dim _{R}\check{H}^{1}({\frak L};R)$ $=\infty. $ 
Hence, $\dim _{R}\Psi (\check{H}^{1}({\frak L};R))=\infty .$ 
Therefore,\ $\dim _{R}\check{H}^{1}(J(G);R)=\infty  .$ 
By the Alexander duality theorem (see \cite[p.296]{Sp}), 
we have $\check{H}^{1}(J(G);R)\cong 
\tilde{H}_{0} (\CCI \setminus J(G); R),$ 
where $\tilde{H}_{0}$ denotes the 
$0$-th reduced homology. 
Hence, $F(G)=\CCI \setminus J(G)$ has infinitely many 
connected components. 
\qed 

\ 

\begin{figure}[htbp]
\caption{The Julia set of $G$ in Proposition~\ref{pbpch1inftyexprop}.}
\ \ \ \ \ \ \ \ \ \ \ \ \ \ \ \ \ \ \ \ \ \ \ \ \ \ \ \ \ \ \ 
\ \ \ \ \ \ \ \ \ \ \ \ \ \ \ \ \ \ \ \ \  \ 
\includegraphics[width=2.5cm,width=2.5cm]{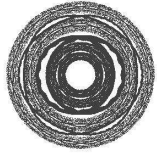}
\label{fig:h1infty}
\end{figure}

We now prove Proposition~\ref{pbpch1inftyexprop}. \\ 
\noindent {\bf Proof of Proposition~\ref{pbpch1inftyexprop}:} 
Let  $a\in \RR $ with $1<a\leq 5$. Let  
 $h_{1}(z)=\frac{1}{a^{2}}z^{3}$ and 
 $h_{2}(z)=z^{2}.$ 
Then $J(h_{1})=\{ z\in \CC \mid |z|=a\}$, 
$J(h_{2})=\{ z\in \CC \mid |z|=1\} $, 
$h_{1}^{-1}(J(h_{2}))=\{ z\in \CC \mid |z|=a^{2/3}\}$, 
and $h_{2}^{-1}(J(h_{1}))=\{ z\in \CC \mid |z|=a^{1/2}\}.$  
Let $c_{1}:=
 (a^{\frac{2}{3}}-a^{\frac{1}{2}})/2.$ 
 Let $g_{3}$ be a 
 polynomial such that 
 $J(g_{3})=\{ z\in \CC \mid  |z-c_{1}|=a^{\frac{2}{3}}-c_{1}\} $ and 
 let $g_{4}$ be a polynomial such that 
 $J(g_{4})=\{ z\in \CC \mid |z+c_{1}|=a^{\frac{2}{3}}-c_{1}\} .$ 
 Take a sufficiently large $n\in \NN $ 
  and 
 let $h_{3}=g_{3}^{n}$ and 
 $h_{4}=g_{4}^{n}.$ 
 Let $G=\langle h_{1},h_{2},h_{3},h_{4}\rangle $ and let 
$K:=\{ z\in \CC \mid 1\leq |z|\leq a\} .$ 
Then, taking a sufficiently large $n$, we have 
$\bigcup _{i=1}^{4}h_{i}^{-1}(K)\subset K.$   
Therefore, by \cite[Corollary 3.2]{HM}, $J(G)\subset K.$ 
(For the figure of $J(G)$, see Figure~\ref{fig:h1infty}.)
Moreover, we can show that $G\in {\cal G}$,  
 the set of all $1$-simplexes of 
 $N_{1}$ is equal to 
 $\{ \{ 1,3\},\ \{ 1,4\},\ \{ 2, 3\} ,\ 
 \{ 2,4\},\ \{ 3,4\} \} $, and 
 there exists no $r$-simplex $S$ of 
 $N_{1}$ for each $r\geq 2.$ 
 Taking a sufficiently large $n$ again, it is easy to show that 
 ${\frak L}=(J(G),(h_{1},h_{2},h_{3},h_{4}))$ 
satisfies all of the conditions 1,...,4 in the assumptions 
 of Theorem~\ref{preicmainthm}.
From Theorem~\ref{icmainthm}, it follows that 
 $\dim _{R}\check{H}^{1}(J(G),(h_{1},\ldots ,h_{4});R)=
 \dim _{R}\Psi (\check{H}^{1}(J(G),(h_{1},\ldots ,h_{4});R))=\infty $ and 
 $F(G)$ has infinitely many connected components. Thus we have completed the proof.
\qed 
\subsection{Proofs of results in section~\ref{Postunbranched}}
\label{ProofsPost}
In this subsection, we prove the results in section~\ref{Postunbranched}. 
We need some lemmas.
\begin{lem}
\label{pufundlem1}
Let ${\frak L}=(L,(h_{1},\ldots ,h_{m}))$ be a forward or backward self-similar system. 
Suppose that ${\frak L}$ is postunbranched. 
Let $r\in \NN .$ Then, for each $r$-simplex $e$ of $N_{1}$, 
there exists a unique $r$-simplex $e_{k+1}$ of $N_{k+1}$ such that 
$\varphi _{k+1,1}(e_{k+1})=e.$  

\end{lem}
\begin{proof}
We will show the conclusion of our lemma when ${\frak L} $ is a backward self-similar 
system (we can show the conclusion of our lemma when ${\frak L}$ is a forward self-similar 
system by using an argument similar to the below). 
The existence of $e_{k+1}$ follows from Lemma~\ref{phionto}. 
We now prove the uniqueness. \\
Case 1: $r=1$. Let $e=\{ i_{1},j_{1}\} $ be a $1$-simplex of $N_{1}.$ 
Then $C_{i_{1},j_{1}}=h_{i_{1}}^{-1}(L)\cap h_{j_{1}}^{-1}(L)\neq \emptyset .$ 
Since ${\frak L}$ is postunbranched, 
there exists a unique $x\in \Sigma _{m}$ such that 
$h_{i_{1}}(C_{i_{1},j_{1}})\subset L_{x}$ and 
such that for each $x'\in \Sigma _{m}$ with $x'\neq x$, 
$h_{i_{1}}(C_{i_{1},j_{1}})\cap L_{x'}=\emptyset .$  
Let $e_{k+1}=\{ (i_{1},\ldots ,i_{k+1}), (j_{1},\ldots ,j_{k+1})\} $ 
be a $1$-simplex of $N_{k+1}$ such that 
$\varphi _{k+1,1}(e_{k+1})=e.$ We will show that 
$(i_{2},\ldots ,i_{k+1})$ and $(j_{2},\ldots ,j_{k+1})$ are uniquely 
determined by the element $(i_{1},j_{1}).$  
Since $e_{k+1}$ is a $1$-simplex of $N_{k+1}$, 
we have 
$h_{i_{1}}^{-1}\cdots h_{i_{k+1}}^{-1}(L)
\cap h_{j_{1}}^{-1}\cdots h_{j_{k+1}}^{-1}(L)\neq \emptyset .$ 
Let $z\in h_{i_{1}}^{-1}\cdots h_{i_{k+1}}^{-1}(L)
\cap h_{j_{1}}^{-1}\cdots h_{j_{k+1}}^{-1}(L)$ be a point. 
Then 
\begin{equation}
\label{pufundlem1eq1}
h_{i_{1}}(z)\in h_{i_{2}}^{-1}\cdots h_{i_{k+1}}^{-1}(L).
\end{equation}
Moreover, since $z\in h_{i_{1}}^{-1}\cdots h_{i_{k+1}}^{-1}(L)
\cap h_{j_{1}}^{-1}\cdots h_{j_{k+1}}^{-1}(L)\subset C_{i_{1},j_{1}}$, 
we have $h_{i_{1}}(z)\in h_{i_{1}}(C_{i_{1},j_{1}})$ and for each 
$x'\in \Sigma _{m}$ with $x'\neq x$, $h_{i_{1}}(z)\not\in L_{x'}.$ 
Furthermore, since $L=\bigcup _{y\in \Sigma _{m}}L_{y}$, (\ref{pufundlem1eq1}) implies that 
there exists an element $y=(y_{1},y_{2},\ldots )\in \Sigma _{m}$ such that 
$h_{i_{1}}(z)\in h_{i_{2}}^{-1}\cdots h_{i_{k+1}}^{-1}(L_{y}).$ Let 
$y'=(i_{2},i_{3},\ldots ,i_{k+1},y_{1},y_{2},\ldots )\in \Sigma _{m}.$ Then 
$h_{i_{1}}(z)\in L_{y'}.$ From the above arguments, it follows that 
$y'=x.$ Therefore, $(i_{2},\ldots ,i_{k+1})=(x_{1},\ldots ,x_{k}).$ 
Thus, $(i_{2},\ldots ,i_{k+1})$ is uniquely determined by 
$(i_{1},j_{1}).$ Similarly, we can show that 
$(j_{2},\ldots ,j_{k+1})$ is uniquely determined by 
$(i_{1},j_{1}).$ Hence, $e_{k+1}$ is uniquely determined by  $e.$ \\ 
Case 2: $r\geq 2.$ The uniqueness immediately follows from Case 1. 

 Thus, we have proved Lemma~\ref{pufundlem1}. 
\end{proof}
\begin{df}
\label{tildephidf}
Let ${\frak L}=(L,(h_{1},\ldots ,h_{m}))$ be a forward or backward self-similar 
system. For each $k\in \NN  $, 
we denote by $S_{k}$ (or $S_{k}({\frak L})$) the 
CW complex $|N_{k}|/|\bigcup _{j=1}^{m}N_{k,j}|.$ 
Furthermore, 
we denote by 
$p_{k}:(|N_{k}|,|\bigcup _{j=1}^{m}N_{k,j}|)\rightarrow (S_{k},\ast )$
the canonical projection. 
Moreover, for all $l,k\in \NN $ with $l>k$, 
we denote by 
$\tilde{\varphi }_{l,k}: S_{l}\rightarrow 
S_{k} $ the cellular map 
 such that 
the following commutes. 
\begin{equation}
\label{tildephidfeq}
\begin{CD}
(|N_{l}|,|\bigcup _{j=1}^{m}N_{l,j}|)@>{p_{l}}>>(S_{l},\ast )\\ 
@V{|\varphi _{l,k}|}VV
@VV{\tilde{\varphi }_{l,k}}V\\ 
(|N_{k}|, |\bigcup _{j=1}^{m}N_{k,j}| )@>>{ p_{k}}>(S_{k} ,\ast )
\end{CD}
\end{equation}
\end{df}
\begin{lem}
\label{pufundlem2}
Let ${\frak L}=(L,(h_{1},\ldots ,h_{m}))$ be a forward or  
backward self-similar system. Suppose that ${\frak L}$ is postunbranched. 
Let $R$ be a $\ZZ $ module. 
Then, we have the following.
\begin{enumerate}
\item \label{pufundlem2-1}
For each $k\in \NN $, the cellular map 
$\tilde{\varphi }_{k+1,1}:(S_{k+1},\ast )
\rightarrow (S_{1} ,\ast )$ is a cellular 
isomorphism and a homeomorphism. In particular, $\tilde{\varphi }_{k+1,1}$ induces 
 isomorphisms on homology and cohomology groups with coefficients $R$.  
\item \label{pufundlem2-2}
For each $k\in \NN $, 
$|\varphi _{k+1,1}|:(|N_{k+1}|,|\bigcup _{j=1}^{m}N_{k+1,j}|)\rightarrow 
(|N_{1}|,\{ 1,\ldots ,m\} )$ induces isomorphisms on homology and cohomology groups with 
coefficient $R.$ 
\end{enumerate}
\end{lem}
\begin{proof}
From Lemma~\ref{pufundlem1}, statement \ref{pufundlem2-1} follows. 
Since $p_{k}$ induces isomorphisms on homology and cohomology groups, 
statement \ref{pufundlem2-2} follows from statement \ref{pufundlem2-1}. 
Thus, we have proved Lemma~\ref{pufundlem2}. 
\end{proof}
\begin{lem}
\label{pufundlem3}
Let ${\frak L}=(L,(h_{1},\ldots ,h_{m}))$ be a forward or  
backward self-similar system. Suppose that 
${\frak L}$ is postunbranched. Let $R$ be a $\ZZ $ module. 
Let $r\in \NN $ and $k\in \NN .$ Then, the connecting 
homomorphism 
$\partial _{\ast }: H_{r+1}(|N_{k}|,|\bigcup _{j=1}^{m}N_{k,j}|;R)\rightarrow 
H_{r}(|\bigcup _{j=1}^{m}N_{k,j}|;R)$ of the homology sequence of the pair 
$(|N_{k}|,|\bigcup _{j=1}^{m}N_{k,j}|)$ is the zero map.

\end{lem}
\begin{proof}
For each $k\in \NN $, let 
$\alpha _{k}:(N_{k},\ast )\rightarrow (N_{k},\bigcup _{j=1}^{m}N_{k,j})$ be 
the canonical embedding. 
Moreover, for each $k\in \NN $, let 
$\gamma _{k}:(|N_{k}|,\ast ) \rightarrow (S_{k},\ast )$ be the canonical projection. 
Then, for each $k\in \NN $, 
$(p_{k})_{\ast }: H_{r+1}(|N_{k}|,|\bigcup _{j=1}^{m}N_{k,j}|;R)\rightarrow 
H_{r+1}(|S_{k}|,\ast ;R)$ is an isomorphism,  
and the following diagram commutes. 
\begin{equation}
\label{pufundlem3eq1}
\begin{CD}
H_{r+1}(|N_{k}|,\ast ;R) @>{(\gamma _{k})_{\ast }}>> H_{r+1}(S_{k},\ast ;R)\\
 @V{(\alpha _{k})_{\ast }}VV
  @VV{Id}V\\
H_{r+1}(|N_{k}|, |\bigcup _{j=1}^{m}N_{k,j}|;R) @>>{(p_{k})_{\ast }}> H_{r+1}(S_{k},\ast ;R) 
\end{CD}
\end{equation}
Hence, we have only to prove that for each $k>1 $, 
$(\gamma _{k})_{\ast }: H_{r+1}(|N_{k}|, \ast ;R)\rightarrow H_{r+1}(S_{k},\ast ;R)$ is 
an epimorphism (if $k=1$, then it is easy to see that Im $\partial _{\ast }=0$). In order to do that,  
let $a=\sum _{i=1}^{t}a_{i}d_{i}\in C_{r+1}(S_{k},\ast ;R)$ be a cycle, 
where for each $i$, $a_{i}\in R$ and $d_{i}$ is an oriented ($r+1$)-cell. 
For each $i$, let 
$d_{i}':=\tilde{\varphi }_{k,1}(d_{i}).$ Then, by Lemma~\ref{pufundlem2}-\ref{pufundlem2-1}, 
 $d_{i}'$ is an ($r+1$)-cell of $S_{1}.$ Let $d_{i}''$ be an oriented ($r+1$)-cell of $|N_{1}|$ such that 
$\gamma _{1}(d_{i}'')=d_{i}'.$ Let $e_{i}''$ be the oriented ($r+1$)-simplex of $N_{1}$ which induces 
$d_{i}''.$ Then, by Lemma~\ref{pufundlem1}, there exists a unique 
oriented ($r+1$)-simplex $\tilde{e}_{i}$ of $N_{k}$ such that 
$\varphi _{k,1}(\tilde{e}_{i})=e_{i}''.$ Let $\tilde{d}_{i}$ be 
the oriented ($r+1$)-cell of $|N_{k}|$ which corresponds $\tilde{e}_{i}.$ 
Let $c:= \sum _{i=1}^{t}a_{i}\tilde{d}_{i}\in C_{r+1}(|N_{k}|,\ast ;R).$ Then 
we have $(\gamma _{k})_{\ast }(c)=a.$ 
We shall prove the following claim:  

\noindent Claim: 
$\tilde{a} = \sum _{i=1}^{t}a_{i}\tilde{e}_{i}\in C_{r+1}(N_{k};R)$ is a cycle. 

 In order to prove the claim, 
let $\{ i_{1}w_{1}^{i},i_{2}w_{2}^{i},\ldots ,i_{r+2}w_{r+2}^{i}\} $ be the set of 
vertices of $\tilde{e}_{i}$, where $i_{s}\in \{ 1,\ldots ,m\} $ and 
$w_{s}^{i}\in \{ 1,\ldots ,m\} ^{k-1}$ for each $s=1,\ldots ,r+2.$  
Then, since $\varphi _{k,1}(\tilde{e}_{i})=e_{i}''$, we have that  
the elements $i_{1},i_{2},\ldots ,i_{r+2}$ are mutually distinct. Moreover, we have 
\begin{equation}
\label{pufundlem3eq2}
(\gamma _{k})_{\ast }(\partial (c))=(\gamma _{k})_{\ast }(\partial (\sum _{i=1}^{t}a_{i}\tilde{d}_{i}))=\partial (\sum _{i=1}^{t}a_{i}d_{i})=0.
\end{equation} 
We now suppose that $\partial (\sum _{i=1}^{t}a_{i}\tilde{e}_{i})=
\sum _{j=1}^{\beta } b_{j}e_{j}\neq 0$, 
where for each $j=1,\ldots ,\beta $, $e_{j}$ is an oriented $r$-simplex of $N_{k}$ 
such that  
$\{ e_{1},\ldots ,e_{\beta}\} $ is linearly independent, 
and $b_{j}\in R$ with $b_{j}\neq 0$ for each $j.$ We will deduce a contradiction. 
Let $\{ j_{1}u_{1}^{j},j_{2}u_{2}^{j},\ldots ,j_{r+1}u_{r+1}^{j}\} $ be the set of 
all vertices of $e_{j}$, where $j_{v}\in \{ 1,\ldots ,m\}$ and 
$u_{v}^{j}\in \{ 1,\ldots ,m\} ^{k-1}$ for each $v=1,\ldots ,r+1.$ 
Then, since the elements $i_{1},i_{2},\ldots ,i_{r+2}$ are mutually distinct, we have that 
the elements $j_{1},j_{2},\ldots ,j_{r+1}$ are mutually distinct. 
In particular, denoting by $c_{j}$ the oriented ($r+1$)-cell of $|N_{k}|$ 
which corresponds $e_{j}$, we have that 
$\gamma _{k}(c_{j})$ is an oriented ($r+1$)-cell for each $j.$ 
Moreover, since $\{ e_{1},\ldots ,e_{\beta }\} $ is linearly independent, 
$\{ \gamma _{k}(c_{1}),\ldots ,\gamma _{k}(c_{\beta })\} $ is linearly independent. 
Hence, $(\gamma _{k})_{\ast }(\partial (c))=(\gamma _{k})_{\ast }(\sum _{j=1}^{\beta }b_{j}c_{j})
\neq 0.$ However, it contradicts (\ref{pufundlem3eq2}). 
Therefore, $\partial (\sum _{j=1}^{t}a_{i}\tilde{e}_{i})=0.$ Thus, we have proved the claim. 

 Since $(\gamma _{k})_{\ast }(c)=a$, the above claim implies that  
$(\gamma _{k})_{\ast }: H_{r+1}(|N_{k}|, \ast ;R)\rightarrow H_{r+1}(S_{k},\ast ;R)$ is 
an epimorphism. Thus, we have proved Lemma~\ref{pufundlem3}.         
\end{proof}

Before proving Theorem~\ref{puthm}, 
we state one of the main idea in the proof. 
 One of the keys to proving Theorem~\ref{puthm} is 
the (co)homology sequence of the pair $(|N_{k}|, |\bigcup _{j=1}^{m}N_{k,j}|).$ 
By the cohomology sequence of the above pair, we obtain the following commutative diagram 
of the cohomology groups (with coefficients $R$): 
\begin{equation}
\label{eq:coseqp1}
\begin{CD}
@>>>\tilde{H}^{r}(S_{k})@>>>H^{r}(N_{k})@>{\eta _{k-1}^{\ast }}>>\bigoplus _{j=1}^{m}H^{r}(N_{k-1} )@>>>\tilde{H}^{r+1}(S_{k})@>>>\\  
@. @VV{\tilde{\varphi }_{k+1,k}^{\ast }}V  @VV{\varphi _{k+1,k}^{\ast }}V @VV{\bigoplus _{j=1}^{m}\varphi _{k,k-1}^{\ast }}V 
@VV{\tilde{\varphi }_{k+1,k}^{\ast }}V \\ 
@>>>\tilde{H}^{r}(S_{k+1})@>>>H^{r}(N_{k+1})@>{\eta _{k}^{\ast }}>>\bigoplus _{j=1}^{m}H^{r}(N_{k})@>>>\tilde{H}^{r+1}(S_{k+1})@>>>
\end{CD}
\end{equation}
in which each row is an exact sequence of groups. By (\ref{eq:coseqp1}), we obtain the following 
commutative diagram of the cohomology groups (with coefficients $R$): 
\begin{equation}
\label{eq:coseqp2}
\begin{CD}
@>>>\tilde{H}^{r}(S_{1})@>>>H^{r}(N_{1})@>>>\bigoplus _{j=1}^{m}H^{r}(\{ j\} )@>>>\tilde{H}^{r+1}(S_{1})@>>>\\  
@. @VV{\tilde{\varphi }_{k+1,1}^{\ast }}V  @VV{\varphi _{k+1,1}^{\ast }}V @VV{\bigoplus _{j=1}^{m}\varphi _{k,1}^{\ast }}V 
@VV{\tilde{\varphi }_{k+1,1}^{\ast }}V \\ 
@>>>\tilde{H}^{r}(S_{k+1})@>>>H^{r}(N_{k+1})@>{\eta _{k}^{\ast }}>>\bigoplus _{j=1}^{m}H^{r}(N_{k})@>>>\tilde{H}^{r+1}(S_{k+1})@>>>\\ 
@. @VVV  @VV{\mu _{k+1,r}}V @VV{\bigoplus _{j=1}^{m}\mu _{k,r}^{\ast }}V 
@VVV \\ 
@>>>\varinjlim _{k}\tilde{H}^{r}(S_{k})@>>>\check{H}^{r}({\frak L})@>{\theta }>>\bigoplus _{j=1}^{m}\check{H}^{r}({\frak L})@>>>
\varinjlim _{k}\tilde{H}^{r+1}(S_{k})@>>>\\ 
\end{CD}
\end{equation}
in which each row is an exact sequence of groups, and 
$\varinjlim _{k}\tilde{H}^{\ast }(S_{k})$ denotes the 
direct limit of $\{ \tilde{H}^{\ast }(S_{k}), \tilde{\varphi }_{l,k}^{\ast }\} _{l,k\in \NN , l>k}.$ 
\ 

\ 

\noindent {\bf Proof of Theorem~\ref{puthm}:}
Under the assumptions of Theorem~\ref{puthm}, let $r,k\in \NN .$ 
Using the homology sequence of the pair $(|N_{k}|, |\bigcup _{j=1}^{m}N_{k,j}|)$, 
we have the following exact sequence: 
\begin{multline} \label{puthmpfeq1}
 \tilde{H}_{r+1}(S_{k};T)
\overset{\alpha _{1}}{\longrightarrow } H_{r}(|\bigcup _{j=1}^{m}N_{k,j}|;T) 
\overset{\alpha _{2}}{\longrightarrow} H_{r}(|N_{k}|;T)
\overset{\alpha _{3}}{\longrightarrow }   \tilde{H}_{r}(S_{k};T)\\ 
\overset{\alpha _{4}}{\longrightarrow} H_{r-1}(|\bigcup _{j=1}^{m}N_{k,j}|;T)
\overset{\alpha _{5}}{\longrightarrow}H_{r-1}(|N_{k}|;T)
\overset{\alpha _{6}}{\longrightarrow}\tilde{H}_{r-1}(S_{k};T),  
\end{multline}
where for each $j$, $\alpha _{j}$ denotes some homomorphism. 
Moreover, by Lemma~\ref{pufundlem2}-\ref{pufundlem2-1}, 
we have 
\begin{equation}
\label{puthmpfeq2}
\tilde{H}_{r}(S_{k};T)\overset{(\tilde{\varphi }_{k,1})_{\ast }}{\cong }
\tilde{H}_{r}(S_{1};T) \mbox{ and } 
\tilde{H}_{r-1}(S_{k};T)
\overset{(\tilde{\varphi }_{k,1})_{\ast }}{\cong }\tilde{H}_{r-1}(S_{1};T).
\end{equation} 
Furthermore, by Lemma~\ref{pufundlem3}, we have that 
\begin{equation}
\label{puthmpfeq2-1}
\mbox{Im}(\alpha _{1})=0.
\end{equation}
We now prove statement \ref{puthm1}. 
Let $r,k\geq 2.$ 
By Lemma~\ref{pufundlem3}, we have 
\begin{equation}
\label{puthmpfeq3}
\mbox{Im}(\alpha _{1})=\mbox{Im}(\alpha _{4})=0.
\end{equation}
Moreover, we have the following commutative diagram:
\begin{equation}
\label{nkskcdeq}
\begin{CD}
H_{r}(|\bigcup _{j=1}^{m}N_{k,j}|;T)@>{\alpha _{2}}>>H_{r}(|N_{k}|;T)@>{\alpha_{3}}>>\tilde{H}_{r}(S_{k};T)\\
@VVV @V{(\varphi _{k,1})_{\ast }}VV
@VV{(\tilde{\varphi }_{k,1})_{\ast }}V\\ 
0@>>> H_{r}(|N_{1}|;T)@>>{(\gamma _{1})_{\ast }}>\tilde{H}_{r}(S_{1};T) 
\end{CD}
\end{equation} 
and $(\gamma _{1})_{\ast }:H_{r}(|N_{1}|;T)\rightarrow \tilde{H}_{r}(S_{1};T)$ is an isomorphism, 
where $\gamma _{1}:|N_{1}|\rightarrow S_{1}$ denotes the canonical projection. 
Combining (\ref{puthmpfeq1}), (\ref{puthmpfeq2}), (\ref{puthmpfeq3}), (\ref{nkskcdeq}), and 
Lemma~\ref{nkwlem}, 
we obtain the following exact sequence: 
\begin{equation}
\label{puthmpfeq4}
0\longrightarrow \bigoplus _{j=1}^{m}H_{r}(N_{k-1};T)
\overset{(\eta _{k-1})_{\ast }}{\longrightarrow } H_{r}(N_{k};T)
\overset{(\varphi _{k,1})_{\ast }}{\longrightarrow } H_{r}(N_{1};T)
\longrightarrow 0.
\end{equation}
By (\ref{puthmpfeq4}), we obtain 
$a_{r,k}=ma_{r,k-1}+a_{r,1}.$ Thus, we have proved statement \ref{puthm1}.
Statement \ref{puthm2} follows easily from statement \ref{puthm1}.

We now prove statement \ref{puthmrgeq2}. 
Let $r\geq 2.$ 
By (\ref{puthmpfeq4}), for each $k\in \NN $, 
we have the following exact sequence 
of cohomology groups:
\begin{equation}
\label{puthmpfrgeq2eq1}
0\longrightarrow \check{H}^{r}(N_{1};R)
\overset{(\varphi _{k+1,1})^{\ast }}{\longrightarrow}\check{H}^{r}(N_{k+1};R)
\overset{\eta _{k}^{\ast }}{\longrightarrow }\bigoplus _{j=1}^{m}\check{H}^{r}(N_{k};R)
\longrightarrow 0.
\end{equation}
Taking the direct limit of (\ref{puthmpfrgeq2eq1}) with respect to $k$, we obtain the 
exact sequence (\ref{puthmrgeq2eq}). Thus, we have proved statement \ref{puthmrgeq2}. 

We now prove statement \ref{puthmmukrinj}. 
If $r=0$, then from Lemma~\ref{icfundlem}, $\mu _{k,r}$ and $\varphi _{k}^{\ast }$ are monomorphisms. 
Let 
$r\geq 2$ and $k\geq 1.$ Let ${\frak L}_{k}=(L,(g_{1},\ldots ,g_{m^{k}}))$ 
be a $k$-th iterate of ${\frak L}$. Then, 
there exist isomorphisms $\zeta _{1}: \check{H}^{r}({\frak L};R)_{k}\cong \check{H}^{r}({\frak L}_{k};R)_{1}$ 
and $\zeta _{2}: \check{H}^{r}({\frak L};R)\cong \check{H}^{r}({\frak L}_{k};R)$ such that the following   
diagram commutes:
\begin{equation}
\label{mucommuteeq}
\begin{CD}
\check{H}^{r}({\frak L};R)_{k} @>{\mu _{k,r}}>>\check{H}^{r}({\frak L};R)\\ 
@V{\zeta _{1}}VV @VV{\zeta _{2}}V\\ 
\check{H}^{r}({\frak L}_{k};R)_{1}@>>{\mu _{1,r}}>\check{H}^{r}({\frak L}_{k};R).
\end{CD}
\end{equation}
Moreover, by Lemma~\ref{puiteratelem}, ${\frak L}_{k}$ is postunbranched. 
Combining it with statement~\ref{puthmrgeq2}, we obtain that 
$\mu _{1,r}: \check{H}^{r}({\frak L}_{k};R)_{1}\rightarrow \check{H}^{r}({\frak L}_{k};R)$ 
is a monomorphism. 
Hence, $\mu _{k,r}: \check{H}^{r}({\frak L};R)_{k}\rightarrow 
\check{H}^{r}({\frak L};R)$ is a monomorphism. 
Therefore, statement \ref{puthmmukrinj} follows.

Statement \ref{puthmcohsupp} easily follows from statement~\ref{puthm1} 
and statement~\ref{puthmrgeq2}.  

 
 We now prove statement~\ref{puthmr1fhom}. 
By (\ref{puthmpfeq1}), (\ref{puthmpfeq2}), (\ref{puthmpfeq2-1}), and Lemma~\ref{nkwlem}, 
we obtain the following commutative diagram of homology groups (with coefficients $T$): 
\begin{equation*}
\begin{CD}
0@>>>\bigoplus _{j=1}^{m}H_{1}(N_{k})@>{(\eta _{k})_{\ast }}>>H_{1}(N_{k+1})@>>>\tilde{H}_{1}(S_{1})\\ 
@VVV @VVV @VV{(\varphi _{k+1,1})_{\ast }}V @VV{Id}V\\ 
0@>>>0@>>>H_{1}(N_{1})@>>>\tilde{H}_{1}(S_{1})
\end{CD}
\end{equation*}
\begin{equation}
\begin{CD}
\label{homologycomeq1}
@>>>\bigoplus _{j=1}^{m}H_{0}(N_{k})@>{(\eta _{k})_{\ast }}>>H_{0}(N_{k+1})@>>> 0\\ 
@. @VVV @VV{(\varphi _{k+1,1})_{\ast }}V @VVV \\ 
@>>>\bigoplus _{j=1}^{m}H_{0}(\{ j\} )@>>>H_{0}(N_{1})@>>>0
\end{CD}
\end{equation}
in which each row is an exact sequence of groups. By (\ref{homologycomeq1}), 
it is easy to see that statement~\ref{puthmr1fhom} holds. 

 We now prove statement \ref{puthmr1f} and statement \ref{puthmr1}. 
By the cohomology sequence of the pair 
$(|N_{k}|, |\bigcup _{j=1}^{m}N_{k,j}|)$, (\ref{eq:coseqp1}), (\ref{eq:coseqp2}), 
(\ref{puthmpfeq2}), (\ref{puthmpfeq2-1}), and Lemma~\ref{nkwlem},  
for each $k\in \NN $, 
we have the following commutative diagram of cohomology groups (with coefficients $R$): 
\begin{equation*}
\begin{CD}
0@>>>H^{0}(N_{1})@>>>\bigoplus _{j=1}^{m}H^{0}(\{ j\} )@>>>\tilde{H}^{1}(S_{1})\\  
@VVV @VV{\varphi _{k+1,1}^{\ast }}V @VVV @VV{Id}V \\ 
0@>>>H^{0}(N_{k+1})@>{\eta _{k}^{\ast }}>>\bigoplus _{j=1}^{m}H^{0}(N_{k})@>>>\tilde{H}^{1}(S_{1})\\ 
@VVV @VV{\mu _{k+1,0}}V @VV{\bigoplus _{j=1}^{m}\mu _{k,0}}V @VV{Id}V \\ 
0@>>>\check{H}^{0}({\frak L})@>{\theta }>>\bigoplus _{j=1}^{m}\check{H}^{0}({\frak L})@>>>\tilde{H}^{1}(S_{1})
\end{CD}
\end{equation*}
\begin{equation}
\label{3bigeq}
\begin{CD}
@>>>H^{1}(N_{1})@>>>0@>>>0\\
@. @VV{\varphi _{k+1,1}^{\ast }}V @VVV @VVV \\
@>>>H^{1}(N_{k+1})@>{\eta _{k}^{\ast }}>>\bigoplus _{j=1}^{m}H^{1}(N_{k})
@>>> 0\\  
@. @VV{\mu _{k+1,1}}V @VV{\bigoplus _{j=1}^{m}\mu _{k,1}}V @VVV \\
 @>>>\check{H}^{1}({\frak L})
@>{\theta }>>\bigoplus _{j=1}^{m}\check{H}^{1}({\frak L})@>>> 0
\end{CD}
\end{equation}
in which each row is an exact sequence of groups. 
By (\ref{3bigeq}), it is easy to see that 
statement \ref{puthmr1f} and statement \ref{puthmr1} hold. 
Thus, we have proved statement~\ref{puthmr1f} and statement~\ref{puthmr1}. 


We now prove statement~\ref{puthmakvalues1}. 
By (\ref{3bigeq}), we have the following exact sequence:
\begin{equation}
0\rightarrow H^{0}(N_{1};R)\rightarrow H^{0}(\{ 1,\ldots ,m\} ;R)\rightarrow \tilde{H}^{1}(S_{1};R)
\rightarrow H^{1}(N_{1};R)\rightarrow 0.
\end{equation}
Hence we have $\dim _{R}H^{1}(S_{1};R)=m-a_{0,1}+a_{1,1}.$ Combining it with 
the exact sequences (\ref{puthmr1feq2}) and (\ref{puthmr1feq1}), 
we can easily obtain that statement~\ref{puthmakvalues1} holds. 
Thus we have proved statement~\ref{puthmakvalues1}.  

 Statement~\ref{puthmlambda} easily follows from the definition of $\lambda _{k}$ and 
$b_{1,\infty }.$ 

  Statement~\ref{puthm3} easily follows from statement~\ref{puthmakvalues1}.   

 Statement~\ref{puthmakvalues2} and statement~\ref{puthm4} easily follow from 
statement~\ref{puthmakvalues1} and statement~\ref{puthmlambda}. 
  
 Statement~\ref{puthmulrgeq1} easily follows from statement~\ref{puthm1} and statement~\ref{puthmakvalues2}.

 Statement~\ref{puthmulr0} easily follows from statement~\ref{puthm4}.  
  
 We now prove statement \ref{puthmavalues}. 
By statements \ref{puthmrgeq2} and \ref{puthmr1}, 
for each $r\in \NN $ with $r\geq 2$ there exists an 
exact sequence $\check{H}^{r}(\frak{L};R)\rightarrow \bigoplus _{j=1}^{m}\check{H}^{r}(\frak{L};R)\rightarrow 0.$ 
Hence, if $m>1$, then either $a_{r,\infty }=0$ or $a_{r,\infty }=\infty .$ 
If $m=1$, then obviously we have $a_{r,\infty }=0.$ Thus, we have proved statement \ref{puthmavalues}.   

We now prove statement~\ref{puthm8}. 
Suppose $a_{0,\infty }<\infty .$ 
Since $a_{0,\infty }=\lim _{k\rightarrow \infty }a_{0,k}$, 
statement~\ref{puthmakvalues1} and statement~\ref{puthmlambda} imply that 
$a_{0,\infty }=ma_{0,\infty }-m+a_{0,1}-a_{1,1}+b_{1,\infty }.$ 
Therefore, statement~\ref{puthm8} follows.  
 
 Statement~\ref{puthm9} easily follows from statement~\ref{puthmavalues} and statement~\ref{puthm8}. 

We now prove statement \ref{puthm5}. 
Suppose that there exists an element $k_{0}\in \NN $ such that 
$a_{0,k_{0}}>\frac{1}{m-1}
(m-a_{0,1}+a_{1,1})$. 
We will show that $a_{0,k+1}>a_{0,k}$ for each $k\geq k_{0}$, by using induction on $k\geq k_{0}.$ 
For the first step, by statement \ref{puthm4} and the assumption 
$a_{0,k_{0}}>\frac{1}{m-1}
(m-a_{0,1}+a_{1,1})$, we have 
$a_{0,k_{0}+1}-a_{0,k_{0}}\geq (m-1)a_{0,k_{0}}-m+a_{0,1}-a_{1,1}>0.$ 
We now suppose that $a_{0,k+1}>a_{0,k}$ for each $k\in \{ k_{0},k_{0}+1,k_{0}+2,\ldots ,t\} .$ 
Then, by statement \ref{puthm4}, we have 
$a_{0,t+2}-a_{0,t+1}\geq (m-1)a_{0,t+1}-m+a_{0,1}-a_{1,1}
\geq (m-1)a_{0,k_{0}}-m+a_{0,1}-a_{1,1}>0.$ Therefore, inductive step is completed. 
Thus, we have proved statement \ref{puthm5}. 

 Statement \ref{puthm6} follows easily from statement \ref{puthm5} (or from statement~\ref{puthm8} and 
statement~\ref{puthmlambda}).

We now prove statement~\ref{puthmless6}. 
Suppose $2\leq m\leq 6$ and $|N_{1}|$ is disconnected. 
In order to show $a_{0,\infty }=\infty $, by Lemmas~\ref{1disconlem} and \ref{ss1disconlem} we may assume that 
each connected component of $|N_{1}|$ has at least two vertices. 
Then it is easy to see that $\frac{m-a_{0,1}+a_{1,1}}{m-1}<2.$ 
Therefore statement~\ref{puthm6} implies that $a_{0,\infty }=\infty .$ 
By Lemma~\ref{icfundlem}-\ref{icfundlem2}, it follows that 
$L$ has infinitely many connected components. Thus we have proved statement~\ref{puthmless6}.

 We now prove statement~\ref{puthmB2zero}. 
Suppose that $B_{2}=0.$ 
Let $C_{k}:=\mbox{Im} 
(\varphi _{k}^{\ast }: \check{H}^{1}({\frak L};R)_{k}\rightarrow \check{H}^{1}({\frak L};R)_{k+1}).$ 
We will show the following claim: \\ 
Claim: For each $k\in \NN $, 
$C_{k}=0.$  

To prove the claim, we will use the induction on $k.$ 
Since $B_{2}=0,$ we have $C_{1}=0.$ 
Moreover, since $B_{2}=0$, the exact sequence (\ref{puthmr1feq2}) implies that 
for each $k\in \NN $, 
$\eta _{k}^{\ast }: \check{H}^{1}({\frak L};R)_{k+1}\rightarrow \bigoplus _{j=1}^{m}\check{H}^{1}({\frak L};R)_{k}$ 
is an isomorphism. 
Furthermore, for each $k\in \NN $, 
we have the following commutative diagram.  
\begin{equation}
\begin{CD}
\check{H}^{1}({\frak L};R)_{k+1}@>{\eta _{k}^{\ast }}>>\bigoplus _{j=1}^{m}\check{H}^{1}({\frak L};R)_{k}\\ 
@V{\varphi _{k+1}^{\ast }}VV @VV{\bigoplus _{j=1}^{m}\varphi _{k}^{\ast }}V\\ 
\check{H}^{1}({\frak L};R)_{k+2}@>>{\eta _{k+1}^{\ast }}> \bigoplus _{j=1}^{m}\check{H}^{1}({\frak L};R)_{k+1}.
\end{CD}
\end{equation}
Hence, if we assume $C_{k}=0$, then $C_{k+1}=0.$ 
Therefore, the induction is completed. 
Thus, we have proved the claim. 

 From the above claim, it is easy to see that 
$\check{H}^{1}({\frak L};R)=0.$ 
Hence, we have proved statement~\ref{puthmB2zero}.   

We now prove statement \ref{puthm7}.
For each $j=1,\ldots, m$, let $c_{j}:N_{k}\rightarrow \{ j\} $ be the 
constant map. 
By (\ref{homologycomeq1}), 
we have the following commutative diagram of homology groups (with coefficients $T$): 
\begin{equation}
\label{puthmpfeq10}
\begin{CD}
0@>>>\bigoplus _{j=1}^{m}H_{1}(N_{k})@>{(\eta _{k})_\ast }>>H_{1}(N_{k+1})
@>>>\check{H}_{1}(S_{1})@>>>\bigoplus _{j=1}^{m}H_{0}(N_{k})\\ 
@VVV @VVV @VV{(\varphi _{k+1,1})_{\ast }}V @VV{Id}V 
@VV{\bigoplus _{j=1}^{m}(c_{j})_{\ast }}V\\ 
0@>>>0@>>>H_{1}(N_{1})@>>>\check{H}_{1}(S_{1})@>>>
\bigoplus _{j=1}^{m}H_{0}(\{ j\})
\end{CD}
\end{equation}
in which each row is an exact sequence of groups.

Suppose that $|N_{1}|$ is connected. Then, by Lemma~\ref{1conall}, $|N_{k}|$ is connected for each $k\in \NN $. 
Hence, $\bigoplus _{j=1}^{m}(c_{j})_{\ast }: \bigoplus _{j=1}^{m}{H}_{0}(N_{k};T)\rightarrow \bigoplus _{j=1}^{m}H_{0}(\{ j\} ;T)$ 
is an isomorphism. Combining it with (\ref{puthmpfeq10}), 
the five lemma implies that 
$(\varphi _{k+1,1})_{\ast }: H_{1}(N_{k+1};T)\rightarrow H_{1}(N_{1};T)$ is an epimorphism. 
Combining it with 
(\ref{puthmr1fhomeq1}),  
we obtain the exact sequence (\ref{puthm7aeq}) in statement \ref{puthm7a}. 
Hence, we have proved statement \ref{puthm7a}.  
Statement \ref{puthm7b} easily follows from statement \ref{puthm7a}. 
We now prove statement \ref{puthm7c}. 
If $a_{1,1}=0$, then statement \ref{puthm7b} implies that for each $k\in \NN $, 
$a_{1,k}=0.$ Therefore, $a_{1,\infty }=0.$ 
If $a_{1,1}\neq 0$, then statement \ref{puthm7b} implies that 
$a_{1,k}\rightarrow \infty $ as $k\rightarrow \infty .$ 
From Lemma~\ref{icfundlem}-\ref{icfundlem3}, it follows that 
$a_{1,\infty }=\infty .$ Therefore, we have proved statement \ref{puthm7c}. 
Statement \ref{puthm7d} follows from statement \ref{puthm7a} and the universal-coefficient theorem. 
We now prove statement \ref{puthm7e}. 
By the exact sequence (\ref{puthm7aeq}), for each $k\in \NN $ we have the following exact sequence: 
\begin{equation}
\label{puthm7epfeq1}
0\longrightarrow H^{1}(N_{1};R)
\overset{(\varphi _{k+1,1})^{\ast }}{\longrightarrow}H^{1}(N_{k+1};R)
\overset{\eta _{k}^{\ast }}{\longrightarrow }\bigoplus _{j=1}^{m}H^{1}(N_{k};R)
\longrightarrow 0.
\end{equation}
Taking the direct limit of (\ref{puthm7epfeq1}) with respect to $k$, we obtain the 
exact sequence (\ref{puthm7eeq}). 

Therefore, we have proved statement \ref{puthm7}. 

Thus, we have proved Theorem~\ref{puthm}.  
\qed 

\ 

We now prove Proposition~\ref{rgeq2exprop}.\\ 
\noindent {\bf Proof of Proposition~\ref{rgeq2exprop}:} 
We first prove statement \ref{rgeq2exprop1}.  
Let $K:=\{ z\in \CC \mid 1\leq |z|\leq 2\} .$ 
It is easy to see that there exists a finite family 
$\{ h_{1},\ldots, h_{n+2}\} $ 
of topological branched covering maps 
on $\CCI $ with the following properties:  
\begin{enumerate}
\item \label{hjp1}
$h_{j}^{-1}(K)\subset K$ for each $j=1,\ldots ,n+2$;  
\item \label{hjp2}
$h_{i}^{-1}(\mbox{int}(K))\cap h_{j}^{-1}(\mbox{int}(K))=\emptyset $ 
for each $(i,j)$ with $i\neq j$; 
\item \label{hjp3}
$h_{1}^{-1}(K)\cap \{ z\in \CC \mid |z|=2\} =\emptyset $ and 
$h_{2}^{-1}(K)\cap \{ z\in \CC \mid |z|=1\} =\emptyset $;   
\item \label{hjp4}
$h_{j}^{-1}(K)\subset \mbox{int}(K)$ for each $j=3,\ldots ,n+2$; 
\item \label{hjp5}
$h_{j}|_{\{ z\in \CC \mid |z|=j\} } =\mbox{Id}$ for each $j=1,2$;  
\item \label{hjp6}
$h_{1}(h_{1}^{-1}(K)\cap h_{k}^{-1}(K))\subset \{ z\in \CC \mid |z|=2\} $ 
for each $k\in \{ 1,\ldots ,n+2\}$ with $k\neq 1$;
\item \label{hjp7}
$h_{2}(h_{2}^{-1}(K)\cap h_{k}^{-1}(K))\subset \{ z\in \CC \mid |z|=1\} $ 
for each $k\in \{ 1,\ldots ,n+2\}$ with $k\neq 2$; 
\item \label{hjp8}
$h_{j}(h_{j}^{-1}(K)\cap h_{k}^{-1}(K))\subset \{ z\in \CC \mid |z|=2\} $ 
for each $j,k\in \{ 1,\ldots ,n+2\}$ with $j\neq k$ and $j\geq 3$; and 
\item \label{hjp9}
$\bigcap _{i=1}^{n+2}h_{i}^{-1}(\partial K)=\emptyset $ and 
for each $j=1,\ldots ,n+2$, 
$\bigcap _{i\in \{ 1,\ldots ,n+2\} \setminus \{ j\}}h_{i}^{-1}(\partial K)
\neq \emptyset .$ 
\end{enumerate}
Let $L:=R_{K,b}(h_{1},\ldots ,h_{n+2})$ and  
let ${\frak L}:=(L, (h_{1},\ldots ,h_{n+2})).$ Then, by Lemma~\ref{Rbsslem}, 
${\frak L}$ is a backward self-similar system. 
From properties \ref{hjp1}, \ref{hjp3}, \ref{hjp4}, and \ref{hjp5},  
we have that for each $j=1,2$,  
$\{ z\in \CC \mid |z|=j\} \subset L_{(j)^{\infty }}\setminus L_{x}$ 
for any $x\in \Sigma _{n+2}$ with $x\neq (j)^{\infty }.$  
Combining it with properties \ref{hjp6}, \ref{hjp7}, and \ref{hjp8}, 
it follows that ${\frak L}$ is postunbranched. 
Moreover, since $\partial K=\bigcup _{j=1}^{2}\{ z\in \CC \mid 
|z|=j\} \subset L$, 
properties \ref{hjp2} and \ref{hjp9} imply that 
$\bigcap _{i=1}^{n+2}h_{i}^{-1}(L)=\emptyset $ and 
for each $j=1,\ldots ,n+2$, 
$\bigcap _{i\in \{ 1,\ldots ,n+2\} \setminus \{ j\} }
h_{i}^{-1}(L)\neq \emptyset .$ 
Hence $H^{n}(N_{1};R)=R$ for each field $R.$ 
From Theorem~\ref{puthm}-\ref{puthmcohsupp}, 
it follows that $\dim _{R}\check{H}^{n}({\frak L};R)=\infty $ 
for each field $R.$  Thus, we have proved 
statement \ref{rgeq2exprop1} of Proposition~\ref{rgeq2exprop}. 

 We now prove statement \ref{rgeq2exprop2} of Proposition~\ref{rgeq2exprop}. 
Let $K':= \{ z\in \CC \mid 1\leq |z|\leq 2\} \times 
[0,1]\subset \RR ^{3}.$ 
We can construct a finite family $\{ h_{j}\} _{j=1}^{n+2}$ of 
continuous and injective maps on $K'$ satisfying properties similar to 
the above properties \ref{hjp1},...,\ref{hjp9} (with ``$-1$'' removed). 
Let $L:= R_{K',f}(h_{1},\ldots ,h_{n+2})$ and let 
${\frak L}:= (L,(h_{1},\ldots ,h_{n+2})).$ Then, 
by the argument similar to that in the previous paragraph, 
we obtain that ${\frak L} $ is postunbranched, 
$H^{n}(N_{1};R)=R$ for each field $R$, and 
$\dim _{R}\check{H}^{n}({\frak L};R)=\infty $ for each field $R.$ 
Thus, we have proved statement \ref{rgeq2exprop2} of 
Proposition~\ref{rgeq2exprop}. 
\qed 

\

\begin{figure}[htbp]
\caption{The invariant set of a postunbranched system such that $\mu _{1,1}$ is not injective. }    
\ \ \ \ \ \ \ \ \ \ \ \ \ \ \ \ \ \ \ \ \ \ \ \ \ \ \ \ \ \ \ \ 
\ \ \ \  \ \ \ \ \ \ \ \ \ \ \ \ \ \ 
\includegraphics[width=2.7cm,width=2.7cm]{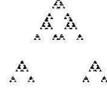}
\label{fig:mu11notinj1}
\end{figure}

 We now prove Proposition~\ref{notinjmuexprop}. \\ 
\noindent {\bf Proof of Proposition~\ref{notinjmuexprop}:} 
Let $p_{1},p_{2},p_{3}\in \CC $ be mutually distinct three points such that 
$p_{1}p_{2}p_{3}$ makes an 
equilateral triangle. 
For each $j=1,2,3,$  
let $g_{j}(z)=\frac{1}{2}(z-p_{j})+p_{j}.$ 
Let $h_{1}:=g_{1}^{2}, h_{2}:=g_{2}^{2}, h_{3}:= g_{3}^{2}, 
h_{4}:=g_{3}\circ g_{1},$ and $h_{5}:=g_{3}\circ g_{2}.$  
Let $L:= M_{\CC }(h_{1},\ldots ,h_{5})$ and 
let ${\frak L} := (L,(h_{1},\ldots ,h_{5}))$ (see Figure~\ref{fig:mu11notinj1}). Then 
${\frak L}$ is a forward self-similar system. 
By Example~\ref{SGex},
${\frak L} $ is postunbranched. 
Since $p_{j}\in L$ for each $j=1,2,3$, 
we have that $h_{3}(L)\cap h_{4}(L)\neq \emptyset $, 
$h_{4}(L)\cap h_{5}(L)\neq \emptyset $, 
$h_{5}(L)\cap h_{3}(L)\neq \emptyset $, 
and $h_{3}(L)\cap h_{4}(L)\cap h_{5}(L)=\emptyset .$ 
Moreover, it is easy to see that for each $r\in \NN $ with $r\geq 2$, 
there exists no $r$-simplex of $N_{1}=N_{1}({\frak L}).$  
Hence $\tau =(3,4)(4,5)(5,3)$ is a closed edge path of $N_{1}=N_{1}({\frak L})$ 
which induces a non-trivial element of $H_{1}(N_{1};R)$ for each field $R.$ 
Hence $\check{H}^{1}({\frak L};R)_{1}\neq 0.$ 
However, considering $N_{2}$, 
it is easy to see that 
$\mbox{Im}((\varphi _{1})_{\ast }: H_{1}(N_{2};R)\rightarrow H_{1}(N_{1};R))=0.$ 
Hence, $B_{2}=0$. From Theorem~\ref{puthm}-\ref{puthmB2zero}, 
it  follows that $\check{H}^{1}({\frak L};R)=0$. 
Moreover, since $\check{H}^{1}({\frak L};R)_{1}\neq 0$, we obtain that 
$\mu _{1,1}:\check{H}^{1}({\frak L};R)_{1}
\rightarrow \check{H}^{1}({\frak L};R)$ is not injective. 
Furthermore, since each $h_{j}:L\rightarrow L$ is a contraction, 
we have that $\Psi: \check{H}^{1}({\frak L}:R)\rightarrow \check{H}^{1}(L;R)$ is an isomorphism. 
Hence $\check{H}^{1}(L;R)=0.$ From the Alexander duality theorem (\cite{Sp}), 
it follows that $\CC \setminus L$ is connected.  
Thus, we have proved Proposition~\ref{notinjmuexprop}. 
\qed  

In order to prove Theorem~\ref{sscij1thm}, we need several lemmas. 
\begin{lem}
\label{sscij1lem1}
Let $\frak{L}=(L,(h_{1},\ldots ,h_{m}))$ be a forward self-similar system 
such that 
for each $j=1,\ldots ,m$, $h_{j} :L\rightarrow L$ is 
injective.  
Suppose that $\sharp C_{i,j}\leq 1$ for each $(i,j)$ with $i\neq j.$ Then, 
for each $r\geq 2$, each $k\geq 1$ and each $\ZZ $ module $T$, we have 
$\check{H}_{r}(\frak{L};T)_{k}=\check{H}_{r}(\frak{L};T)=0.$
\end{lem}
\begin{proof}
 Let $a=\sum _{i=1}^{t}a_{i}d_{i}\in C_{r}(N_{k};T)$ be a cycle, where 
 for each $i$, $a_{i}\in T$ and $d_{i}$ is an oriented $r$-simplex. We may assume 
that $\{ d_{1},\ldots ,d_{t}\} $ is linearly independent.  
Let $\Omega $ be the graph such that the vertex set is equal to 
$\{d_{1},\ldots ,d_{t} \} $ and such that $\{ d_{i},d_{j}\} $ is an edge 
if and only if there exists a $1$-simplex $e$ of $N_{k}$ with $|e|\subset |d_{i}|\cap |d_{j}|.$ 
Let $\{ \Omega _{1},\Omega _{2},\ldots ,\Omega _{p}\} $ be the set of 
all connected components of $|\Omega |.$ 
Then we have 
$\sum _{i=1}^{t}a_{i}d_{i}=\sum _{l=1}^{p}\sum _{d_{i}\in \Omega _{l}}a_{i}d_{i}.$ 
We now show the following claim: \\ 
\noindent Claim 1: 
For each $l$, $\partial (\sum _{d_{i}\in \Omega _{l}}a_{i}d_{i})=0$ in $C_{r-1}(N_{k};T).$ 

 In order to show claim 1, suppose that there exists an $l$ such that 
 $\partial (\sum _{d_{i}\in \Omega _{l}}a_{i}d_{i})
=\sum _{j=1}^{\beta }b_{j}e_{j}\neq 0$, 
where $e_{j}$ is an oriented $r-1$ simplex of $N_{k}$ for each $j$, 
$\{ e_{1},\ldots ,e_{\beta }\} $ is linearly independent, 
and $b_{j}\in T$ with $b_{j}\neq 0$ for each $j.$  
 Since $\partial (\sum _{i=1}^{t}a_{i}d_{i})=0$, 
there exists an $l'$ with $l'\neq l$ and an element $d_{q}\in \Omega _{l'}$ 
such that $|d_{q}|\supset |e_{1}|.$ However, it implies that 
$d_{q}\in \Omega _{l}$ and this is a contradiction since 
$\Omega _{l}\cap \Omega _{l'}=\emptyset .$ Hence, we have proved claim 1.

  We now prove the following claim: \\ 
\noindent Claim 2: 
Let $l\in \{ 1,\ldots ,p\} $ be a number. 
Let $\{ v_{0},\ldots ,v_{s}\} $ be the union 
$\bigcup _{d_{i}\in \Omega _{l}} \{ \mbox{all vertices of } d_{i}\} .$ 
Then, $M_{l}:= \{ v_{0},\ldots ,v_{s}\} $ is an $s$-simplex of $N_{k}.$ 

 In order to prove claim 2, let  
$d_{i},d_{j}\in \Omega _{l}$ be two elements 
such that there exists a $1$-simplex 
$e=\{ u_{1},u_{2}\} $ of $N_{k}$ with $|d_{i}|\cap |d_{j}|\supset |e|$, 
where $u_{1},u_{2}\in \{ 1,\ldots ,m\} ^{k}$. 
Let $\{ w_{0},\ldots ,w_{r}\} $ be the set of 
all vertices of $d_{i}$ and let 
$\{ w_{0}',\ldots ,w_{r}'\} $ be the set of all 
vertices of $d_{j}.$ 
Then we have 
\begin{equation}
\label{sscijlem1pfeq1}
\emptyset \neq \bigcap _{j=0}^{r}h_{\overline{w_{j}}}(L)\subset h_{\overline{u_{1}}}(L)\cap h_{\overline{u_{2}}}(L)\ 
\mbox{ and } 
\emptyset \neq \bigcap _{j=0}^{r}h_{\overline{w_{j}'}}(L)\subset h_{\overline{u_{1}}}(L)\cap h_{\overline{u_{2}}}(L). 
\end{equation}
Since $\sharp C_{i,j}\leq 1 $ for each $(i,j)$ with $i\neq j$ and 
$h_{j}: L\rightarrow L $ is injective for each $j$,  
we have $\sharp (h_{\overline{u_{1}}}(L)\cap h_{\overline{u_{2}}}(L))\leq 1.$ 
Combining it with (\ref{sscijlem1pfeq1}), it follows that 
there exists a point $z\in L$ such that 
$\bigcap _{j=0}^{r}h_{\overline{w_{j}}}(L)=\bigcap _{j=0}^{r}h_{\overline{w_{j}'}}(L)=\{ z\} .$
The above argument implies that $\bigcap _{j=0}^{s}h_{\overline{v_{j}}}(L)=\{ z\} .$ 
Hence, $M_{l}=\{ v_{0},\ldots ,v_{s}\} $ is an $s$-simplex of $N_{k}.$      
Therefore, we have proved claim 2.

 By Claim 2, we obtain that for each $l$, 
 $\sum _{d_{i}\in \Omega _{l}}a_{i}d_{i}\in C_{r}(M_{l};T).$ 
Combining it with claim 1, we get that for each $l$, 
$\sum _{d_{i}\in \Omega _{l}}a_{i}d_{i}$ is a cycle of 
$C_{r}(M_{l};T).$ Since $H_{r}(M_{l};T)=0$, 
it follows that for each $l$, $\sum _{d_{i}\in \Omega _{l}}a_{i}d_{i}$ 
is a boundary element of $C_{r}(N_{k};T).$ Hence, 
we get that $H_{r}(N_{k};T)=0.$ Therefore, we have proved 
Lemma~\ref{sscij1lem1}. 
\end{proof}
By the same method, we can prove the following lemma.
\begin{lem}
Let ${\frak L}=(L,(h_{1},\ldots ,h_{m}))$ be a backward self-similar system. 
Suppose that $\sharp C_{i,j}\leq 1$ for each $(i,j)$ with $i\neq j.$ 
Let $T$ be a $\ZZ $ module. Then, for each $r\in \NN $ with $r\geq 2$, 
we have $\check{H}_{r}({\frak L};T)_{1}=0.$ 
\end{lem}

\begin{lem}\label{sscij1lem2}
Let $\frak{L}=(L,(h_{1},\ldots ,h_{m}))$ be a forward self-similar 
system such that for each $j=1,\ldots ,m$, 
$h_{j}: L\rightarrow L$ is injective. 
 Let $T$ be a $\ZZ $ module. Suppose that 
$\sharp C_{i,j}\leq 1$ for each $(i,j)$ with $i\neq j.$ 
Then, for each $k\in \NN $, 
$H_{2}(S_{k};T)=0.$ 
\end{lem}
\begin{proof}
Let $a=\sum _{i=1}^{r}a_{i}d_{i}\in C_{2}(S_{k};T)$ be a cycle, 
where for each $i$, $a_{i}\in T$ and $d_{i}$ is an oriented $2$-cell of $S_{k}.$ 
We will show that $a$ is a boundary. 
Let $\gamma _{k}: |N_{k}|\rightarrow S_{k}$ be the canonical projection. 
For each $i$, let $\tilde{d}_{i}$ be an oriented $2$-simplex of $N_{k}$ such that 
$\gamma _{k}(|\tilde{d}_{i}|)=d_{i}.$  
Let $\Omega $ be the graph such that the vertex set is equal to $\{ d_{1},\ldots ,d_{r}\} $ 
and such that $\{ d_{i},d_{j}\} $ is an edge of $\Omega $ if and only if 
there exists an $1$-cell $e$ of $S_{k}$ such that 
$d_{i}\cap d_{j}\supset e.$ 
Let $\{ \Omega _{1},\ldots ,\Omega _{p}\} $ be the set of 
all connected components of $|\Omega |.$ 
Then we have $a=\sum _{l=1}^{p}\sum _{d_{i}\in \Omega _{l}}a_{i}d_{i}.$ 
We now prove the following claim.\\ 
Claim 1: For each $l$, $\partial (\sum _{d_{i}\in \Omega _{l}}a_{i}d_{i})=0$ in 
$C_{1}(S_{k};T).$ 

 In order to prove claim 1, suppose that the statement is false. 
Then, there exists an $l$ such that 
$\partial (\sum _{d_{i}\in \Omega _{l}}a_{i}d_{i})=\sum _{j=1}^{\beta }b_{j}e_{j}\neq 0$, 
where for each $j$, $b_{j}\in T$ and $e_{j}$ is an oriented $1$-cell of $S_{k}$ such that 
$\{ e_{1},\ldots ,e_{\beta }\} $ is linearly independent. 
Since $\partial (\sum _{i=1}^{r}a_{i}d_{i})=0$, it follows that 
there exists an $l'$ with $l'\neq l$ and an element $d_{i}\in \Omega _{l'}$ such that 
$d_{i}\supset e_{1}.$ It implies that $d_{i}\in \Omega _{l}.$ 
However, this is a contradiction, since $l'\neq l.$ Therefore, we have proved claim 1. 

 We now prove the following claim.\\ 
\noindent Claim 2: For each $l$, there exists an $s\in \NN $ with $s\geq 2$ and an $s$-simplex $M$ of $N_{k}$ 
such that $\bigcup _{d_{i}\in \Omega _{l}}d_{i}\subset \gamma _{k} (|M|).$ 

 In order to prove claim 2, let $d_{i}\in \Omega _{l}$ be an element. 
Let $d_{j}\in \Omega _{l}$ be another element such that 
$\{ d_{i},d_{j}\} $ is an edge of $\Omega _{l}.$ 
Then there exist four vertices $w_{1},w_{2},w_{3},w_{4}$ of $N_{k}$ such that 
the set of vertices of $\tilde{d}_{i}$ is equal to $\{ w_{1},w_{2},w_{3}\} $ and 
the set of vertices of 
$\tilde{d}_{j}$ is equal to $ \{ w_{2},w_{3},w_{4}\} .$ 
We have  
$\bigcap _{j=1}^{3}h_{\overline{w_{j}}}(L)\neq \emptyset $ and 
$\bigcap _{j=2}^{4}h_{\overline{w_{j}}}(L)\neq \emptyset .$ 
Since $\sharp C_{i,j} \leq 1$ for each $(i,j)$ with $i\neq j$ and 
$h_{j}: L\rightarrow L$ is injective for each $j$, 
there exists a point $z\in L$ such that 
$h_{\overline{w_{2}}}(L)\cap h_{\overline{w_{3}}}(L)=\{ z\} .$ 
Therefore, $\bigcap _{j=1}^{4}h_{\overline{w_{j}}}(L)=\{ z\} .$ 
This argument implies that denoting by $\{ v_{1},\ldots ,v_{s}\} $ the set of 
all vertices of $\bigcup _{d_{i}\in \Omega _{l}}\tilde{d}_{i}$, 
we have $\bigcap _{j=1}^{s}h_{v_{j}}(L)=\{ z\} .$ 
Let $M=\{ v_{1},\ldots ,v_{s}\} .$ Then $M$ is an $s$-simplex of $N_{k}$ and 
$\bigcup _{d_{i}\in \Omega _{l}}d_{i}\subset \gamma _{k}(|M|).$ Thus, we have proved claim 2. 

 Since $\gamma _{k} (|M|)$ is a subcomplex of $S_{k}$ and 
$\sum _{d_{i}\in \Omega _{l}}a_{i}d_{i}$ is a cycle of $C_{2}(S_{k};T)$, 
we obtain that $\sum _{d_{i}\in \Omega _{l}}a_{i}d_{i}$ is a cycle of 
$C_{2}(\gamma _{k} (|M|);T).$  We now prove the following claim.\\ 
\noindent Claim 3: 
$H_{2}(\gamma _{k} (|M|);T)=0.$ 

 In order to prove claim 3, let $\tilde{\gamma }_{k}:|M|/(|M\cap \bigcup _{j=1}^{m}N_{k,j}|)\rightarrow \gamma _{k}(|M|)$ be 
 the cellular map induced by $\gamma _{k}.$ Then, $\tilde{\gamma }_{k}$ is a homeomorphism. 
Moreover, we have the following homology sequence of the pair  
$(|M|, |\bigcup _{j=1}^{m}N_{k,j}|)$: 
\begin{equation}
\cdots \rightarrow H_{2}(|M|;T)\rightarrow H_{2}(|M|/|M\cap \bigcup _{j=1}^{m}N_{k,j}|;T)\rightarrow 
H_{1}(|M\cap \bigcup _{j=1}^{m}N_{k,j}|;T)\rightarrow \cdots .
\end{equation} 
Since $M$ is an $s$-simplex, $H_{2}(|M|;T)=0.$ 
Moreover, 
$$H_{1}(|M\cap \bigcup _{j=1}^{m}N_{k,j}|;T)\cong \bigoplus _{j=1}^{m}H_{1}(|M\cap N_{k,j}|;T).$$  
Let $\{ u_{1},\ldots ,u_{t}\} $ be the set of all vertices of $M\cap N_{k,j}.$ 
Then, $u=\{ u_{1},\ldots ,u_{t}\} $ is a $(t-1)$-simplex of $M$.  
Since $M$ is a subcomplex of $N_{k}$, we obtain that $u$ is a simplex of $N_{k}.$ 
Moreover, since each $u_{j}$ is a vertex of $N_{k,j}$, it follows that 
$u$ is a simplex of $N_{k,j}.$ Therefore, 
$u$ is a simplex of $M\cap N_{k,j}.$ 
Hence, 
$H_{1}(|M\cap N_{k,j}|;T)=0.$ Combining these arguments, 
we obtain that $H_{2}(\gamma _{k} (|M|);T)=0.$ Thus, we have proved claim 3. 

 By claim 3, the cycle $\sum _{d_{i}\in \Omega _{l}}a_{i}d_{i}\in 
C_{2}(\gamma _{k}(|M|);T)$ is a boundary element of $C_{2}(\gamma _{k}(|M|);T).$ 
Therefore, $\sum _{d_{i}\in \Omega _{l}}a_{i}d_{i}$ is a boundary element of 
$C_{2}(S_{k};T).$ Hence, $a=\sum _{j=1}^{r}a_{i}d_{i}$ is a boundary element of 
$C_{2}(S_{k};T).$ Thus, we have proved  Lemma~\ref{sscij1lem2}.     
\end{proof}
By the same method, we can prove the following lemma.
\begin{lem}
Let ${\frak L}=(L,(h_{1},\ldots ,h_{m}))$ be a backward self-similar system. 
Suppose that $\sharp C_{i,j}\leq 1$ for each $(i,j)$ with $i\neq j.$ 
Let $T$ be any $\ZZ $ module. Then, 
$H_{2}(S_{1};T)=0.$ 
\end{lem}

 We now prove Theorem~\ref{sscij1thm}.\\ 
{\bf Proof of Theorem~\ref{sscij1thm}:} 
From Lemma~\ref{sscij1lem1}, statement \ref{sscij1thm-1} follows. 

We now prove statement \ref{sscij1thm-2}. 
Let $k\in \NN .$ 
By the homology sequence of the pair $(|N_{k+1}|, |\bigcup _{j=1}^{m}N_{k+1,j}|)$, 
we have the following exact sequence:
\begin{equation}
\cdots \rightarrow H_{2}(S_{k+1};T)\rightarrow H_{1}(\bigcup _{j=1}N_{k+1,j};T)\rightarrow H_{1}(N_{k+1};T)\rightarrow \cdots .
\end{equation}
By Lemma~\ref{sscij1lem2}, we have $H_{2}(S_{k+1};T)=0.$ Moreover, 
$H_{1}(\bigcup _{j=1}N_{k+1,j};T)\cong \bigoplus _{j=1}^{m}H_{1}(N_{k};T).$ 
Therefore, it follows that $ma_{1,k}\leq a_{1,k+1}.$ 
Thus, we have proved statement \ref{sscij1thm-2}. 

 We now prove statement \ref{sscij1thm-3}. 
 Suppose $|N_{1}|$ is connected and 
 $\check{H}^{1}({\frak L};R )\neq 0.$ Then, 
 there exists a $k\in \NN $ such that 
 $a_{1,k}\neq 0.$ From statement \ref{sscij1thm-2}, 
 it follows that $\lim _{k\rightarrow \infty } a_{1,k}=\infty . $ 
By Lemma~\ref{icfundlem}-\ref{icfundlem3}, we obtain that 
$a_{1,\infty }=\infty . $ 
Therefore, we have proved statement \ref{sscij1thm-3}. 

 Thus, we have proved Theorem~\ref{sscij1thm}.  
\qed  

\ 

 We now prove Proposition~\ref{checkprop}.\\ 
\noindent {\bf Proof of Proposition~\ref{checkprop}:} 
For each $i=1,2$, let $U_{i}$ be an open neighborhood of $h_{i}(L).$ 
Then, by the Mayer-Vietoris sequence, we have the following exact sequence:
\begin{equation}
\cdots \rightarrow H^{n+1}(U_{1}\cup U_{2};R)\rightarrow H^{n+1}(U_{1};R)\oplus H^{n+1}(U_{2};R)\rightarrow 
H^{n+1}(U_{1}\cap U_{2};R)\rightarrow \cdots .
\end{equation}
We take the direct limit $\varinjlim _{U_{1},U_{2}}$ of this sequence, where 
$U_{i}$ runs over all open neighborhoods of $h_{i}(L).$ 
Then, by \cite[p 341, Corollary 9 and p 334, Corollary 8]{Sp}, 
we obtain the following exact sequence: 
\begin{equation}
\label{chekcpropeq1}
\check{H}^{n+1}(h_{1}(L)\cup h_{2}(L);R)\rightarrow \check{H}^{n+1}(h_{1}(L);R)\oplus \check{H}^{n+1}(h_{2}(L);R)
\rightarrow \check{H}^{n+1}(C_{1,2};R)\rightarrow \cdots .
\end{equation}
By the assumption, we have $\check{H}^{n+1}(C_{1,2};R)=0.$ 
Similarly, we obtain the following exact sequence: 
\begin{multline}
\label{checkpropeq2}
\check{H}^{n+1}(\bigcup _{j=1}^{3}h_{j}(L);R)\rightarrow \check{H}^{n+1}(\bigcup _{j=1}^{2}h_{j}(L);R)\oplus 
\check{H}^{n+1}(h_{3}(L);R)\\ 
\rightarrow \check{H}^{n+1}((\bigcup _{j=1}^{2}h_{j}(L))\cap h_{3}(L);R)\rightarrow \cdots.
\end{multline}
By the assumption, we have 
$\check{H}^{n+1}((\bigcup _{j=1}^{2}h_{j}(L))\cap h_{3}(L);R)=\check{H}^{n+1}(C_{1,3}\cup C_{2,3};R)=0.$  
From these arguments, 
it follows that there exists an exact sequence: 
\begin{equation}
\check{H}^{n+1}(\bigcup _{j=1}^{3}h_{j}(L);R)\rightarrow \bigoplus _{j=1}^{3}\check{H}^{n+1}(h_{j}(L);R)\rightarrow 0.
\end{equation}
Continuing this method, we obtain the following exact sequence: 
\begin{equation}
\check{H}^{n+1}(\bigcup _{j=1}^{m}h_{j}(L);R)\rightarrow \bigoplus _{j=1}^{m}\check{H}^{n+1}(h_{j}(L);R)\rightarrow 0. 
\end{equation}
Since for each $j=1,2$, $h_{j}:L\rightarrow h_{j}(L)$ is a homeomorphism, we obtain the following exact sequence: 
\begin{equation}
\check{H}^{n+1}(L;R)\rightarrow \check{H}^{n+1}(L;R)\oplus \check{H}^{n+1}(L;R)\oplus \bigoplus _{j=3}^{m}
\check{H}^{n+1}(h_{j}(L);R)\rightarrow 0.
\end{equation} 
From this exact sequence, it follows that either $\check{H}^{n+1}(L;R)=0$ or 
$\dim _{R}\check{H}^{n+1}(L;R)=\infty .$ Thus, we have proved Proposition~\ref{checkprop}. 
\qed

\end{document}